\documentclass[a4paper,11pt]{amsart}
\usepackage{amssymb}
\usepackage[T1]{fontenc}
\usepackage[ansinew]{inputenc}
\usepackage[english,
francais]{babel}%%Francais par defaut
\usepackage{color}
\usepackage{graphicx}
\usepackage[active]{srcltx} % SRC Specials: DVI [Inverse] Search
\usepackage{amssymb,amsthm,amsfonts,amssymb,euscript,mathrsfs
}
\usepackage{amsxtra}
\usepackage{amscd}
\usepackage{float}

\ifx \undefined \cftil \def \cftil#1{\~#1}\fi

\renewcommand{\and}{et}
\renewcommand{\andify}{%
 \nxandlist{\unskip, }{\unskip{} et~}{\unskip{ } et~}}

\renewcommand{\a}{\mathfrak{a}}

\newcommand{\g}{\mathfrak{g}}
\newcommand{\gl}{\mathfrak{gl}}
\newcommand{\h}{\mathfrak{h}}

\newcommand{\n}{{\mathfrak n}}

\newcommand{\p}{{\mathfrak p}}
\newcommand{\s}{{\mathfrak s}}

\newcommand{\so}{{\mathfrak s\mathfrak o}}

\newcommand{\z}{{\mathfrak z}}

\newcommand{\R}{{\mathbb R}}

\newcommand{\K}{{\mathbb K}}
 \DeclareMathAlphabet{\got}{U}{euf}{m}{n}     %% les gothiques
\DeclareMathAlphabet{\mat}{U}{msb}{m}{n}     %% les lettres Bbb=20
\DeclareMathAlphabet{\mathbold}{OML}{cmm}{bx}{it} 

\newtheorem{theo}{Th\'eor\`eme}[subsection]
\newtheorem{pr}{Proposition}[subsection]
\newtheorem{lem}{Lemme}[subsection]
\newtheorem{defi}{D\'efinition}[subsection]
\newtheorem{defis}{D\'efinitions}[subsection]
\newtheorem{co}{Corollaire}[subsection]
\newtheorem{rem}{Remarque}[subsection]
\newtheorem{definition}{Definition}[section]

\newcommand{\cur}{\mathscr}

\newenvironment{theo*}{\vskip 1em\bf Th\'eor\`eme.\it}{\par\rm}
\newenvironment{dem}{\vskip 1em{\it D\'emonstration} :}%
{\unskip\hfill\null\nobreak\hfill\carre\vskip1em\par}
\newcommand{\carre}{\rule{1ex}{1ex}} 

\makeatletter
\providecommand*{\diff}%
{\@ifnextchar^{\DIfF}{\DIfF^{}}}
\def\DIfF^#1{%
\mathop{\mathrm{\mathstrut d}}%
\nolimits^{#1}\gobblespace}
\def\gobblespace{%
\futurelet\diffarg\opspace}
\def\opspace{%
\let\DiffSpace\!%
\ifx\diffarg(%
\let\DiffSpace\relax
\else
\ifx\diffarg[%
\let\DiffSpace\relax
\else
\ifx\diffarg\{%
\let\DiffSpace\relax
\fi\fi\fi\DiffSpace}

\providecommand*{\eu}%
{\ensuremath{\mathrm{e}}}

\providecommand*{\iu}%
{\ensuremath{\mathrm{i}}}

\newcommand{\Ad}{\mathop{\mathrm{Ad}}}

\newcommand{\ad}{\mathop{\mathrm{ad}}}

\newcommand{\Pfaff}{\mathop{\mathrm{Pfaff}}}
\newcommand{\ind}{\mathop{\mathrm{ind}}}
\newcommand{\rang}{\mathop{\mathrm{rang}}}

\newcommand{\Tr}{\mathop{\mathrm{Tr}}}

\begin{document}
%\title[Algèbres de Lie quasi-réductives]{\footnote{\footnotesize
%K\lowercase{eywords :} }}

%\selectlanguage{french}
\title{Algèbres de Lie quasi-réductives}

\date{\today}
\selectlanguage{french}
\author{Michel DUFLO}
\address{Université Denis Diderot-Paris~7,  Institut de Math\'ematiques de
Jussieu\\
 UFR de Mathématiques, Case 7012\\    F-75205 Paris~cedex~13}
\email{duflo@math.jussieu.fr}

\author{Mohamed Salah Khalgui}
\address{Faculté des Sciences de Tunis, Département de Mathématiques\\
Campus universitaire\\
1060 Tunis, Tunisie}
\email{mohamedsalah.khalgui@fst.rnu.tn}

\author{Pierre Torasso}
\address{Université de Poitiers, Laboratoire de Mathématiques et
  Applications, UMR  CNRS 6086\\
SP2MI, BP30179\\
F-86962 Chasseneuil cedex}
\email{pierre.torasso@math.univ-poitiers.fr}

\selectlanguage{english}
\begin{abstract}
We study the class of algebraic Lie algebras for which the generic
stabilizer of the coadjoint action is reductive modulo the center.
\end{abstract}

\selectlanguage{french}
\begin{abstractf}
Nous étudions la classe des algèbres de Lie dont le stabilisateur
générique de la représentation coadjointe est réductif modulo le centre.
\end{abstractf}

%\selectlanguage{french}
\maketitle

%\selectlanguage{english}
\section{Introduction}\label{0}
We fix a field $\K$ of characteristic $0$. Let $\g$ be an algebraic
Lie algebra (that is the Lie algebra of some affine algebraic group
$\mathbold{G}$ defined over $\K$). Let $\mathfrak{z}_\g$ be the center
of $\g$, and $\g^*$ the dual space, endowed with the coadjoint
representations of $\g$ and $\mathbold{G}$.  For a linear form $g\in
\g^*$, we denote by $\g(g)\subset \g$ its centralizer. We identify $
\g(g)/\mathfrak{z}_{\mathfrak{g}}$ with its image in $\gl(\g)$. It is
the Lie algebra of the group $\mathbold{G}(g)/\mathbold{Z}_\g$, where
$\mathbold{G}(g)$ is the centralizer of $g$ in $\mathbold{G}$, and
$\mathbold{Z}_\g$ the centralizer of $\g$ in $\mathbold{G}$.

\begin{definition}\label{def:typered}
A form    $g\in \g^*$ is said  to be \emph{of reductive type} if $
\g(g)/\mathfrak{z}_\g$ is a reductive Lie subalgebra of $\gl(\g)$ whose
center consists of semi-simple elements of $\gl(\g)$.
\end{definition}

Equivalently, it means that the algebraic group
$\mathbold{G}(g)/\mathbold{Z}_\g$ is reductive. The notion of form of
reductive type may be also useful for non algebraic Lie algebras (see
section \ref{sec:suppl} below), and for Lie algebras over a field of
characteristic $p\geq 2$ (see \cite{premet-skryabin-1999}), but this
is outside the scope of this paper.

%
%Rappelons, pour éviter toute ambiguïté, ce que nous entendons par
%algèbre de Lie  réductive. Par définition, un algèbre de Lie
%algébrique $\h$ est l'algèbre de Lie d'un groupe algébrique affine
%$H$ défini sur $\K$. Un élément $X \in \h$ est dit semi-simple
%(resp. nilpotent) si son image $\rho(X) \in \End(V)$ est
%semi-simple (resp. nilpotente)  pour toute représentation
%rationnelle $\rho$  de $H$ dans un espace vectoriel $V$.%  On sait
%%que chaque élément $X$ de $\g$ s'écrit de manière unique sous la
%%forme $X=X_s+X_n$ d'un élément semi-simple $X_s \in \g$ et d'un
%%élément nilpotent $X_n\in \g$ qui commutent.
%
%On note  ${}^u\h$ le radical unipotent de $\h$, c'est-à-dire le
%plus grand ideal de $\h$ formé d'éléments nilpotents. On dit que
%$\h$ est réductive si ${}^u\h$ est nul. Il est équivalent de
%demander que $\h$ soit une somme directe  $\h=\s\oplus \a$, d'un
%idéal semi-simple $\s$, et d'un idéal abélien $\a$ formé
%d'éléments semi-simples.
%
%
%
%\medskip
%Donc un élément $g\in \g^*$ est de type réductif  si et seulement
%si  le radical unipotent ${}^u(\g(g))$  de $\g(g)$ est contenu
%dans $\mathfrak{z}$.

\bigskip

Since Kirillov \cite{kirillov-1962}, it is known that the coadjoint
representation of $\mathbold{G}$ in $\g^*$ plays an important role in
representation theory. It is an experimental fact that representations
associated to coadjoint orbits of reductive type may be specially
interesting. Let us give just one example, which in fact was our
motivating example : when $\K=\R$, all unitary irreducible square
integrable representations of $\mathbold{G}_\R$ (where
$\mathbold{G}_\K$ is the group of $\K$-points of $\mathbold{G}$) are
associated to coadjoint orbits of reductive type (see
e.g. \cite[II.16]{duflo-1982} and \cite{anh-1974}).

\bigskip

The purpose of this paper is to study the classification of
coadjoint orbits of reductive type of $\mathbold{G}_\K$ in $\g^*$. The
classification of all coadjoint orbits is usually a ``wild''
problem (see   \cite{drozd-ip-1992}), but the case of coadjoint orbits
of reductive type  is much better behaved.  The set of coadjoint
orbits of reductive type may be empty. So we introduce the
following definition.

\bigskip

\begin{definition}\label{def:quasired}
We say that the  Lie algebra $\g$ is \emph{quasi-reductive} if
there exists a form $g\in \g^*$ of reductive type.
\end{definition}

It is easy to see that if $\g$ is quasi-reductive,  then the set
of linear forms of reductive type contains a non empty open set of
$\g^*$.

\bigskip
For instance, if $\g$ is reductive (that is, if $G$ is reductive),
$\g$ is quasi-reductive because the linear form $0\in \g^*$ is of
reductive type.  Assume moreover that $\K$ is algebraically closed and
let $\h\subset \g$ a Cartan subalgebra.  Then it is well known that
coadjoint orbits of reductive type are parameterized by the orbits in
$\h^*$ of the normalizer $N_{\mathbold{G}_\K}(\h)$ of $\h$ in
$\mathbold{G}_\K$. \emph{In this paper, we present similar results for
  all quasi-reductive Lie algebras}.

\medskip
In the rest of this introduction, we assume $\K$ algebraically
closed.

\bigskip
When $\g$ is an unimodular quasi-reductive Lie algebra, it has a
simple structure (see section \ref{2} below) which reduces the
classification of $\mathbold{G}_\K$-coadjoint orbits of reductive type
to the case of a Levi factor of $\mathbold{G}$. This was known for
long (see \cite{anh-1974}, \cite{duflo-1982}). The main purpose of
this paper is to extend this theory to non unimodular algebraic
quasi-reductive Lie algebras : we reduce the classification of
coadjoint orbits of reductive type to that of particular reductive
subgroups of $\mathbold{G}_\K$, whose conjugacy class depends only on
the restriction of the linear forms to the unipotent part of
$\mathfrak{z}_\g$ (theorems \ref{theo1f1} and \ref{theo1f2}). These
subgroups are the so called \emph{canonical reductive subgroups} of
$\mathbold{G}_{\K}$ and their Lie algebras the so called
\emph{canonical Lie subalgebras} of $\g$ (see Definitions
\ref{defi1f1}); it turns out that the canonical Lie subalgebras are
all conjugated under $\mathbold{G}_{\K}$. The proofs use tools on the
behavior of coadjoint orbits under induction which where introduced in
\cite{duflo-1982}, in particular the notion of linear form on $\g$ of
unipotent type.

\bigskip
``Most'' algebraic Lie algebras are not quasi-reductive. For
instance, a solvable Lie algebra $\g$ for which the center $\mathfrak{z}_\g$
is zero is quasi-reductive if and only $\mathbold{G}_\K$ has an open orbit in
the dual $\n^*$ of its unipotent radical $\n$. However, many
interesting classes of Lie algebras are quasi-reductive :
reductive Lie algebras and many of their parabolic subalgebras
(see below), Frobenius Lie algebras, etc...

\medskip
As an illustration of the notion of quasi-reductive Lie algebra, we
consider the specially interesting case of parabolic subalgebras of a
simple Lie algebra $\s$.  A striking fundamental old result is the
fact that Borel subalgebras of $\s$ are quasi-reductive (Kostant
\cite{kostant} --unpublished, see Joseph \cite{joseph-1977}). For $\s$
of type $A$ and $C$, all parabolic subalgebras are quasi-reductive
(Panyushev \cite{panyushev-2005}).  However, this is not longer true
for $\s=\so(n,\K)$ for $n\geq 7$ (see Tauvel-Yu
\cite{tauvel-yu-2004-b} ) or for exceptional $\s$ (see Baur-Moreau
\cite{baur-moreau-2009}). In theorems \ref{theo3q1}, \ref{theo3q2} and
\ref{theo3q3} \emph{we give the list of quasi-reductive parabolic
  subalgebras $\so(n,\K)$}.  Many results about parabolics of
$\so(n,\K)$ are explicitly or implicitly given in the literature,
in particular, Dvorsky \cite{dvorsky-1995}, \cite{dvorsky-1996}
and Panyushev \cite{panyushev-2005}. However, we think it is
useful to have a clear full picture. In
\cite{moreau-yakimova-2011}, Moreau and Yakimova describe the
conjugacy class of canonical reductive Lie subalgebras for
quasi-reductive parabolic subalgebras of simple Lie algebras,
and bi-parabolic subalgebras of $\gl(n,\K)$.

\bigskip
Tauvel and Yu (see \cite{tauvel-yu-2005}, ch. 40) studied a related
class of algebraic Lie algebras: let us say that $\g$ is
``stable'' if there is a regular form $g \in \g^*$ such that
$\g(g) \subset \g$ has a supplementary $\g(g)$-invariant subspace.
This is equivalent to the fact  that the generic stabilizers
$G(g)$ of  the coadjoint representation are conjugate by $G$.
Quasi-reductive Lie algebras are  stable. However, there are stable
Lie algebras which are not quasi-reductive (for instance, most
unipotent Lie algebras are not quasi-reductive, but the first
example of a non stable unipotent Lie algebra is given in
\cite{kosmann-stern-1974}). For a parabolic subalgebra of a simple Lie
algebra, we do not know if the notions are equivalent.

\selectlanguage{french}
\section{Introduction}\label{01}

Dans la suite $\mathbb{K}$ désigne un corps de caractéristique
nulle. Soit $\g$ une algèbre de Lie algébrique, 
c'est-à-dire l'algèbre de Lie d'un groupe algébrique affine
$\mathbold{G}$ défini sur $\mathbb{K}$. On désigne par
$\mathfrak{z}_{\g}$ le centre de $\g$ et par $\g^{*}$ l'espace
dual de $\g$ muni des actions coadjointes de $\g$ et
$\mathbold{G}$. \'Etant donnée une forme linéaire $g\in\g^{*}$, on
note $\g(g)\subset\g$ son centralisateur. On identifie
$\g(g)/\mathfrak{z}_{\g}$ avec son image dans $\mathfrak{gl}(\g)$.
C'est l'algèbre de Lie du groupe
$\mathbold{G}(g)/\mathbold{Z}_{\g}$, où $\mathbold{G}(g)$ est le
centralisateur de $g$ dans $\mathbold{G}$ et $\mathbold{Z}_{\g}$
est le centralisateur de $\g$ dans $\mathbold{G}$.
\begin{defi}
Une forme linéaire $g\in\g^{*}$ est dite \emph{de type réductif} si
$\g(g)/\mathfrak{z}_{\g}$ est une sous-algèbre de Lie réductive de
$\mathfrak{gl}(\g)$ dont le centre est formé d'éléments semi-simples de
$\mathfrak{gl}(\g)$.
\end{defi}
De manière équivalente, cela revient à demander que le groupe
algébrique $\mathbold{G}(g)/\mathbold{Z}_{\g}$ soit réductif. Il peut
également être utile d'introduire la notion de forme linéaire de type
réductif pour des algèbres de Lie non algébriques (voir la section
\ref{sec:suppl} ci-après), et pour les algèbres de Lie sur un corps de
caractéristique $p\geq2$ (voir \cite{premet-skryabin-1999}), mais ceci
reste hors du champ du présent article.

\bigskip
Depuis les travaux de Kirillov \cite{kirillov-1962}, il est bien connu
que la représentation coadjointe de $\mathbold{G}$ dans $\g^{*}$ joue
un rôle important en théorie des représentations. C'est un fait
expérimental que les représentations associées aux orbites de type
réductif sont particulièrement intéressantes. Par exemple, lorsque
$\mathbb{K}=\mathbb{R}$, chaque série discrète de
$\mathbold{G}_{\mathbb{R}}$ (où $\mathbold{G}_{\mathbb{K}}$ désigne le
groupe des $\mathbb{K}$-points de $\mathbold{G}$) est associée à une
orbite de type réductif (voir \cite[II.16]{duflo-1982} et
\cite{anh-1974}) et ceci a été notre principale motivation pour
introduire et étudier les formes linéaires de type réductif.

\bigskip
Le but de cet article est d'étudier la classification des orbites de
type réductif de $\mathbold{G}_{\mathbb{K}}$ dans $\g^{*}$. Si la
classification de toutes les orbites coadjointes est en général un
problème \flqq sauvage\frqq~(voir \cite{drozd-ip-1992}), il n'en est
pas de même pour les orbites coadjointes de type réductif. Il se peut
que l'ensemble des orbites coadjointes de type réductif soit
vide. C'est pourquoi nous introduisons la définition suivante.
\begin{defi}
Nous disons que l'algèbre de Lie $\g$ est \emph{quasi-réductive} s'il
existe une forme de type réductif $g\in\g^{*}$.
\end{defi}
Il est facile de voir que si $\g$ est quasi-réductive, l'ensemble des
formes linéaires de type réductif contient un ouvert non vide de
$\g^{*}$.

\bigskip
Par exemple, si $\g$ est réductive (c'est à dire, $\mathbold{G}$ est
réductif), $\g$ est quasi-réductive puisque la forme linéaire
$0\in\g^{*}$ est de type réductif. Supposons de plus $\mathbb{K}$
algébriquement clos et soit $\h\subset\g$ une sous-algèbre de Cartan.
Il est alors bien connu que les orbites coadjointes de type
réductif sont paramétrées par les orbites dans $\h^{*}$ du
normalisateur $N_{\mathbold{G}_{\mathbb{K}}}(\h)$ de $\h$ dans
$\mathbold{G}_{\mathbb{K}}$. \emph{ Dans cet article, nous présentons des
  résultats similaires pour toutes les algèbres de Lie quasi-réductives.}

\bigskip
Dans la suite de cette introduction, on suppose que $\mathbb{K}$ est
algébriquement clos.

\bigskip
Lorsque $\g$ est une algèbre de Lie quasi-réductive unimodulaire, elle
a une structure particulièrement simple (voir la section \ref{2}
ci-après) permettant de ramener la classification des
$\mathbold{G}_{\mathbb{K}}$-orbites coadjointes de type réductif au
même problème pour un de ses facteurs de Levi. Ceci était connu depuis
longtemps (voir \cite{anh-1974} et \cite{duflo-1982}). Dans cet
article, notre but principal est d'étendre cette théorie aux algèbres
de Lie quasi-réductives non unimodulaires~: nous ramenons la
classification des orbites coadjointes de type réductif au même
problème pour certains sous-groupes réductifs de
$\mathbold{G}_{\mathbb{K}}$ dont la classe de conjugaison ne dépend
que de la restriction des formes linéaires au radical unipotent de
$\mathfrak{z}_{\g}$ (voir les théorèmes \ref{theo1f1} et
\ref{theo1f2}). Ces sous-groupes sont appelés les \emph{sous-groupes
  réductifs canoniques} de $\mathbold{G}_{\K}$ et leurs algèbres de
Lie, les \emph{sous-algèbres réductives canoniques} de $\g$ (voir
Définitions \ref{defi1f1}) ; il s'avère que ces dernières sont toutes
conjuguées entre elles sous l'action de $\mathbold{G}_{\K}$.  Les
démonstrations utilisent des outils adaptés au comportement des
orbites coadjointes par rapport à l'induction, qui ont été introduits
dans \cite{duflo-1982}, en particulier la notion de forme linéaire de
type unipotent.

\bigskip
\flqq Beaucoup\frqq~ d'algèbres de Lie algébriques ne sont pas
quasi-réductives. Par exemple, une algèbre de Lie résoluble $\g$ dont
le centre $\mathfrak{z}_{\g}$ est réduit à $0$ est quasi-réductive si
et seulement si $\mathbold{G}_{\mathbb{K}}$ a une orbite ouverte dans
le dual $\n^{*}$ de son radical unipotent $\n$. Cependant, de
nombreuses classes intéressantes d'algèbres de Lie sont
quasi-réductives~: les algèbres de Lie réductives et beaucoup de leurs
sous-algèbres paraboliques (voir ci-après), les algèbres de Lie de
Frobenius, etc...

\bigskip
Pour illustrer la notion d'algèbre de Lie quasi-réductive, nous
considérons le cas particulièrement intéressant des sous-algèbres
paraboliques d'une algèbre de Lie simple $\mathfrak{s}$. Un résultat
fondamental, connu de longue date et particulièrement frappant est le
fait que les sous-algèbres de Borel de $\mathfrak{s}$ sont
quasi-réductives (Kostant \cite{kostant}, non publié, voir Joseph
\cite{joseph-1977}). Pour $\mathfrak{s}$ de type $A$ et $C$, toutes
les sous-algèbres paraboliques sont quasi-réductives (Panyushev
\cite{panyushev-2005}). Cependant, ceci n'est plus vrai pour
$\mathfrak{s}=\mathfrak{so}(n,\mathbb{K})$ avec $n\geq7$ (voir
Tauvel-Yu \cite{tauvel-yu-2004-b}), ou pour $\mathfrak{s}$ de type
exceptionnel (voir Baur-Moreau \cite{baur-moreau-2009}). Dans l'énoncé
des théorèmes \ref{theo3q1}, \ref{theo3q2} et \ref{theo3q3} \emph{nous
  donnons la liste des sous-algèbres paraboliques de
  $\mathfrak{so}(n,\mathbb{K})$ qui sont quasi-réductives}. De
nombreux résultats concernant les sous-algèbres paraboliques de
$\mathfrak{so}(n,\mathbb{K})$ se trouvent de manière explicite ou
implicite dans la littérature, en particulier, Dvorsky
\cite{dvorsky-1995}, \cite{dvorsky-1996} et Panyushev
\cite{panyushev-2005}. Cependant, nous pensons utile d'en avoir un
tableau complet et clair. Dans \cite{moreau-yakimova-2011}, Moreau et
Yakimova décrivent la classe de conjugaison des sous-algèbres
réductives canoniques pour les sous-algèbres paraboliques des algèbres
de Lie simples, et des sous-algèbres bi-paraboliques de $\gl(n,\K)$.

\bigskip
Tauvel et Yu (voir \cite[ch. 40]{tauvel-yu-2005}) ont étudié une
classe d'algèbres de Lie reliée à celle des algèbres de Lie
quasi-réductives~: disons qu'une algèbre de Lie $\g$ est \flqq
stable\frqq~ s'il existe une forme linéaire régulière $g\in\g^{*}$
telle que $\g(g)\subset\g$ admet un sous-espace supplémentaire
$\g(g)$-invariant. Ceci est équivalent au fait que les
stabilisateurs génériques $G(g)$ de la représentation coadjointe
sont conjugués sous $G$. Les algèbres de Lie quasi-réductives sont
stables. Cependant, il existe des algèbres de Lie stables qui ne
sont pas quasi-réductives par exemple,  beaucoup
d'algèbres de Lie unipotentes ne sont pas quasi-réductives, mais
le premier exemple d'algèbre de Lie unipotente non stable est
donné dans \cite{kosmann-stern-1974}). Pour ce qui concerne les
sous-algèbres paraboliques d'une algèbre de Lie simple, nous ne
savons pas si ces notions sont équivalentes.

\section{Algèbres de Lie quasi-réductives% et formes de type réductif :
%définition et premières propriétés
.}\label{1}

Dans la suite et sauf mention du contraire, $\mathbb{K}$ est un corps
algébriquement clos de caractéristique nulle. Par algèbre de Lie
algébrique, nous entendons une algèbre de Lie qui soit l'algèbre de
Lie d'un groupe algébrique affine défini sur $\mathbb{K}$.

Si $\mathbold{G}$ est un groupe algébrique affine, on l'identifie à
l'ensemble de ses points rationnels, on note ${}^u\!\mathbold{G}$ le
radical unipotent de $\mathbold{G}$, $\g$ son algèbre de Lie,
${}^u\!\g$ l'algèbre de Lie de ${}^u\!\mathbold{G}$. Rappelons que
$\mathbold{G}$ admet une décomposition de Levi : il existe un
sous-groupe réductif $\mathbold{R}\subset\mathbold{G}$, appelé
\emph{sous-groupe de Levi ou facteur réductif de $\mathbold{G}$}, dont
la classe de conjugaison modulo ${}^u\!\mathbold{G}$ est uniquement
déterminée, tel que
$\mathbold{G}=\mathbold{R}\,{}^u\!\mathbold{G}$. Au niveau des
algèbres de Lie, on a la décomposition de Levi $\g=\mathfrak{r}\oplus
{}^u\!\g$. On dit que $\mathfrak{r}$ est un facteur réductif de $\g$
et que ${}^u\!\g$ est son radical unipotent.

Pour une action rationnelle de $\mathbold{G}$ sur une variété
algébrique $\mathbold{X}$ et $x \in \mathbold{X}$, on note
$\mathbold{G}(x)$ son stabilisateur et $\g(x)$ l'algèbre de Lie de
$\mathbold{G}(x)$.

Si $V$ est un espace vectoriel, on note $V^*$ l'espace vectoriel
dual. Si $W\subset V$ est un sous-espace vectoriel, et $\mu \in
W^*$, on note $V^*_\mu \subset V^*$ l'espace affine formé des $g
\in V^*$ dont la restriction $g|_W$ à $W$ est égale à $\mu$.

Si $\g$ est une algèbre de Lie, on note $\mathfrak{z}_{\g}$ ou plus simplement
$\mathfrak{z}$ son centre. Si $\g$ est une algèbre algébrique et commutative, on
note $^{r}\!\g$ son unique facteur réductif.

Si $\mathbold{G}$ est un groupe algébrique, on note $\mathbold{G}_{0}$
sa composante neutre et on désigne par $\mathbold{Z}_{\g}$ ou plus
simplement $\mathbold{Z}$ le centralisateur dans $\mathbold{G}$ de
l'algèbre de Lie $\g$. L'algèbre de Lie de $\mathbold{Z}$ est le
centre $\mathfrak{z}$ de l'algèbre de Lie $\g$ de $\mathbold{G}$.
%\begin{defi}

%\end{defi}

\subsection{Formes de type réductif.}\label{1a}

\begin{defis}
Soit $\mathfrak{g}$ une algèbre de Lie algébrique.

(i) une forme linéaire $g\in\mathfrak{g}^{*}$ est dite \emph{de type
  réductif} si le radical unipotent de $\mathfrak{g}(g)$ est central
dans $\mathfrak{g}$,

(ii) l'algèbre de Lie $\mathfrak{g}$ est dite \emph{quasi-réductive}
s'il existe une forme linéaire sur $\mathfrak{g}$ qui soit de type
réductif.
\end{defis}

\begin{rem}\label{rem1a1}
Soit $\mathfrak{g}$ une algèbre de Lie algébrique et $\mathfrak{z}$
son centre. Il est clair qu'une forme linéaire $g\in\mathfrak{g}^{*}$
est de type réductif si et seulement si
$^{u}\!(\mathfrak{g}(g))=\,^{u}\!\mathfrak{z}$ et que, dans ce cas,
$\mathfrak{g}(g)$ possède un unique facteur réductif, que l'on note
$\mathfrak{r}_{g}$.

Soit $\mathbold{G}$ un groupe algébrique d'algèbre de Lie
$\mathfrak{g}$. Si une forme linéaire $g$ est de type réductif, il en
est de même de toutes celles qui sont contenues dans son orbite sous
l'action coadjointe de $\mathbold{G}$~: l'orbite sous $\mathbold{G}$
d'une telle forme est dite de type réductif.

On voit alors qu'une forme linéaire $g\in\g^{*}$ est de type réductif
si et seulement si l'une des propriétés équivalentes suivantes est
vérifiée~:
\begin{gather}
\mathbold{G}(g)/\mathbold{Z} \mbox{ est réductif}, \label{eq:red1}\\
\g(g)/\z \mbox{ est réductive}, \label{eq:red2}\\
{}^u\!(\mathbold{G}(g))\subset\mathbold{Z} , \label{eq:red3}\\
{}^u\!(\g(g))\subset \z . \label{eq:red4}
\end{gather}
\end{rem}

\subsection{Formes linéaires  de type unipotent.}\label{1b}

Cette notion est introduite dans \cite[I.10]{duflo-1982}  pour
encoder le résultat final de la théorie de Mackey.

Rappelons quelques définitions. Soit $\g$ une algèbre de Lie et $g\in
\g^*$. On note $\beta_g$ la forme bilinéaire $(X,Y) \mapsto g([X,Y])$ sur
$\g \times \g$. Soit $\mathfrak{b}$ un sous-espace de $\g$. On note
$\mathfrak{b}^g \subset \g$ son orthogonal dans $\g$.

\begin{defi}\label{def:coisotrope}
Soit $\mathfrak{b}$ un sous-espace de $\g$. On dit que $\mathfrak{b}$ est
\emph{coisotrope relativement à $g$} s'il contient $\mathfrak{b}^g$.
\end{defi}

\begin{defi}\label{def:typeunipotent}
Soit $\g$ une algèbre de Lie algébrique et $g\in \g^*$. On dit que $g$
est de \emph{type unipotent} si et seulement s'il vérifie les
conditions (U1) et (U2) qui suivent.

(U1) Soit $\g(g)=\mathfrak{r}_{g}+{}^u(\g(g))$ une décomposition de Levi de
$\g(g)$. Alors la restriction de $g$ à $\mathfrak{r}_{g}$ est nulle (cette
condition ne dépend pas du choix de la décomposition de Levi).

(U2) Il existe une sous-algèbre $\mathfrak{b}$ de $\g$ telle que
\begin{enumerate}
\item $\mathfrak{b}$ est algébrique,
\item  $\mathfrak{b}$ est coisotrope relativement à $g$,
\item $\mathfrak{b}=\g(g)+{}^u\!\mathfrak{b}$.
\end{enumerate}
\end{defi}

On note $\g^{*}_{u}$ l'ensemble des formes de type unipotent. Il est
clair qu'il est invariant sous l'action du groupe des automorphismes
de $\g$. L'orbite d'un élément de $g$ de $\g^{*}_{u}$ sera dite de
type unipotent.

La condition (U2) n'est peut-être pas très intuitive. Nous verrons
(théorème \ref{theo1d1})
que si $g$ est une forme de type réductif,
alors elle est de type unipotent si et seulement si la condition
(U1) est vérifiée.

\begin{rem}\label{rem1b1}
Pour une forme linéaire $g\in\g^{*}$, la condition (U1) est
équivalente à la condition

(U'1) il existe une sous-algèbre de Cartan de $\mathfrak{r}_{g}$ sur
laquelle $g$ s'annule.
\end{rem}

%\begin{dem}
%En effet, les sous-algèbres de Cartan de $\mathfrak{r}_{g}$ sont
%toutes conjuguées par le sous-groupe algébrique connexe d'algèbre de
%Lie $\mathfrak{r}_{g}$, lequel est contenu dans $\mathbold{G}(g)$, et
%leur réunion contient l'ouvert de Zariski dense des éléments
%semi-simples réguliers.
%\end{dem}

Dans \cite[I.20]{duflo-1982}, on construit pour tout $g\in\g^{*}$ une
sous-algèbre coisotrope particulière, dite canonique relativement à
$g$, notée $\mathfrak{b}_{g}$. Elle est déterminée par récurrence sur
la dimension de $\g$ par les conditions suivantes~: on note $u$ la
restriction de $g$ à $^{u}\!\g$, $\h$ le stabilisateur de $u$ dans
$\g$ et $h$ la restriction de $g$ à $\h$. Alors, on a
\begin{enumerate}
\item $\mathfrak{b}_{g}=\mathfrak{b}_{h}+\,^{u}\!\g$, si $\h\neq\g$,
\item $\mathfrak{b}_{g}=\g$, si $\h=\g$.
\end{enumerate}

D'après \cite[fin de I.20]{duflo-1982}, on a~:

\begin{lem}\label{lem1b1}
Soit $g\in\g^{*}_{u}$. Alors $\mathfrak{b}_{g}$ vérifie la condition
\emph{(U2)} relativement à $g$.
\end{lem}

Si $g\in\mathfrak{g}^{*}_{u}$, on désigne par $L(g)$ l'ensemble des
$\lambda\in\mathfrak{g}(g)^{*}$ tels que $g$ et $\lambda$ aient même
restriction à $^{u}\!(\mathfrak{g}(g))$. On note $\cur{D}$ l'ensemble
des couples $(g,\lambda)$ avec $g\in\mathfrak{g}^{*}$ de type
unipotent et $\lambda\in L(g)$. Le groupe $\mathbold{G}$ agit de
manière naturelle dans $\cur{D}$.

 D'après
\cite[I.26]{duflo-1982}, on a~:
\begin{pr}\label{pr1b1}
(i) Soit $(g,\lambda)\in\cur{D}$. Alors, la $\mathbold{G}$-orbite
  d'une forme $g'\in\mathfrak{g}^{*}$ telle que
  $g'_{\mid\mathfrak{g}(g)}=\lambda$ et
  $g'_{\mid^{u}\!\mathfrak{b}_{g}}=g_{\mid^{u}\!\mathfrak{b}_{g}}$ ne
  dépend que de $g$ et de $\lambda$~; on la note
  $\cur{O}_{g,\lambda}$.

(ii) L'application $(g,\lambda)\mapsto\cur{O}_{g,\lambda}$
induit une bijection de $\mathbold{G}\backslash\cur{D}$ sur
$\mathbold{G}\backslash\mathfrak{g}^{*}$.
\end{pr}

\begin{defi}\label{defi1b1}
Deux formes linéaires $g_{j}$, $j=1,2$, sur $\g$ sont dites
\emph{$u$-équivalentes} s'il existe $g\in\g^{*}_{u}$ et $\lambda_{j}\in
L(g)$, $j=1,2$, tels que $g_{j}\in\cur{O}_{g,\lambda_{j}}$, $j=1,2$.
\end{defi}

\begin{rem}\label{rem1b2}
Le fait pour deux formes linéaires sur $\g$ d'être $u$-équivalentes,
définit une relation d'équivalence sur $\g^{*}$ dont l'ensemble des
classes d'équivalences est naturellement en bijection avec l'ensemble
des $\mathbold{G}$-orbites de type unipotent.
\end{rem}

%La proposition suivante rassemble des résultats contenus dans
%\cite[I.28]{duflo-1982}~:
%\begin{pr}\label{pr1b2} Soit $g\in\g^{*}$.

%(i) Il existe $g_{u}\in\g^{*}_{u}$ telle que $g$ et $g_{u}$ aient même
%  restriction à $^{u}\!\mathfrak{b}_{g}$.

%On choisit un facteur réductif $\mathfrak{r}$ de $\g(g_{u})$ et on note
%$b$ (resp. $\lambda$, $\mu$) la restriction de $g$ à
%$\mathfrak{b}_{g}$ (resp. $\g(g_{u})$, $\mathfrak{r}$). Alors, on a

%(ii) $\mathfrak{b}_{g}=\mathfrak{b}_{g_{u}}$, si bien que
%$g\in\cur{O}_{g_{u},\lambda}$,

%(iii) $\g(g)+\,^{u}\mathfrak{b}_{g}(b)=
%\g(g)(\lambda)+\,^{u}\mathfrak{b}_{g}(b)$,

%(iv)
%$\dim\g/\g(g)=\dim\g/\g(g_{u})+\dim\mathfrak{r}/\mathfrak{r}(\mu)$.
%\end{pr}

\subsection{Formes fortement régulières.}\label{1c}

Dans la suite $\g$ désigne une algèbre de Lie algébrique et
$\mathbold{G}$ un groupe algébrique d'algèbre de Lie $\g$. On note
$2e_{\g}$ ou plus simplement $2e$ la dimension des orbites coadjointes
régulières.

Rappelons qu'une forme linéaire $g\in\got{g}^{*}$ est dite fortement
régulière si elle est régulière, auquel cas $\got{g}(g)$ est une
algèbre de Lie algébrique commutative, et si de plus le tore
$\got{j}_{g}$, unique facteur réductif de $\got{g}(g)$, est de
dimension maximale lorsque $g$ parcourt l'ensemble des formes
régulières. Il est bien connu que l'ensemble des formes fortement
régulières, noté $\got{g}^{*}_{r}$, est un ouvert de Zariski
$\mathbold{G}$-invariant non vide de $\got{g}^{*}$.

Les tores $\got{j}_{g}$, $g\in\got{g}^{*}_{r}$ sont appelés les
sous-algèbres de Cartan-Duflo de $\got{g}$ dans
\cite[1.18]{khalgui-torasso-1993}. Ils sont deux à deux conjugués sous
l'action du groupe adjoint connexe.

Soit $\got{j}$ une sous-algèbre de Cartan-Duflo de $\g$. On désigne
par $\mathbold{H}=Z_{\mathbold{G}}(\got{j})$ et
$\mathbold{H}'=N_{\mathbold{G}}(\got{j})$ son centralisateur et son
normalisateur dans $\mathbold{G}$. Ce sont deux sous-groupes presque
algébriques de $\mathbold{G}$ d'algèbre de Lie $\h=\got{g}^{\got{j}}$,
le centralisateur de $\mathfrak{j}$ dans $\g$. Par suite le groupe quotient
$W(\mathbold{G},\got{j})=\mathbold{H}'/\mathbold{H}$ est un groupe
fini.

Le dual de l'algèbre de Lie $\h$ s'identifie à l'orthogonal de
$[\mathfrak{j},\mathfrak{g}]$ dans $\g^{*}$ qui n'est autre que
$\g^{*\mathfrak{j}}$, l'ensemble des points fixes de $\mathfrak{j}$
dans $\g^{*}$.
%Si $g\in\h^{*}$, on définit la forme bilinéaire alternée $\beta_{g}$
%sur $[\mathfrak{j},\mathfrak{g}]$ en posant $\beta_{g}(X,Y)=\langle
%g,[X,Y]\rang$, $X,Y\in[\mathfrak{j},\mathfrak{g}]$.
On montre alors qu'une forme linéaire $g\in\h^{*}$ est fortement
régulière relativement à $\g$ si et seulement si~:

\begin{enumerate}
\item elle est régulière relativement à $\h$,
\item la restriction à $[\mathfrak{j},\mathfrak{g}]$ de la forme
  bilinéaire $\beta_{g}$ est non dégénérée.
\end{enumerate}
On voit en particulier que $[\mathfrak{j},\mathfrak{g}]$ est toujours
de dimension paire, disons $2d$.

Choisissons un élément non nul
$\eta\in\wedge^{2d}([\mathfrak{j},\mathfrak{g}]^{*})$. On note alors
$\pi_{\mathfrak{j},\eta}$ (ou plus simplement $\pi_{\mathfrak{j}}$) la
fonction polynomiale sur $\h^{*}$ telle que~:
\begin{equation}\label{eq1c1}
\frac{1}{d!}\wedge^{d}((\beta_{g})_{\vert[\mathfrak{j},\mathfrak{g}]})
=\pi_{\mathfrak{j}}(g)\eta.
\end{equation}
C'est un élément non nul de $S(\h)$ qui est
$\mathbold{H}'$-semi-invariant de poids
$(\det\Ad_{\g/\h})^{-1}$.

On a alors le résultat suivant~:
\begin{theo}\label{theo1c1}
(i) L'application $\cur{O}\mapsto\cur{O}\cap\h^{*}$ induit une
  bijection de l'ensemble des orbites coadjointes fortement régulières
  de $\mathbold{G}$ sur l'ensemble des orbites coadjointes de
  $\mathbold{H}'$ contenues dans $(\g^{*}_{r})^{\mathfrak{j}}$.

(ii) Une forme linéaire $g\in\h^{*}$ est fortement
régulière relativement à $\g$ si et seulement si~:
\begin{enumerate}
\item elle est régulière relativement à $\h$,
\item $\pi_{\mathfrak{j}}(g)\neq0$.
\end{enumerate}
\end{theo}

Pour ce qui précède, voir \cite[numéro
  11]{andler-1985}.

\subsection{Algèbres de Lie quasi-réductives.}\label{1d}

Dans ce numéro, on établit les premières propriétés des algèbres de
Lie quasi-réductives.

Si $\mathfrak{g}$ est une algèbre de Lie algébrique, on désigne par
$\mathfrak{g}^{*}_{red}$ (resp.  $\mathfrak{g}^{*}_{red,u}$) le
sous-ensemble de $\mathfrak{g}^{*}$ constitué des formes linéaires de
type réductif (resp. de type réductif et unipotent).

\begin{rem}\label{rem1d1}
Une algèbre de Lie algébrique réductive $\mathfrak{g}$ est
 quasi-réductive. Elle possède une unique forme de type réductif et
 unipotent, savoir $0$ et l'on a $\mathfrak{r}_{0}=\mathfrak{g}$. De
 plus si l'on identifie $\mathfrak{g}$ et $\mathfrak{g}^{*}$ au moyen
 d'une forme bilinéaire symétrique non dégénérée invariante, les
 formes de type réductif (resp. régulières de type réductif) sont
 celles qui s'identifient aux éléments semi-simples
 (resp. semi-simples réguliers de $\mathfrak{g}$).
\end{rem}

\begin{theo}\label{theo1d1}
Soit $\mathfrak{g}$ une algèbre de Lie algébrique et $\mathfrak{z}$ son centre.

a) Soit $g\in\mathfrak{g}^{*}$ une forme de type unipotent. Alors, les
assertions suivantes sont équivalentes~:

(i) $g$ est de type réductif,

(ii) il existe une forme $g'\in\mathfrak{g}^{*}$ $u$-équivalente à $g$
telle que $\mathfrak{g}(g')=\mathfrak{j}\oplus\,^{u}\!\mathfrak{z}$,
où $\mathfrak{j}$ est un tore,

(iii) il existe une forme $g'\in\mathfrak{g}^{*}$ $u$-équivalente à $g$
qui soit de type réductif.

%Lorsque les conditions équivalentes sont satisfaites, on peut choisir
%la forme linéaire $g'$ $u$-équivalente à $g$ des conditions (ii) et
%(iii) telle que
%$g_{\vert\,^{u}\!\mathfrak{b}_{g}}=g'_{\vert\,^{u}\!\mathfrak{b}_{g}}$
%et $\mathfrak{r}_{g'}$ est une sous-algèbre réductive de
%$\mathfrak{r}_{g}$ de même rang.

b) Soit $g\in\g^{*}$ une forme de type réductif et unipotent, $g'$ une
forme $u$-équivalente à $g$ et $\lambda\in L(g)$ telle que
$g'\in\cur{O}_{g,\lambda}$. Alors, $g'$ est de type réductif si et
seulement si $\lambda_{\vert\mathfrak{r}_{g}}$ l'est. Dans ce cas,
$\mathfrak{r}_{g'}$ est $\mathbold{G}$-conjuguée à la sous-algèbre
réductive de rang maximal
$\mathfrak{r}_{g}(\lambda_{\vert\mathfrak{r}_{g}})$ de
$\mathfrak{r}_{g}$.

c) Les assertion suivantes sont équivalentes~:

(i) $\mathfrak{g}$ est quasi-réductive,

(ii) il existe $g\in\mathfrak{g}^{*}$ une forme de type unipotent et
de type réductif,

(iii) il existe $g\in\mathfrak{g}^{*}$ une forme fortement régulière et
de type réductif,

d) On suppose $\mathfrak{g}$ quasi-réductive. Alors,

(i) l'ensemble des formes
fortement régulières sur $\mathfrak{g}$ est l'ensemble des
$g\in\mathfrak{g}^{*}$ telles que
$\mathfrak{g}(g)=\mathfrak{j}\oplus\,^{u}\!\mathfrak{z}$, avec
$\mathfrak{j}$ un tore. En particulier, l'ouvert de Zariski
$\mathfrak{g}^{*}_{r}$ est contenu dans $\mathfrak{g}^{*}_{red}$,

(ii) l'ensemble des formes fortement régulières est égal à l'ensemble
des formes régulières de type réductif.

e) Si $g\in\mathfrak{g}^{*}_{red}$, tout tore maximal de
$\mathfrak{g}(g)$ est une sous-algèbre de Cartan-Duflo de
$\mathfrak{g}$.

f) Si $g\in\mathfrak{g}^{*}_{red}$, $g$ est de type unipotent si et
seulement si elle vérifie la condition \emph{(U1)}
\end{theo}

\begin{dem}
a) Montrons (i) $\Rightarrow$ (ii). Soit $g\in\mathfrak{g}^{*}$ une
forme de type unipotent et réductif, $\mathfrak{b}_{g}$ la
sous-algèbre acceptable canonique associée. D'après le lemme
\ref{lem1b1}, $\mathfrak{r}_{g}$ est également un facteur réductif de
$\mathfrak{b}_{g}$. Soit $\lambda\in\mathfrak{g}(g)^{*}$ et
$g'\in\mathfrak{g}^{*}$ telle que $g'_{\mid\mathfrak{g}(g)}=\lambda$
et
$g'_{\mid^{u}\!\mathfrak{b}_{g}}=g_{\mid^{u}\!\mathfrak{b}_{g}}$. Alors,
il résulte de \cite[I.26 et proposition I.28]{duflo-1982} et
\cite[lemme 1.21]{khalgui-torasso-1993} que la forme $g'$ est
$u$-équivalente à $g$, que $\mathfrak{b}_{g}=\mathfrak{b}_{g'}$ et
que, si $b=g_{\mid\mathfrak{b}_{g}}$, on a~: \begin{eqnarray}
  \mathfrak{g}(g')+(^{u}\!\mathfrak{b}_{g})(b) & = &
  \mathfrak{r}_{g}(\lambda_{\vert\mathfrak{r}_{g}})
  \oplus(^{u}\!\mathfrak{b}_{g})(b)\label{eq1a1}\\ \dim\mathfrak{g}/\mathfrak{g}(g')
  & = & \dim\mathfrak{g}/\mathfrak{g}(g)+
  \dim\mathfrak{r}_{g}/\mathfrak{r}_{g}(\lambda_{\vert\mathfrak{r}_{g}})\label{eq1a2}
\end{eqnarray}
la relation \ref{eq1a2} étant équivalente à~:
\begin{equation}\label{eq1a3}
\dim\mathfrak{g}(g')=\dim\mathfrak{r}_{g}(\lambda_{\vert\mathfrak{r}_{g}})
+\dim\,^{u}\!(\mathfrak{g}(g)).
\end{equation}

Supposons $\lambda$ choisie telle que
$\mathfrak{r}_{g}(\lambda_{\vert\mathfrak{r}_{g}})$ soit une
sous-algèbre de Cartan de $\mathfrak{r}_{g}$. Alors, il résulte de
\ref{eq1a1} que $\mathfrak{r}_{g}(\lambda_{\vert\mathfrak{r}_{g}})$
est conjugué par un élément de $\exp((^{u}\!\mathfrak{b}_{g})(b))$à un
facteur réductif $\mathfrak{j}$ de $\mathfrak{g}(g')$. Comme on a
$^{u}\!(\mathfrak{g}(g))=
\,^{u}\!\mathfrak{z}\subset\,^{u}\!(\mathfrak{g}(g'))$, il suit de
\ref{eq1a3} que
$\mathfrak{g}(g')=\mathfrak{j}\oplus\,^{u}\!\mathfrak{z}$,
$\mathfrak{j}$ étant un tore.

L'implication (ii) $\Rightarrow$ (iii) est évidente.

Montrons (iii) $\Rightarrow$ (i). Soit donc $g'\in\mathfrak{g}^{*}$
une forme linéaire de type réductif et $u$-équivalente à $g$. Quitte à
remplacer $g'$ par un élément de son orbite sous l'action d'un groupe
algébrique connexe d'algèbre de Lie $\mathfrak{g}$, on peut supposer
que $g$ et $g'$ ont même restriction à $^{u}\!\mathfrak{b}_{g}$,
auquel cas $\mathfrak{b}_{g}=\mathfrak{b}_{g'}$ (voir encore
\cite[I.28 et sa démonstration]{duflo-1982} ou
\cite[1.16]{khalgui-torasso-1993}). Soit
$\lambda=g'_{\mid\mathfrak{g}(g)}$ et
$b=g_{\mid\mathfrak{b}_{g}}$. Alors, les relations \ref{eq1a1} et
\ref{eq1a3} s'appliquent. Il suit de \ref{eq1a1} et du fait que
$^{u}\!(\mathfrak{g}(g'))=
\,^{u}\!\mathfrak{z}\subset\,^{u}\!(\mathfrak{b}_{g})(b)$ que
$\mathfrak{r}_{g}(\lambda_{\vert\mathfrak{r}_{g}})$ est conjugué au
facteur réductif de $\mathfrak{g}(g')$. La relation \ref{eq1a3} montre
alors que $\dim\,^{u}\!(\mathfrak{g}(g))=
\dim\,^{u}\!(\mathfrak{g}(g'))=\dim\,^{u}\!\mathfrak{z}$ et donc que
$^{u}\!\mathfrak{z}=\,^{u}\!(\mathfrak{g}(g))$.

b) Soit $b$ la restriction de $g$ à $\mathfrak{b}$. D'après
\cite[lemme 1.21]{khalgui-torasso-1993}, on peut supposer que $g$ et
$g'$ ont même restriction à $^{u}\!\mathfrak{b}_{g}$, que
$\lambda=g'_{\vert\mathfrak{r}_{g}}$ et que les relations \ref{eq1a1}
et \ref{eq1a3} sont satisfaites. Mais alors, les facteurs réductifs de
$\g(g')$ d'une part et de
$\mathfrak{r}_{g}(\lambda_{\vert\mathfrak{r}_{g}})$ d'autre part sont
conjugués par un élément de $\exp((^{u}\!\mathfrak{b}_{g})(b))$. On
voit donc que $g'$ est de type réductif si et seulement si
$\lambda_{\vert\mathfrak{r}_{g}}$ l'est.

Supposons donc que $g'$ est de type réductif. Comme
$^{u}(\!\mathfrak{g}(g'))=\,^{u}\!(\mathfrak{g}(g))=\,^{u}\!\mathfrak{z}$,
on en déduit que $\mathfrak{r}_{g}(\lambda_{\vert\mathfrak{r}_{g}})$ est
un facteur réductif de $\mathfrak{g}(g')$. Notre assertion est alors
claire.

c) Compte tenu de a), il suffit de montrer que toute forme linéaire
$g\in\mathfrak{g}^{*}$ telle que
$\mathfrak{g}(g)=\mathfrak{j}\oplus\,^{u}\!\mathfrak{z}$, avec
$\mathfrak{j}$ un tore, est fortement régulière. Soit donc $g$ une
telle forme linéaire. On a
$\mathfrak{g}=\mathfrak{g}^{\mathfrak{j}}\oplus[\mathfrak{j},\mathfrak{g}]$,
d'où résulte la décomposition duale
$\mathfrak{g}^{*}=
\mathfrak{g}^{*\mathfrak{j}}\oplus[\mathfrak{j},\mathfrak{g}]^{*}$. Comme
$\mathfrak{g}(g)\subset\mathfrak{g}^{\mathfrak{j}}$, la forme
bilinéaire alternée $(\beta_{g})_{\mid[\mathfrak{j},\mathfrak{g}]}$
est non dégénérée, de sorte que l'application $X\mapsto X.g$ induit un
isomorphisme linéaire de $[\mathfrak{j},\mathfrak{g}]$ sur
$[\mathfrak{j},\mathfrak{g}]^{*}$. Il en résulte que, si
$\mathbold{G}$ désigne un groupe algébrique connexe d'algèbre de Lie
$\mathfrak{g}$, l'application
$\psi:\mathbold{G}\times\mathfrak{g}^{*\mathfrak{j}}
\rightarrow\mathfrak{g}^{*}$
définie par $\psi(x,g)=x.g$, dont la différentielle en $(1,g)$ est
donnée par $\diff\psi_{(1,g)}(X,h)=X.g+h$, est étale en ce point, de
sorte que son image contient un voisinage de $g$ et donc des éléments
fortement réguliers de $\mathfrak{g}^{*}$. On en déduit qu'il existe
une forme fortement régulière $g'\in\mathfrak{g}^{*\mathfrak{j}}$, de
sorte que
$\mathfrak{g}(g)=
\mathfrak{j}\oplus\,^{u}\!\mathfrak{z}\subset\mathfrak{g}(g')$
et donc que $\mathfrak{g}(g)=\mathfrak{g}(g')$. D'où le fait que $g$
est fortement régulière.

d)  (i)  résulte de la démonstration de c) et du fait que ni
$\dim\mathfrak{g}(g)$, ni $\dim\,^{u}\!(\mathfrak{g}(g))$ ne dépend du
choix de la forme fortement régulière $g\in\mathfrak{g}^{*}$.

(ii) est conséquence de (i) et du fait que, si $g$ est une forme
régulière, $\g(g)$ est une algèbre de Lie commutative.

e) est conséquence de (i) $\Rightarrow$ (ii) du a), de b) et de d).

f) Soit $g\in\g^{*}_{red}$ telle que
$g_{\vert\mathfrak{r}_{\g}}=0$. D'après \cite[lemme
  1.20]{khalgui-torasso-1993}, il existe une forme linéaire $g'$ de
type unipotent $u$-équivalente à $g$ telle que
$\mathfrak{b}_{g}=\mathfrak{b}_{g'}$, $g$ et $g'$ aient même
restriction à $^{u}\!\mathfrak{b}_{g}$ et
$\mathfrak{r}_{g}\subset\mathfrak{r}_{g'}$. Compte tenu de b),
$\mathfrak{r}_{g}$ est une sous-algèbre réductive de rang maximal de
$\mathfrak{r}_{g'}$ de sorte que $\mathfrak{r}_{g'}=\mathfrak{r}_{g}+
[\mathfrak{r}_{g},\mathfrak{r}_{g'}]$. On en déduit que $g$ s'annule
sur $\mathfrak{r}_{g'}$, de sorte que
$g\in\cur{O}_{(0,g')}=\mathbold{G}.g'$.
\end{dem}

Choisissons une réalisation de $\mathfrak{g}$ comme une sous-algèbre
de Lie algébrique de l'algèbre de Lie des endomorphismes d'un
$\mathbb{K}$-espace vectoriel de dimension finie $V$ et désignons par
$\Tr X$ la trace d'un endomorphisme $X$ de $V$. Soit $\mathfrak{j}$
une sous-algèbre de Cartan-Duflo de $\g$ et $\h$ son centralisateur
dans $\g$. On désigne par $\h_{1}$ l'idéal de $\h$, orthogonal de
$\mathfrak{j}$ pour la forme bilinéaire $(X,Y)\mapsto\Tr XY$. Alors,
$\mathfrak{h}_{1}$ est une sous-algèbre de Lie algébrique de
$\mathfrak{h}$ (voir \cite[lemme 15]{khalgui-torasso-2006}) et l'on a
  $\h=\mathfrak{j}\oplus\h_{1}$. On en déduit la décomposition en
  somme directe duale $\h^{*}=\mathfrak{j}^{*}\oplus\h_{1}^{*}$, où
  l'on a identifié $\mathfrak{j}^{*}$ (resp $\mathfrak{h}_{1}^{*}$) à
  l'orthogonal de $\mathfrak{h}_{1}$ (resp. $\mathfrak{j}$).

Soit $g\in\g^{*}$ une forme de type réductif. D'après l'assertion e)
du théorème précédent, tout tore maximal de $\mathfrak{r}_{g}$ est une
sous-algèbre de Cartan-Duflo de $\g$.

\begin{theo}\label{theo1d2}
Soit $g\in\g^{*}$ une forme de type réductif et
$\mathfrak{j}\subset\mathfrak{r}_{g}$ une sous-algèbre de Cartan-Duflo
de $\g$. \'Ecrivons $g=\lambda+g_{1}$ avec
$\lambda\in\mathfrak{j}^{*}$ et $g_{1}\in\h^{*}_{1}$. Alors,

(i) $g_{1}$ est une forme de type réductif et unipotent
$u$-équivalente à $g$ telle que
$\mathfrak{b}_{g}=\mathfrak{b}_{g_{1}}$ et
$\mathfrak{r}_{g}\subset\mathfrak{r}_{g_{1}}$. En particulier,
$\mathfrak{j}$ est une sous-algèbre de Cartan de
$\mathfrak{r}_{g_{1}}$.

(ii) Les orbites coadjointes des formes de type réductif
$u$-équivalentes à $g$ sont les orbites des éléments de la forme
$\mu+g_{1}$, avec $\mu\in\mathfrak{j}^{*}$.

(iii) Soit $\mu\in\mathfrak{j}^{*}$. La forme linéaire $\mu+g_{1}$ est
fortement régulière si et seulement si $\mu$, vu comme un élément de
$\mathfrak{r}_{g_{1}}^{*}$, est une forme linéaire régulière.
\end{theo}
\begin{dem}
(i) Lorsque $g$ est fortement régulière, ce point résulte de
  \cite[lemme 1.26 et proposition 1.16]{khalgui-torasso-1993}. La
  démonstration dans le cas général est identique.

(ii) Soit $g'$ une forme de type réductif $u$-équivalente à $g$. Soit
$\mathfrak{j}'\subset\mathfrak{r}_{g'}$ une sous-algèbre de
Cartan-Duflo de $\g$. \'Ecrivons $g'=\lambda'+g'_{1}$ avec
$\lambda'\in\mathfrak{j}'^{*}$ et $g'_{1}\in\h'^{*}_{1}$. Comme $g$ et
$g'$ sont $u$-équivalentes, $g_{1}$ et $g'_{1}$ sont deux formes de
type unipotent $u$-équivalentes, donc situées sur la même orbite. On
peut donc supposer que $g'_{1}=g_{1}$. Mais alors d'après (i), on a
$\mathfrak{j}'\subset\mathfrak{r}_{g'}\subset\mathfrak{r}_{g_{1}}$. Quitte,
à faire agir un élément de $\mathbold{G}(g_{1})$, on peut supposer que
$\mathfrak{j}'=\mathfrak{j}$ ce qui montre que l'orbite de $g'$ est
bien de la forme indiquée.

Réciproquement, soit $\mu\in\mathfrak{j}^{*}$. D'après le résultat
précédent, il existe $\mu'\in\mathfrak{j}^{*}$ tel que la forme
$\mu'+g_{1}$ soit fortement régulière et de type réductif. On peut
donc appliquer \cite[lemme 1.28]{khalgui-torasso-1993} selon lequel
$\g(\mu+g_{1})=\mathfrak{r}_{g_{1}}(\mu)\oplus\,^{u}\!(\g(g_{1}))$. Il
est alors clair que la forme $\mu+g_{1}$ est de type réductif.

(iii) est une conséquence immédiate de ce qui précède.
\end{dem}

\subsection{Formes quasi-réductives et orbites de type réductif et
  unipotent.}\label{1e}

Soit $\mathfrak{g}$ une algèbre de Lie algébrique
et $\mathfrak{a}$ un idéal central de $\mathfrak{g}$. \'Etant donné
$\mu\in\mathfrak{a}^{*}$, on désigne par $\mathfrak{g}^{*}_{\mu}$
l'ensemble des $g\in\mathfrak{g}^{*}$ tels que
$g_{\mid\mathfrak{a}}=\mu$.

\begin{defi}
Soit $\mathfrak{g}$ une algèbre de Lie algébrique et $\mathfrak{a}$ un
idéal central de $\mathfrak{g}$. On dit qu'une forme linéaire
$\mu\in\mathfrak{a}^{*}$ est \emph{quasi-réductive}, sous-entendu
relativement à $\mathfrak{g}$, s'il existe une forme
$g\in\mathfrak{g}^{*}_{\mu}$ qui soit de type réductif.
\end{defi}

Si $\mathfrak{a}$ est un idéal central de l'algèbre de Lie algébrique
$\mathfrak{g}$, on désigne par $\mathfrak{a}^{*}_{\mathfrak{g},red}$
le sous-ensemble de $\mathfrak{a}^{*}$ constitué des formes
quasi-réductives.

\begin{theo}\label{theo1e1} Soit $\mathfrak{g}$ une algèbre de Lie
algébrique, $\mathfrak{z}$ son centre.

(i) si $\g$ est quasi-réductive et si $\mathfrak{a}$ est un idéal
central de $\mathfrak{g}$, l'ensemble
$\mathfrak{a}^{*}_{\mathfrak{g},red}$ est un ouvert de Zariski de
$\mathfrak{a}^{*}$, image par la projection naturelle de
$\mathfrak{g}^{*}$ sur $\mathfrak{a}^{*}$ de l'ouvert de Zariski
$\mathfrak{g}^{*}_{r}$ des formes fortement régulières sur
$\mathfrak{g}$.

Soit $\mu\in\,^{u}\!\mathfrak{z}^{*}$.

(ii) la forme $\mu$ est quasi-réductive si et seulement s'il existe une
forme $g\in\mathfrak{g}^{*}_{\mu}$ de type réductif et unipotent,

(iii) si $\mu$ est quasi-réductive, l'ensemble des
formes linéaires de type réductif et unipotent contenues dans
$\mathfrak{g}^{*}_{\mu}$ est une orbite sous l'action du groupe adjoint
connexe de $\mathfrak{g}$,

(iv) si $\mu$ est quasi-réductive, l'ensemble des éléments de type
unipotent et réductif de $\mathfrak{g}^{*}_{\mu}$ est égal à
l'ensemble des $g\in\mathfrak{g}^{*}_{\mu}$ de type réductif tels que
$\mathfrak{g}(g)$ soit de dimension maximale et
$g_{\vert\,^{r}\!\mathfrak{z}}=0$.
\end{theo}
\begin{dem}
(i) est une conséquence des assertions a) et d) (ii) du
  théorème \ref{theo1d1} et du fait que deux formes
  $u$-équivalentes ont même restriction à $^{u}\!\mathfrak{z}$.

(ii) résulte du fait déjà signalé que deux formes linéaires
  $u$-équivalentes ont même restriction à $^{u}\!\mathfrak{z}$ et du
  théorème \ref{theo1d1} a).%, il
  %est immédiat qu'il existe
  %$g\in\mathfrak{g}^{*}_{\mu}$ de type unipotent.

(iii) soit donc $g\in\mathfrak{g}^{*}_{\mu}$ de type unipotent et soit
  $\mathfrak{j}$ une sous-algèbre de Cartan-Duflo de $\g$. D'après le
  théorème \ref{theo1d2}, dont on reprend les notations, et comme
  manifestement $^{u}\!\mathfrak{z}$ est contenu dans $\h_{1}$,
  l'orbite de $g$ rencontre $(\h_{1}^{*})_{\mu}$ et toute forme de
  type réductif contenue dans $(\h_{1}^{*})_{\mu}$ est de type
  unipotent.

%Soit $\mathfrak{g}^{\mathfrak{j},1}$ l'idéal de
%$\mathfrak{g}^{\mathfrak{j}}$ orthogonal de $\mathfrak{j}$ pour la
%forme $(X,Y)\mapsto\Tr XY$. Alors $\mathfrak{g}^{\mathfrak{j},1}$
%contient $^{u}\!\mathfrak{z}$, de sorte que
%$(\mathfrak{g}^{\mathfrak{j},1})^{*}_{\mu}$ est défini. On a
%$\mathfrak{g}^{\mathfrak{j}}=\mathfrak{j}\oplus\mathfrak{g}^{\mathfrak{j},1}$
%ainsi que la décomposition duale $\mathfrak{g}^{*\mathfrak{j}}=
%\mathfrak{j}^{*}\oplus(\mathfrak{g}^{\mathfrak{j},1})^{*}$ et, il
%résulte de \cite[lemme 1.26]{khalgui-torasso-1993} que si l'on écrit
%$g'=\lambda+g_{1}$ avec $\lambda\in\mathfrak{j}^{*}$ et
%$g_{1}\in(\mathfrak{g}^{\mathfrak{j},1})^{*}$, $g_{1}$ est une forme
%de type unipotent $u$-équivalente à $g'$ et donc située dans l'orbite
%de $g$ sous l'action du groupe adjoint algébrique connexe de
%$\mathfrak{g}$. On peut donc supposer que
%$g\in(\mathfrak{g}^{\mathfrak{j},1})^{*}$ et

Il suffit donc de montrer que, si $\mathbold{G}$ est un groupe
algébrique connexe d'algèbre de Lie $\mathfrak{g}$, l'ensemble des
formes de type réductif situées dans $(\h_{1}^{*})_{\mu}$
constitue une orbite sous l'action du centralisateur $\mathbold{H}$
de $\mathfrak{j}$ dans $\mathbold{G}$, ce qui sera établi si l'on
montre que la $\mathbold{H}$-orbite d'une telle forme est ouverte dans
$(\mathfrak{h}_{1}^{*})_{\mu}$.

Soit donc $g\in(\mathfrak{h}_{1}^{*})_{\mu}$ une forme de type
réductif. Comme $\h^{*}$ est le sous-espace de $\g^{*}$ formé des
vecteurs invariants par l'action du tore $\mathfrak{j}$, $\h(g)$ est
la partie de $\g(g)=\mathfrak{r}_{g}\oplus\,^{u}\!\mathfrak{z}$
centralisée par $\mathfrak{j}$. Comme $\mathfrak{j}$ est une
sous-algèbre de Cartan de $\mathfrak{r}_{g}$, on a
$\h(g)=\mathfrak{j}\oplus\,^{u}\!\mathfrak{z}$. Il est alors clair que
que la $\mathbold{H}$-orbite de $g$ est ouverte dans
$(\mathfrak{h}_{1})^{*}_{\mu}$.

(iv) résulte de l'assertion (iii) et des assertions b) et f) du
théorème \ref{theo1d1}.
\end{dem}

\begin{co}\label{co1e1}
Soit $\mathfrak{g}$ une algèbre de Lie algébrique et $\mathfrak{z}$
son centre.  Par restriction à $^{u}\!\mathfrak{z}$, l'ensemble des
$\mathbold{G}$-orbites de type réductif et unipotent est en bijection
avec l'ensemble des $\mathbold{G}$-orbites dans
$^{u}\!\mathfrak{z}^{*}_{\mathfrak{g},red}$. En particulier, si
$\mu\in\!^{u}\!\mathfrak{z}^{*}_{\mathfrak{g},red}$, il existe une
unique $\mathbold{G}$-orbite coadjointe de type réductif et unipotent
rencontrant $\g^{*}_{\mu}$~: on la note $O_{\mathbold{G},\mu}$.
\end{co}
\begin{dem}
C'est une conséquence du théorème \ref{theo1e1} (iii).
\end{dem}

\subsection{Sous-algèbres et sous-groupes réductifs
  canoniques.}\label{1f}

Dans ce numéro, on suppose que $\mathfrak{g}$ est une algèbre de Lie
algébrique quasi-réductive de centre $\mathfrak{z}$. On rappelle que
si $g\in\mathfrak{g}^{*}$ est de type réductif, $\mathfrak{r}_{g}$
désigne le facteur réductif de $\mathfrak{g}(g)$~; on désignera alors
par $\mathbold{R}_{g}$ un facteur réductif de $\mathbold{G}(g)$. Les
facteurs réductifs $\mathbold{R}_{g}$ ont tous même composante neutre
et sont tous conjugués sous l'action de $^{u}\!\mathbold{Z}$.

%Supposons que $^{u}\!\mathbold{Z}$ soit un sous-groupe central de
%$\mathbold{G}$. Alors les sous-espaces affines $\g^{*}_{\mu}$,
%$\mu\in\,^{u}\!\mathfrak{z}^{*}$ sont $\mathbold{G}$-invariants et il
%suit du théorème \ref{theo1e1} (ii) que, si
%$\mu\in\,^{u}\!\mathfrak{z}^{*}_{\mathfrak{g},red}$, $\g^{*}_{\mu}$
%contient une unique orbite de type réductif et unipotent que l'on note
%$\cur{O}_{\mu}$.

\begin{theo}\label{theo1f1}
(i) les stabilisateurs $\g(g)$, pour $g\in\mathfrak{g}^{*}$ de type
  unipotent et réductif sont deux à deux conjugués par le groupe
  adjoint algébrique connexe.

(ii) il existe un ouvert de Zariski non vide $U$ de
  $^{u}\!\mathfrak{z}^{*}_{\mathfrak{g},red}$, tel que si $g$ et $g'$
  soient deux formes de type unipotent dont les restrictions à
  $^{u}\!\mathfrak{z}$ soient dans $U$, alors les stabilisateurs
  $\mathbold{G}(g)$ et $\mathbold{G}(g')$ sont
  $\mathbold{G}$-conjugués.

(iii) on suppose que $^{u}\!\mathbold{Z}$ est un sous-groupe central
  de $\mathbold{G}$. Alors, l'ensemble des
  $\mu\in\,^{u}\!\mathfrak{z}^{*}_{\mathfrak{g},red}$ tels que le
  nombre de composantes connexes de $\mathbold{G}(g)$ pour
  $g\in\g^{*}_{\mu}$ de type unipotent et réductif soit maximal est un ouvert
  de Zariski de $^{u}\!\mathfrak{z}^{*}$. Cet ouvert de Zariski
  vérifie l'assertion (ii).

\end{theo}
\begin{dem}
(i) Il est équivalent de montrer que les facteurs réductifs
  $\mathfrak{r}_{g}$, pour $g\in\mathfrak{g}^{*}$ de type unipotent et
  réductif sont deux à deux conjugués. Pour ce faire, il suffit de
  montrer que si $g$ est de type unipotent et réductif et si
  $\mu=g_{\mid^{u}\!\mathfrak{z}}$, il existe un ouvert de Zariski $U$
  de $^{u}\!\mathfrak{z}^{*}_{\mathfrak{g},red}$ contenant $\mu$ et
  tel que, pour tout $\mu'\in U$, il existe $g'$ de type unipotent et
  réductif tel que $\mu'=g'_{\mid^{u}\!\mathfrak{z}}$ et que
  $\mathfrak{r}_{g}$ et $\mathfrak{r}_{g'}$ soient conjugués par le
  groupe adjoint algébrique connexe.

En effet, si l'on suppose ceci établi, étant donnés $g_{i}$, $i=1,2$,
de type unipotent et réductif, il existe des ouverts de Zariski
$U_{i}$ de $^{u}\!\mathfrak{z}^{*}_{\mathfrak{g},red}$ contenant
$\mu_{i}=g_{i\mid^{u}\!\mathfrak{z}}$, $i=1,2$, et ayant la propriété
établie. Mais alors $U_{1}\cap U_{2}$ est non vide et, si $\mu'\in
U_{1}\cap U_{2}$, il existe pour $i=1,2$ $g'_{i}$ de type unipotent et
réductif tel que $\mu'=g'_{i\mid^{u}\!\mathfrak{z}}$ et que
$\mathfrak{r}_{g_{i}}$ et $\mathfrak{r}_{g'_{i}}$ soient conjugués par
le groupe adjoint algébrique connexe. Comme, d'après le théorème
\ref{theo1e1} $\mathfrak{r}_{g'_{1}}$ et $\mathfrak{r}_{g'_{2}}$ sont
conjugués par le groupe adjoint algébrique connexe, on voit que
$\mathfrak{r}_{g_{1}}$ et $\mathfrak{r}_{g_{2}}$ le sont.

Soit donc $g$ de type unipotent et réductif,
$\mu=g_{\mid^{u}\!\mathfrak{z}}$ et $g'$ une forme fortement régulière
$u$-équivalente à $g$, de sorte que
$\mathfrak{g}(g')=\mathfrak{j}\oplus\,^{u}\!\mathfrak{z}$ avec
$\mathfrak{j}$ une sous-algèbre de Cartan-Duflo de
$\mathfrak{g}$. Reprenant les notations de la démonstration du
théorème \ref{theo1e1}, on écrit
$\mathfrak{g}^{\mathfrak{j}}=\mathfrak{j}\oplus\mathfrak{g}^{\mathfrak{j},1}$
et $\mathfrak{g}^{*\mathfrak{j}}=
\mathfrak{j}^{*}\oplus(\mathfrak{g}^{\mathfrak{j},1})^{*}$. On sait
alors que $g'$ se décompose sous la forme $g'=\lambda+g_{1}$, avec
$\lambda\in\mathfrak{j}^{*}$ et
$g_{1}\in(\mathfrak{g}^{\mathfrak{j},1})^{*}$, $g_{1}$ étant de type
unipotent $u$-équivalente à $g'$ et donc dans l'orbite de $g$ sous le
groupe adjoint algébrique connexe. On peut donc supposer que
$g=g_{1}$, auquel cas $\mathfrak{j}$ est une sous-algèbre de Cartan de
$\mathfrak{r}_{g}$.

Soit $\mathfrak{m}$ un supplémentaire $\mathfrak{r}_{g}$-invariant de
$^{u}\!\mathfrak{z}$ dans $\mathfrak{g}$. %contenant
%$\mathfrak{r}_{g}$.
La décomposition
$\mathfrak{g}=\mathfrak{m}\oplus\,^{u}\!\mathfrak{z}$ permet
d'identifier $^{u}\!\mathfrak{z}^{*}$ à un sous-espace de
$\mathfrak{g}^{*}$ fixé point par point par $\mathfrak{r}_{g}$. De
plus, on a $\mathfrak{g}^{\mathfrak{j}}=
\mathfrak{m}^{\mathfrak{j}}\oplus\,^{u}\!\mathfrak{z}$ et il est clair
que $^{u}\!\mathfrak{z}$ est contenu dans
$\mathfrak{g}^{\mathfrak{j},1}$ de sorte que
$\mathfrak{g}^{\mathfrak{j},1}=
(\mathfrak{m}^{\mathfrak{j}}\cap\mathfrak{g}^{\mathfrak{j},1})
\oplus\,^{u}\!\mathfrak{z}$. Par suite, $^{u}\!\mathfrak{z}^{*}$ est
contenu dans $(\mathfrak{g}^{\mathfrak{j},1})^{*}$.

Maintenant si $\nu\in\,^{u}\!\mathfrak{z}^{*}$, on voit que
$\mathfrak{r}_{g}\oplus\,^{u}\!\mathfrak{z}\subset\mathfrak{g}(g+\nu)$
et $\mathfrak{j}\oplus\,^{u}\!\mathfrak{z}\subset\mathfrak{g}(g'+\nu)$
avec, dans les deux cas, égalité pour $\nu$ dans un ouvert de Zariski
$V$ de $^{u}\!\mathfrak{z}^{*}$ contenant $0$. Mais alors, si $\nu\in
V$, $g'+\nu$ est une forme fortement régulière et de type réductif
appartenant à $\mathfrak{g}^{*\mathfrak{j}}$ et $g'+\nu=\lambda+g+\nu$
avec $\lambda\in\mathfrak{j}^{*}$ et
$g+\nu\in(\mathfrak{g}^{\mathfrak{j},1})^{*}$, de sorte que $g+\nu$
est une forme de type unipotent et réductif $u$-équivalente à $g'$
dont le stabilisateur est
$\mathfrak{r}_{g}\oplus\,^{u}\!\mathfrak{z}$. Il est alors clair que
l'ouvert de Zariski $U=\mu+V$ de $^{u}\!\mathfrak{z}^{*}$ convient.

(ii) dans \cite[1.29, 1.30]{khalgui-torasso-1993} on a attaché à
chaque forme linéaire $g\in\g^{*}$ un sous-groupe algébrique
$\mathbold{B}_{g}$ de $\mathbold{G}$ d'algèbre de Lie
$\mathfrak{b}_{g}$, dont la classe de conjugaison ne dépend que de la
classe de $u$-équivalence de $g$ et qui vérifie
$\mathbold{B}_{g}=\mathbold{G}(g)(\mathbold{B}_{g})_{0}$ lorsque $g$
est de type unipotent. De plus, on a montré
(\cite[2.20]{khalgui-torasso-1993}) qu'il existe un ouvert de Zariski
non vide $V$ de $\g^{*}$ tel que si $g, g'\in V$ les facteurs
réductifs de $\mathbold{B}_{g}$ et $\mathbold{B}_{g'}$ sont
conjugués. Or, si la forme linéaire $g$ est de type réductif et
unipotent, $\mathbold{G}(g)=\mathbold{R}_{g}\,^{u}\!\mathbold{Z}$ et
$\mathbold{R}_{g}$ est un facteur réductif de $\mathbold{B}_{g}$. Il
est alors clair que l'ouvert de Zariski $U$, image de
$V\cap\g^{*}_{r}$ par la projection naturelle de $\g^{*}$ sur
$^{u}\!\mathfrak{z}^{*}$ convient.

(iii) soit $\mu\in\!^{u}\!\mathfrak{z}^{*}_{\mathfrak{g},red}$ telle
que le nombre de composantes connexes de $\mathbold{G}(g)$ pour
$g\in\g^{*}_{\mu}$ de type unipotent et réductif soit maximal. Soit
$\mathbold{R}$ un facteur réductif de $\mathbold{G}$ contenant
$\mathbold{R}_{g}$ et soit $\mathfrak{m}$ un supplémentaire
$\mathbold{R}$-invariant de $^{u}\!\mathfrak{z}$ dans $\g$. La
décomposition $\mathfrak{g}=\mathfrak{m}\oplus\,^{u}\!\mathfrak{z}$
permet d'identifier $^{u}\!\mathfrak{z}^{*}$ à un sous-espace de
$\mathfrak{g}^{*}$ fixé point par point par $\mathbold{R}$. Il suit de
la démonstration de (i), qu'il existe un ouvert de Zariski $V$ de
$^{u}\!\mathfrak{z}^{*}$ contenant $0$ tel que, pour tout $\nu\in V$,
la forme linéaire $g'=g+\nu$ est une forme de type réductif et
unipotent. Il est clair que, pour une telle forme, on a
$\mathbold{G}(g)\subset\mathbold{G}(g')$, ces deux groupes ayant même
composante neutre. Comme le nombre de composantes connexes de
$\mathbold{G}(g)$ est maximal, on a
$\mathbold{G}(g)=\mathbold{G}(g')$. Notre assertion est alors claire.
\end{dem}

\begin{defis}\label{defi1f1} Soit $\g$ une algèbre de Lie
algébrique quasi-réductive et $\mathbold{G}$ un groupe algébrique
d'algèbre de Lie $\mathfrak{g}$.

(i) les sous-algèbres $\mathfrak{r}_{g}$, pour $g\in\g^{*}_{red,u}$,
sont appelées les sous-algèbres \emph{réductives canoniques de} $\g$,

(ii) si $\mu\in\!^{u}\!\mathfrak{z}^{*}_{\mathfrak{g},red}$, les
sous-groupes $\mathbold{R}_{g}$, $g\in O_{\mathbold{G},\mu}$, sont appelés les
\emph{sous-groupes réductifs canoniques relatifs à} $\mu$ de $\mathbold{G}$.
\end{defis}

D'après le théorème \ref{theo1e1}, la classe de conjugaison, sous le
groupe adjoint connexe de $\g$, des sous-algèbres réductives
canoniques est uniquement déterminée et ne dépend que de l'algèbre de
Lie quasi-réductive $\g$.

Par contre la classe de conjugaison des sous-groupes
réductifs canoniques relatifs à $\mu$ peut dépendre de la
$\mathbold{G}$-orbite de $\mu$ dans
$^{u}\!\mathfrak{z}^{*}_{\mathfrak{g},red}$~: voir le numéro \ref{1i6}
ci-après.

%\subsection{}\label{1g}
Le résultat suivant montre que, pour
$\mu\in\!^{u}\!\mathfrak{z}^{*}_{\mathfrak{g},red}$ fixé, les
$\mathbold{G}$-orbites des formes de type réductif rencontrant le
sous-espace $\g^{*}_{\mu}$ sont en bijection avec les orbites de type
réductif d'un sous-groupe réductif canonique relatif à $\mu$.

On reprend les notations du numéro \ref{1b}, proposition \ref{pr1b1}.

\begin{theo}\label{theo1f2}
Soit $\mathbold{G}$ un groupe algébrique sur $\mathbb{K}$ d'algèbre de
Lie $\mathfrak{g}$ supposée quasi-réductive. Soit
$\mu\in\,^{u}\!\mathfrak{z}^{*}_{\g, red}$, $g\in\mathfrak{g}^{*}_{\mu}$
une forme de type réductif de type unipotent et $\mathbold{R}_{g}$ un
facteur réductif de $\mathbold{G}(g)$. Alors, l'application
$\lambda\mapsto\cur{O}_{g,\lambda+\mu}$ induit une bijection de
l'ensemble des $\mathbold{R}_{g}$-orbites de type réductif dans
$\mathfrak{r}_{g}^{*}$ sur l'ensemble des $\mathbold{G}$-orbites de type
réductif rencontrant $\mathfrak{g}^{*}_{\mu}$.
\end{theo}
\begin{dem}
Désignons par $\cur{U}_{g}$ le sous-ensemble,
$\mathbold{G}$-invariant, de $\mathfrak{g}^{*}$ constitué des formes
linéaires $u$-équivalentes à $g$. On a la décomposition en somme
directe $\mathfrak{g}(g)=\mathfrak{r}_{g}\oplus\,^{u}\!\mathfrak{z}$
et la décomposition duale $\mathfrak{g}(g)^{*}=
\mathfrak{r}_{g}^{*}\oplus\,^{u}\!\mathfrak{z}^{*}$. Il est immédiat
que l'application $\lambda\mapsto\lambda+\mu$ induit une bijection de
$\mathfrak{r}_{g}^{*}$ sur $L(g)$. Utilisant la proposition
\ref{pr1b1}, on en déduit facilement que l'application
$\lambda\mapsto\cur{O}_{g,\lambda+\mu}$ induit une bijection de
$\mathbold{R}_{g}\backslash\mathfrak{r}_{g}^{*}$ sur
$\mathbold{G}\backslash\cur{U}_{g}$.

Soit $\lambda\in L(g)$ et $g'\in\mathfrak{g}^{*}$ ayant même
restriction que $\lambda$ à $\mathfrak{r}_{g}$ et que $g$ à
$^{u}\!\mathfrak{b}_{g}$. Il résulte de l'assertion b) du théorème
\ref{theo1d1} que $g'$ est de type réductif, si et seulement si
$\lambda$ l'est. Le théorème est alors conséquence immédiate de ce
qui  précède  et du fait que, d'après le corollaire
\ref{co1e1}, il existe une et une seule orbite de type réductif et
unipotent rencontrant $\g^{*}_{\mu}$.
\end{dem}

\subsection{Exemples.}\label{1g}

Dans ce numéro, nous donnons quelques exemples d'algèbres de Lie
quasi-réductives et de formes de type réductif.

\subsubsection{}\label{1g1}
 Une algèbre de Lie algébrique réductive $\mathfrak{g}$ est
 quasi-réductive (voir
 la remarque \ref{rem1d1}).

 %Elle possède une unique forme de type réductif et
 %unipotent, savoir $0$ et l'on a $\mathfrak{r}_{0}=\mathfrak{g}$. De
 %plus si l'on identifie $\mathfrak{g}$ et $\mathfrak{g}^{*}$ au moyen
 %d'une forme bilinéaire symétrique non dégénérée invariante, les
 %formes de type réductif (resp. régulières de type réductif) sont
 %celles qui s'identifient aux éléments semi-simples
 %(resp. semi-simples réguliers de $\mathfrak{g}$).

\subsubsection{}\label{1g2}
Soit $V$ un espace vectoriel de dimension finie sur $\mathbb{K}$ muni
d'une forme symplectique $B$ et $\mathfrak{h}_{V}$ l'algèbre de
Heisenberg construite sur $V$~: $\mathfrak{h}_{V}=V\oplus\mathbb{K}z$,
$\mathfrak{z}=\,^{u}\!\mathfrak{z}=\mathbb{K}z$ est le centre de
$\mathfrak{h}_{V}$ et si $v,w\in V$, leur crochet est donné par
$[v,w]=B(v,w)z$. Alors $\mathfrak{h}_{V}$ est une algèbre de Lie
quasi-réductive, les formes linéaires de type réductif sont celles
dont la restriction au centre est non nulle et pour toute forme de
type réductif $g$ on a $\mathfrak{r}_{g}=0$. De plus, on a
$\mathfrak{z}^{*}_{\mathfrak{z},red}=\{\nu z^{*} :
\nu\in\mathbb{K}^×\}$ et si $g\in(\mathfrak{h}_{V})_{\nu z^{*}}$,
$\nu\in\mathbb{K}^×$, l'orbite de $g$ sous le groupe adjoint
connexe est $(\mathfrak{h}_{V})_{\nu z^{*}}$.

\subsubsection{}\label{1g3}
On garde les notations du numéro précédent et on désigne par
$\mathfrak{sp}(V)$ l'algèbre de Lie du groupe symplectique de $(V,B)$~:
l'action naturelle de $\mathfrak{sp}(V)$ sur $V$ se prolonge en une
action par dérivation sur $\mathfrak{h}_{V}$, triviale sur le centre. Soit
alors $\mathfrak{g}$ le produit semi-direct de $\mathfrak{sp}(V)$ par
$\mathfrak{h}_{V}$ de sorte que
$\mathfrak{g}=\mathfrak{sp}(V)\oplus\mathfrak{h}_{V}$, le crochet étant
donné par
\begin{equation*}
[X_{1}+h_{1},X_{2}+h_{2}]=[X_{1},X_{2}]+X_{1}.h_{2}-X_{2}.h_{1}+[h_{1},h_{2}]
\mbox{, }X_{1},X_{2}\in\mathfrak{sp}(V)\mbox{, }h_{1},h_{2}\in\mathfrak{h}_{V}.
\end{equation*}
Soit $\mathfrak{j}$ une sous-algèbre de Cartan de $\mathfrak{sp}(V)$. On
désigne par $g_{z^{*}}$ la forme linéaire sur $\mathfrak{g}$ nulle sur
$\mathfrak{sp}(V)\oplus V$ et telle que $g_{z^{*}}(z)=1$ et on identifie un
élément $\lambda$ de $\mathfrak{j}^{*}$ avec la forme linéaire sur
$\mathfrak{g}$ ayant même restriction à $\mathfrak{j}$ et nulle sur
$[\mathfrak{j},\mathfrak{sp}(V)]\oplus\mathfrak{h}_{V}$.

Alors $\mathfrak{g}$ est une algèbre quasi-réductive. Un système de
représentants des orbites de type réductif sous l'action du groupe
adjoint connexe est constitué des formes s'écrivant $\lambda+\nu
g_{z^{*}}$, $\lambda$ parcourant un système de représentants des orbites
du groupe de Weyl dans $\mathfrak{j}^{*}$ et
$\nu\in\mathbb{K}^×$, celles de type unipotent correspondant à
$\lambda=0$ et celles régulières correspondant à $\lambda$ régulier
dans $\mathfrak{sp}(V)^{*}$. Si $\nu\in\mathbb{K}^×$, on a
$\mathfrak{g}(\nu g_{z^{*}})=\mathfrak{sp}(V)\oplus\mathbb{K}z$,
$\mathfrak{r}_{\nu g_{z^{*}}}=\mathfrak{sp}(V)$ et si
$\lambda\in\mathfrak{j}^{*}$ est régulier, on a
$\mathfrak{g}(\lambda+\nu g_{z^{*}})=\mathfrak{j}\oplus\mathbb{K}z$,
$\mathfrak{r}_{\lambda+\nu g_{z^{*}}}=\mathfrak{j}$.

\subsubsection{}\label{1g4}
On garde les notations du numéro précédent, sauf pour
$\mathfrak{g}$. On fait agir $t\in\mathbb{K}$ au travers d'une
dérivation sur $\mathfrak{h}_{V}$ en décidant que $t.v=tv$, $v\in V$ et
$t.z=2tz$. Alors cette action commute à celle de $\mathfrak{sp}(V)$.  On
désigne maintenant par $\mathfrak{g}$ l'algèbre de Lie produit
semi-direct de $\mathfrak{sp}(V)\times\mathbb{K}$ par $\mathfrak{h}_{V}$,
qui contient $\mathfrak{sp}(V)\oplus\mathfrak{h}_{V}$, et on prolonge
$g_{z^{*}}$ en une forme linéaire sur $\mathfrak{g}$, encore notée
$g_{z^{*}}$, nulle sur le facteur $\mathbb{K}$ de
$\mathfrak{sp}(V)\times\mathbb{K}$.

Alors $\mathfrak{g}$ est une algèbre quasi-réductive et son centre est
trivial. L'ensemble des formes de type réductif et unipotent est
l'orbite de $g_{z^{*}}$ sous l'action du groupe adjoint connexe. Un
système de représentants des orbites sous l'action du groupe adjoint
connexe des formes de type réductif est constitué des formes
s'écrivant $\lambda+ g_{z^{*}}$, $\lambda$ parcourant un système de
représentants des orbites du groupe de Weyl dans
$\mathfrak{j}^{*}$. On a
$\mathfrak{r}_{g_{z^{*}}}=\mathfrak{g}(g_{z^{*}})=\mathfrak{sp}(V)$ et si
$\lambda\in\mathfrak{j}^{*}$ est régulier, on a $\mathfrak{g}(\lambda+
g_{z^{*}})=\mathfrak{j}$.

%\subsubsection{}\label{1g5}

\subsection{Algèbres de Lie préhomogènes.}\label{1h}

Une forme $g \in \g^*$ telle que $\g(g)=\z$ est évidemment de type
réductif. Il revient au même de demander que l'orbite $\mathbold{G}_{0}.g$ soit
un ouvert de l'espace affine $g+ \z^\bot$, où $\z^\bot$ désigne
l'orthogonal de $\z$ dans $\g^*$. Nous verrons que ces formes
jouent un rôle particulier dans la théorie des algèbres de Lie
quasi-réductives, ce qui justifie une nouvelle définition.

\begin{defis}\label{defi1h1}
(i) Une algèbre de Lie algébrique $\mathfrak{g}$ de centre
$\mathfrak{z}$ est dite \emph{préhomogène} s'il existe une forme linéaire
$g\in\mathfrak{g}^{*}$ telle que $\mathfrak{g}(g)=\mathfrak{z}$.

(ii) Si $\mathfrak{a}$ est un idéal central de l'algèbre de Lie
  préhomogène $\mathfrak{g}$, une forme $\mu\in\mathfrak{a}^{*}$ est
  dite \emph{préhomogène}, s'il existe $g\in\mathfrak{g}^{*}_{\mu}$
  telle que $\g(g)=\mathfrak{z}$.

(iii) Une algèbre de Lie préhomogène de centre nul est dite de
  Frobenius.
\end{defis}

\begin{rem}\label{rem1h1}
Lorsque $\z$ est nul, l'algèbre $\g$ est préhomogène si et
seulement si $\mathbold{G}$ a une orbite coadjointe ouverte, c'est-à-dire si
$\g^*$ est une espace préhomogène pour $\mathbold{G}$ au sens de  Kimura-Sato
(voir \cite{kimura-b-2003}). C'est pourquoi nous avons choisi cette
terminologie. Cependant, les algèbres de Lie préhomogènes sont
intéressantes à plusieurs points de vue et apparaissent dans la
littérature sous d'autres noms.  Voir notamment
%\citelist{\cite{ooms-1980} \cite{elashvili-1982}
%\cite{elashvili-ooms-2003}}
\cite{ooms-1980,elashvili-1982,elashvili-ooms-2003}
où elles sont appelées \emph{algèbres de
Frobenius} lorsque $\z$ est nul et \emph{de carré intégrable} en
général. Dans
%\citelist{\cite{khakimdjanov-goze-medina-2004} \cite{burde-2006}}
\cite{khakimdjanov-goze-medina-2004,burde-2006}, les algèbres de
Lie de la forme $\g/\z$ (où $\g$ est préhomogène) sont
 appelées \emph{symplectiques} ou bien
\emph{quasi-Frobenius}.
\end{rem}

Les algèbres de Lie préhomogènes sont quasi-réductives. Si
$\mathfrak{a}$ est un idéal central de l'algèbre de Lie préhomogène
$\g$, une forme linéaire $\mu\in\mathfrak{a}^{*}$ est préhomogène si
et seulement si elle est quasi-réductive. Plus
précisément, on a le résultat suivant~:
\begin{theo}\label{theo1h1}
Soit $\mathfrak{g}$ une algèbre de Lie préhomogène de centre $\mathfrak{z}$.

(i) les formes de type réductif sur $\mathfrak{g}$ sont les formes
régulières et elles sont toutes fortement régulières. Celles qui sont
de type unipotent sont celles qui s'annulent sur le facteur réductif
$^{r}\!\mathfrak{z}$de $\mathfrak{z}$,% et, pour une telle forme $g$,
%$\mathfrak{r}_{g}=\,^{r}\!\mathfrak{z}$,

(ii) si $\mu\in\mathfrak{z}^{*}$ est préhomogène, l'ensemble des
formes de type réductif contenues dans $\mathfrak{g}^{*}_{\mu}$ est
une orbite sous l'action du groupe adjoint connexe et c'est l'unique
orbite ouverte contenue dans $\mathfrak{g}^{*}_{\mu}$,

(iii) l'ensemble des $\mu\in\mathfrak{z}^{*}$ qui sont préhomogènes
est un ouvert de Zariski, image réciproque par la projection naturelle
de $\mathfrak{z}^{*}$ sur $^{u}\!\mathfrak{z}^{*}$ de
$^{u}\!\mathfrak{z}^{*}_{\mathfrak{g},red}$,

(iv) par restriction à $\mathfrak{z}$, l'ensemble des orbites
coadjointes de type réductif est en bijection avec les orbites de
$\mathbold{G}$ dans $\mathfrak{z}^{*}_{\g,red}$.
\end{theo}
\begin{dem}
Il est clair que $\mathfrak{g}$ admet $^{r}\!\mathfrak{z}$ pour unique
sous-algèbre de Cartan-Duflo. Il résulte alors de l'assertion e) du
théorème \ref{theo1d1} que toute forme de type réductif et unipotent
est telle que $\mathfrak{g}(g)=\mathfrak{z}$, donc
régulière. L'assertion (i) découle facilement de cette remarque et du
fait que d'après l'assertion a) du théorème \ref{theo1d1} si une
forme de type unipotent est $u$-équivalente à une forme de type
réductif, elle est elle-même de type réductif.

Soit $\mu\in\,^{u}\!\mathfrak{z}^{*}$ et $g\in\g^{*}_{\mu}$ une forme
de type réductif. Alors l'orbite de $g$ sous le groupe adjoint connexe
est contenue dans $\g^{*}_{\mu}$ et son espace tangent en $g$ est
$\mathfrak{z}^{\perp}$. L'assertion (ii) est alors claire.

Les assertions (iii) et (iv) sont des reformulations des résultats du
théorème \ref{theo1e1} et du corollaire \ref{co1e1}.
\end{dem}

\begin{co}\label{co1h1}
Soit $\g$ une algèbre de Lie préhomogène et nilpotente. Pour toute
forme linéaire quasi-réductive $\mu\in\,^{u}\!\mathfrak{z}^{*}$,
$\g^{*}_{\mu}$ est une orbite de type réductif sous l'action du groupe
adjoint connexe.
\end{co}
\begin{dem}
C'est une conséquence immédiate de l'assertion (ii) du théorème
précédent et du fait que les orbites coadjointes d'un groupe unipotent
sont fermées.
\end{dem}

Soit $\g$ une algèbre de Lie préhomogène. Alors $2e=2e_{\g}$ est la
dimension de $\g/\mathfrak{z}$. Fixons un élément non nul $\omega$ de
$\wedge^{2e}((\g/\mathfrak{z})^{*})$. Si $g\in\g^{*}$, la forme
bilinéaire $\beta_{g}$ passe au quotient à $\g/\mathfrak{z}$ en une
forme encore notée de même. On note $\Psi_{\g,\omega}$ (ou plus
simplement $\Psi_{\g}$) la fonction polynomiale sur $\g^{*}$ telle que
$\frac{1}{e!}\wedge^{e}(\beta_{g})=\Psi_{\g}(g)\omega$,
$g\in\g^{*}$. C'est un élément non nul de $S([\g,\g])$,
$\mathbold{G}$-semi-invariant de poids $\det\Ad_{\g}$.

On écrit
\begin{equation}\label{eq1h1}
  \Psi_{\g}=\Phi_{\g}\Psi_{\g,1}
\end{equation}
où $\Phi_{\g}$ est le produit de tous les facteurs premiers de
$\Psi_{\g}$ qui sont contenus dans $S(\mathfrak{z})$.

%{\color{red} Je pense que cela ne détermine pas la décomposition.
%Il faut dire que $\Phi_{\g}$ est le produit de tous les facteurs
%premiers de $\Psi_{\g}$ qui sont contenus dans $S(\mathfrak{z})$ }

On a alors le corollaire évident suivant du théorème \ref{theo1h1}~:
\begin{co}\label{co1h2}
Soit $\g$ une algèbre de Lie algébrique préhomogène. Alors,

(i) un élément $g\in\g^{*}$ est régulier si et seulement si
$\Psi_{\g}(g)\neq 0$

(ii) un élément $\mu$ de $\mathfrak{z}^{*}$ est préhomogène si et
seulement si $\Phi_{\g}(\mu)\neq 0$,

(iii) $\Phi_{\g}$ est un élément de $S(^{u}\!\mathfrak{z})$.
\end{co}

\begin{rem}\label{rem1h2}
  Lorsque $\g$ est une algèbre de Frobenius, on a $\Phi_{\g}=1$.
\end{rem}

\subsection{Exemples.}\label{1i}

Voici quelques exemples d'algèbres de Lie préhomogènes.

\subsubsection{}\label{1i1}
L'algèbre de Heisenberg $\mathfrak{h}_{V}$ du numéro \ref{1g2} est
préhomogène et l'on a (pour un choix convenable de $\omega$)
$\Psi_{\mathfrak{h}_{V}}=\Phi_{\mathfrak{h}_{V}}=z^{e}$, où $2e=\dim V$.

\subsubsection{}\label{1i2}
On fait agir $\mathbb{K}$ par dérivations sur $\mathfrak{h}_{V}$ comme au
numéro \ref{1g4} et on considère le produit semi-direct $\mathfrak{g}$
de $\mathbb{K}$ par $\mathfrak{h}_{V}$. Alors, $\mathfrak{g}$ est une
algèbre de Frobenius, l'ensemble des formes de type réductif étant
l'unique orbite ouverte du groupe adjoint connexe~: elle est
constituée des formes linéaires qui ne sont pas nulles sur le centre de
$\mathfrak{h}_{V}$. Ici, on a $\Psi_{\g}=z^{e+1}$.

% {\color{red} \medskip Il me semble que l'on peut enlever cet exemple (et
% aussi la référence à {\color{red} Raïs.}

%\subsubsection{}\label{1i3}
%Soit $V$ un espace vectoriel de dimension finie $r$ sur $\mathbb{K}$ et
%$\mathfrak{a}$ l'algèbre de Lie produit semi-direct $\mathfrak{gl}(V)$
%par $V$. Alors, $\mathfrak{a}$ est une algèbre de Frobenius,
%l'ensemble des formes de type réductif étant l'unique orbite ouverte du
%groupe adjoint connexe.

%Si $n$ est une forme linéaire sur $V$, on désigne par
%$\phi_{n}:V\rightarrow\mathfrak{a}$ l'application linéaire telle
%que, pour $x,y\in V$, $\phi_{n}(x)=(n(x)Id_{V},x)$. Soit $a\in\a^{*}$. On
%définit les formes linéaires $n_{k}\in V^{*}$, $1\leq k\leq r$ en
%décidant que~:
%\begin{gather}
%n_{1}=a_{\vert V}\\
%n_{k+1}=a\circ\phi_{n_{k}}\mbox{, }1\leq k\leq r-1,
%\end{gather}
%de sorte que, pour $1\leq k\leq r$, $n_{k}$ est une fonction
%polynomiale homogène de degré $k$ et $\mathrm{GL}(V)$-équivariante de
%$a$. On choisit une base de $V^{*}$ et on note $\det$ la fonction
%déterminant par rapport à cette base. Alors, on a
%$\Psi_{\a}(a)=\det(n_{1},\ldots,n_{r})$ (voir \cite[Théorème
%  3.8]{rais-1978}).}

\subsubsection{}\label{1i4}
Dans ce numéro et les deux suivants, $V=\mathbb{K}e\oplus\mathbb{K}f$
est un $\mathbb{K}$-espace vectoriel de dimension 2. On désigne par
$\mathfrak{b}$ la sous-algèbre de Borel de $\mathfrak{sl}(V)$ constituée
des endomorphismes ayant une matrice triangulaire supérieure dans la
base $(e,f)$ et par $\mathbold{B}$ le sous-groupe de Borel
correspondant de $\mathrm{SL}(V)$.  On note $h$ (resp. $x$)
l'endomorphisme de $V$ tel que $h.e=e$ et $h.f=f$ (resp. $x.e=0$ et
$x.f=e$), de sorte que $(h,x)$ est une base de $\mathfrak{b}$
vérifiant $[h,x]=2x$. Alors $\mathfrak{b}$ est une algèbre de
Frobenius. On a $\Psi_{\mathfrak{b}}=x$.

\subsubsection{}\label{1i5}
Soit $\mathfrak{f}$ le produit semi-direct de $\mathfrak{b}$ et $V$~:
$\mathfrak{f}=\mathfrak{b}\oplus V$. Outre celui déjà donné, les
crochets non nuls des éléments de la base $(h,x,e,f)$ sont~:
$[h,e]=e$, $[h,f]=-f$ et $[x,f]=e$. Alors $\mathfrak{f}$ est une
algèbre de Frobenius, l'orbite ouverte étant constituée des formes
linéaires $g$ telles que $g(e)\neq0$. On a $\Psi_{\mathfrak{f}}=e^{2}$.

\subsubsection{}\label{1i6}
On munit $V$ de la forme symplectique $B$ telle que $B(e,f)=1$ et on
considère l'algèbre de Heisenberg $\mathfrak{h}_{V}=V\oplus\mathbb{K}z$
correspondante. On désigne par $\mathbold{H}_{V}$ le groupe unipotent
d'algèbre de Lie $\mathfrak{h}_{V}$. On a $\mathfrak{sp}(V)=\mathfrak{sl}(V)$
et $\mathrm{Sp}(V)=\mathrm{SL}(V)$. On peut donc considérer le groupe
$\mathbold{G}$ produit semi-direct de $\mathbold{B}$ par
$\mathbold{H}_{V}$ dont l'algèbre de Lie est
$\mathfrak{g}=\mathfrak{b}\oplus\mathfrak{h}_{V}$. Alors outre ceux déjà
mentionés dans $\mathfrak{b}$ et $\mathfrak{h}_{V}$, les crochets non nuls
des éléments de la base $(h,x,e,f,z)$ de $\mathfrak{g}$ sont les
suivants~: $[h,e]=e$, $[h,f]=-f$, $[x,f]=e$. L'algèbre de Lie
$\mathfrak{g}$ est préhomogène de centre
$\mathfrak{z}=\mathbb{K}z$. De plus, le quotient
$\mathfrak{g}/\mathfrak{z}$ est isomorphe à l'algèbre
$\mathfrak{f}$. On a $\Psi_{\g}=e^{2}-2xz$ et $\Phi_{\g}=1$.

Le centre de $\mathbold{G}$ est le sous-groupe unipotent
$\mathbold{Z}$ d'algèbre de Lie $\mathfrak{z}$. Si $g\in\g^{*}$ est
une forme régulière, $\mathbold{G}(g)/\mathbold{Z}$ est un groupe
fini. \`A automorphisme extérieur près, il y a deux types de formes
régulières~: celles telles que $g(z)\neq0$, par exemple
$g=x^{*}+z^{*}$, et celles telles que $g(z)=0$, par exemple
$g=e^{*}$.  On vérifie que l'on a
$\mathbold{G}(x^{*}+z^{*})/\mathbold{Z}\simeq \mathbb{Z}/2\mathbb{Z}$
et $\mathbold{G}(e^{*})/\mathbold{Z}\simeq \{1\}$.

Ce dernier exemple illustre le fait que la classe d'isomorphie du
stabilisateur d'une forme de type réductif et unipotent peut dépendre
de cet élément et montre que les résultats du théorème \ref{theo1f1}
concernant ces classes d'isomorphie sont les meilleurs que l'on puisse
obtenir.

\subsection{Classification des orbites de type réductif.}\label{1j}

Dans ce numéro on ramène le problème de la classification des orbites
de type réductif pour une algèbre de Lie quasi-réductive au même
problème pour le centralisateur d'une sous-algèbre de Cartan-Duflo,
qui se trouve être une algèbre de Lie préhomogène.

Soit $\g$ une algèbre de Lie quasi-réductive et $\mathbold{G}$ un
groupe algébrique d'algèbre de Lie $\g$. Soit $\mathfrak{j}$ une
sous-algèbre de Cartan-Duflo de $\g$. On pose $\h=\g^{\mathfrak{j}}$,
$\mathbold{H}=Z_{\mathbold{G}}(\mathfrak{j})$ et
$\mathbold{H}'=N_{\mathbold{G}}(\mathfrak{j})$.  Rappelons
  que l'on désigne par $\h^{1}$ l'idéal de $\h$ orthogonal de
  $\mathfrak{j}$ pour la forme $(X,Y)\mapsto\Tr XY$, et que $
  \pi_{\mathfrak{j}} \in S(\h)$ a été défini au numéro \ref{1d},
  formule \ref{eq1c1}.

D'après \cite[2.28]{khalgui-torasso-1993}, le système de racines
$\Delta_{\mathfrak{j}}$ de $\mathfrak{j}$ dans $\mathfrak{r}$ ne
dépend pas de la sous-algèbre réductive canonique $\mathfrak{r}$ contenant
$\mathfrak{j}$, de même que, pour tout $\alpha\in\Delta_{\mathfrak{j}}$,
la coracine $H_{\alpha}$ dans la partie semi-simple de
$\mathfrak{r}$. On choisit un ensemble de racines positives
$\Delta^{+}_{\mathfrak{j}}\subset\Delta_{\mathfrak{j}}$ et on pose alors
\begin{equation}\label{eq1j1}
\pi_{\mathfrak{j},s}=\prod_{\alpha\in\Delta^{+}_{\mathfrak{j}}}H_{\alpha}.
\end{equation}
Alors, toujours d'après \cite[2.28]{khalgui-torasso-1993}, il existe
un élément $\pi_{\mathfrak{j},1}$ de $S(\h^{1})$ tel que
\begin{equation}\label{eq1j2}
  \pi_{\mathfrak{j}}=\pi_{\mathfrak{j},s}\pi_{\mathfrak{j},1}.
\end{equation}
On écrit alors
\begin{equation}\label{eq1j3}
\pi_{\mathfrak{j},1}=\pi_{\mathfrak{j},z}\pi_{\mathfrak{j},2}
\end{equation}
où $\pi_{\mathfrak{j},z}$ est le produit de tous les facteurs premiers
de $\pi_{\mathfrak{j},1}$ qui appartiennent à
$S(\mathfrak{z}_{\mathfrak{h}_{1}})$.

%{\color{red} Je pense que cela ne détermine pas la décomposition.
% }

Avec ces notations on
a~:
\begin{theo}\label{theo1j1}
(i) $\h$ est une algèbre de Lie préhomogène de centre
  $\mathfrak{z}_{\h}=\mathfrak{j}\oplus\,^{u}\!\mathfrak{z}$ et l'on a
  $\g^{*}_{red}\cap\h^{*}\subset\h^{*}_{red}$.

(ii) Si $g\in\g^{*}$ est une forme linéaire de type réductif,
  $\mathbold{G}.g\cap\h^{*}$ est une $\mathbold{H}'$-orbite de type
  réductif.

(iii) L'application $O\mapsto O\cap\h^{*}$ induit une bijection de
  l'ensemble des $\mathbold{G}$-orbites de type réductif (resp. de
  type réductif et unipotent) dans $\g^{*}$ sur l'ensemble des
  $\mathbold{H}'$-orbites de type réductif (resp. de type réductif et
  unipotent) relativement à $\h$ contenues dans
  $\h^{*}\cap\g^{*}_{red}$,

(iv) Une forme linéaire $g\in\h^{*}$ est de type réductif relativement
  à $\g$ si et seulement si elle est régulière relativement à $\h$ et
  elle vérifie $\pi_{\mathfrak{j},1}(g)\neq0$. Autrement dit, on a
\begin{equation}\label{eq1j4}
\begin{split}
\g^{*}_{red}\cap\h^{*}&=
\{g\in\h^{*}_{red}\vert\pi_{\mathfrak{j},1}(g)\neq0\}\\
&=\{g\in\h^{*}\vert\Psi_{\h}(g)\neq0\mbox{ et }\pi_{\mathfrak{j},1}(g)\neq0\}.
\end{split}
\end{equation}
Une telle forme est de type unipotent si et seulement si sa
restriction à $\mathfrak{j}$ est nulle.

(v) Une forme linéaire $\mu\in\,^{u}\!\mathfrak{z}^{*}$ est
quasi-réductive si et seulement si elle est préhomogène relativement à
$\h$ et elle vérifie $\pi_{\mathfrak{j},z}(\mu)\neq0$. Autrement dit, on a
\begin{equation}\label{eq1j5}
\begin{split}
^{u}\!\mathfrak{z}^{*}_{\g,red}&=
\{\mu\in^{u}\!\mathfrak{z}^{*}_{\h,red}\vert\pi_{\mathfrak{j},z}(\mu)\neq0\}\\
&=\{\mu\in^{u}\!\mathfrak{z}^{*}
\vert\Phi_{\h}(\mu)\neq0\mbox{ et }\pi_{\mathfrak{j},z}(\mu)\neq0\}.
\end{split}
\end{equation}

(vi) Soit $g\in\g^{*}$ une forme de type réductif et unipotent et
$\mu$ la restriction de $g$ à $^{u}\!\mathfrak{z}$. Alors
$\mathbold{G}.\mu=\mathbold{H}'.\mu$ et $\mathbold{G}.g\cap\h^{*}$ est
la réunion des orbites $O_{\mathbold{H}_{0},\nu}$, pour $\nu$
parcourant $\mathbold{G}.\mu$.

(vii) On suppose que $^{u}\!\mathbold{Z}$ est un sous-groupe central
de $\mathbold{G}$. Alors, pour toute forme linéaire $g\in\h^{*}$ de
type réductif et unipotent relativement à $\g$, on a
$\mathbold{H}'=\mathbold{N}_{\mathbold{R}_{g}}(\mathfrak{j})\mathbold{H}_{0}$.

(viii) Dans le cas général, pour toute forme linéaire $g\in\h^{*}$ de
type réductif et unipotent relativement à $\g$, on a
$\mathbold{H}'=
(\mathbold{N}_{\mathbold{G}}(\mathfrak{r}_{g})\cap\mathbold{H}')\mathbold{H}_{0}$.
\end{theo}
\begin{dem}
(i) Il est clair que $\mathfrak{j}\oplus\,^{u}\!\mathfrak{z}$ est une
  sous-algèbre de Lie algébrique centrale de $\h$. Cela dit, d'après
  le numéro \ref{1c}, $\h^{*}\cap\g^{*}_{r}$ est non vide tandis que,
  d'après le théorème \ref{theo1d1} d) (i),
  $\g^{*}_{r}\subset\g^{*}_{red}$. Comme $\h$ stabilise l'orthogonal
  $[\mathfrak{j},\g]$ de $\h^{*}$, on voit que $\h(g)\subset\g(g)$,
  pour tout $g\in\h^{*}$. Si $g\in\h^{*}\cap\g^{*}_{r}$, on a
  donc $\mathfrak{j}\oplus\,^{u}\!\mathfrak{z}\subset\h(g)
  \subset\g(g)=\mathfrak{j}\oplus\,^{u}\!\mathfrak{z}$. D'où
  l'assertion.

(ii) Soit $g\in\g^{*}_{red}$. Il résulte de l'assertion e) du théorème
  \ref{theo1d1} que $\mathbold{G}.g$ rencontre $\h^{*}$. On peut donc
  supposer que $g\in\h^{*}\cap\g^{*}_{red}\subset\h^{*}_{red}$. Soit
  $x\in\mathbold{G}$ tel que $g'=x.g\in\h^{*}$. Il suit encore de
  l'assertion e) du théorème \ref{theo1d1} que $\mathfrak{j}$ est une
  sous-algèbre de Cartan de $\mathfrak{r}_{g'}$. Mais alors, $\Ad
  x^{-1}(\mathfrak{j})$ est une sous-algèbre de Cartan de
  $\mathfrak{r}_{g}$ et il existe $y\in(\mathbold{R}_{g})_{0}$ tel que
  $\Ad y(\mathfrak{j})=\Ad x^{-1}(\mathfrak{j})$. On voit donc que
  $xy\in\mathbold{H}'$ et vérifie $g'=xy.g$. D'où notre assertion.

(iii) C'est une conséquence immédiate de l'assertion (ii) et du fait
  qu'une forme linéaire de type réductif $g$ est de type unipotent si
  et seulement si il existe une sous-algèbre de Cartan de
  $\mathfrak{r}_{g}$ sur laquelle elle s'annule (voir la remarque
  \ref{rem1b1}).

(iv) Soit $g\in\h^{*}_{red}$ que l'on écrit $g=\lambda+g_{1}$, avec
  $\lambda\in\mathfrak{j}^{*}$ et $g_{1}\in\h_{1}^{*}$.

Supposons que $g$ soit de type réductif relativement à $\g$. Alors, il
résulte du théorème \ref{theo1d2} que $g_{1}$ est de type réductif et
unipotent, $u$-équivalente à $g$ et qu'il existe
$\lambda'\in\mathfrak{j}^{*}$ tel que $g'=\lambda'+g_{1}$ soit une forme
fortement régulière relativement à $\g$. D'après le théorème
\ref{theo1c1}, on a $\pi_{\mathfrak{j}}(g')\neq0$, ce qui, d'après
\ref{eq1j2}, entraîne que $\pi_{\mathfrak{j},1}(g')\neq0$. Comme
$\pi_{\mathfrak{j},1}$ est un élément de $S(\h_{1})$, on voit que
$\pi_{\mathfrak{j},1}(g)=\pi_{\mathfrak{j},1}(g_{1})=\pi_{\mathfrak{j},1}(g')\neq0$.

Réciproquement, supposons que $\pi_{\mathfrak{j},1}(g)\neq0$ et soit
$\lambda'\in\mathfrak{j}^{*}$ tel que
$\pi_{\mathfrak{j},s}(\lambda')\neq0$. Alors, la forme
$g'=\lambda'+g_{1}$ est régulière relativement à $\h$ et vérifie
$\pi_{\mathfrak{j}}(g')\neq0$. D'après le théorème \ref{theo1c1}, elle
est fortement régulière relativement à $\g$. Il suit alors du théorème
\ref{theo1d2} que $g$ est de type réductif relativement à $\g$.

(v) C'est une conséquence immédiate de (iv) et des définitions.

(vi) C'est une conséquence de (iii) et du théorème \ref{theo1e1} (iii).

(vii) C'est une conséquence immédiate de (vi).
%Soit $g\in\h^{*}$ une forme linéaire de type réductif et
%unipotent relativement à $\g$. Comme $^{u}\!\mathbold{Z}$ est central
%dans $\mathbold{G}$, il suit de (vi) que
%$\mathbold{H}'.g=\mathbold{H}_{0}.g$. D'où le résultat.

(viii) Soit $g\in\h^{*}$ une forme linéaire de type réductif et
unipotent relativement à $\g$. Il suit de (vii) que
$\mathbold{H}'\cap\mathbold{G}_{0}\subset
\mathbold{N}_{\mathbold{R}_{g}}(\mathfrak{j})\mathbold{H}_{0}$. Il
suffit donc de montrer que
$\mathbold{H}'=(\mathbold{N}_{\mathbold{G}}(\mathfrak{r}_{g})\cap\mathbold{H}')
(\mathbold{H}'\cap\mathbold{G}_{0})$. Soit donc
$x\in\mathbold{H}'$. Comme les sous-algèbres réductives canoniques
sont toutes conjuguées par le groupe adjoint connexe, il existe
$y\in\mathbold{G}_{0}$ tel que
$yx\in\mathbold{N}_{\mathbold{G}}(\mathfrak{r}_{g})$. Quitte à
multiplier $y$ à gauche par un élément de $(\mathbold{R}_{g})_{0}$, on
peut supposer que $yx$ normalise aussi $\mathfrak{j}$, de sorte que
$yx\in\mathbold{N}_{\mathbold{G}}(\mathfrak{r}_{g})\cap\mathbold{H}'$. Notre
résultat est alors clair.
\end{dem}

\subsection{Exemples.}\label{1k}

Donnons quelques exemples pour illustrer les notions et résultats du
numéro précédent. Dans ce qui suit $\mathfrak{j}$ désigne une
sous-algèbre de Cartan-Duflo de l'algèbre de Lie quasi-réductive $\g$
et $\h$ son centralisateur dans $\g$.

\subsubsection{}\label{1k1}
Si $\g$ est réductive, une sous-algèbre de Cartan-Duflo est une
sous-algèbre de Cartan. Alors, $\h=\mathfrak{j}$ est une algèbre de
Lie commutative réductive et
$\pi_{\mathfrak{j}}=\pi_{\mathfrak{j},s}$. Le quotient
$\mathbold{H}'/\mathbold{H}$ est le groupe de Weil
$W_{\mathbold{G}}(\mathfrak{j})$ de $\mathbold{G}$ dans $\mathfrak{j}$
et l'ensemble des $\mathbold{G}$-orbites de type réductif dans $\g$
s'identifie, au moyen de l'application $O\mapsto
O\cap\mathfrak{j}^{*}$, à l'ensemble des
$W_{\mathbold{G}}(\mathfrak{j})$-orbites dans $\mathfrak{j}^{*}$.

\subsubsection{}\label{1k2}
Si $\g$ est nilpotente ou plus généralement préhomogène,
$\mathfrak{j}=\{0\}$, $\h=\g$ et $\pi_{j}=1$.

\subsubsection{}\label{1k3}
On reprend l'exemple du numéro \ref{1g4}~:
$\g=\mathfrak{sp}(V)\oplus\mathfrak{h}_{V}$, où $V$ est un espace
vectoriel de dimension $2r$ sur $\mathbb{K}$ muni de la forme symplectique
$B$. Dans ce cas, on peut prendre pour $\mathfrak{j}$ une sous-algèbre
de Cartan de $\mathfrak{sp}(V)$. Alors
$\h=\mathfrak{j}\oplus\mathbb{K}z$ est une algèbre de Lie
commutative. Alors, $\Phi_{\h}=\Psi_{\h}=1$, $\pi_{\mathfrak{j},s}$
est le produit des coracines d'un système de racines positives de
$\mathfrak{j}$ dans $\mathfrak{sp}(V)$,
$\pi_{\mathfrak{j}}=z\pi_{\mathfrak{j},s}$ et
$\pi_{\mathfrak{j},1}=\pi_{\mathfrak{j},z}=z$. On a
$^{u}\!\mathfrak{z}^{*}_{\g,red}=\mathbb{K}^×z^{*}
\subsetneqq\,^{u}\!\mathfrak{z}^{*}_{\h,red}=\mathbb{K}z^{*}$.

Soit $\mathbold{S}$ le groupe des automorphismes linéaires $x$ de $V$ qui
conservent la forme symplectique $B$ au signe près~:
\begin{equation*}
B(x.v,x.w)=\epsilon(x)B(v,w)\mbox{, }v,w\in V.
\end{equation*}
Alors $\epsilon$ est un caractère d'ordre $2$ de $\mathbold{S}$ dont
le noyau est $\mathbold{S}_{0}=\mathrm{Sp}(V)$. Le groupe
$\mathbold{S}$ opère par automorphismes dans l'algèbre de Lie de
Heisenberg $\mathfrak{h}_{V}$ de la manière suivante~:
\begin{equation*}
x.(v+\lambda z)=x.v+\epsilon(x)\lambda z\mbox{, }x\in\mathbold{S},v\in
V,\lambda\in\mathbb{K}.
\end{equation*}

On peut donc considérer le groupe algébrique $\mathbold{G}$, produit
semi-direct de $\mathbold{S}$ par le groupe unipotent
$\mathbold{H}_{V}$ d'algèbre de Lie $\mathfrak{h}_{V}$.

Soit $(e_{1},\ldots,e_{r},f_{1},\ldots,f_{r})$ une base symplectique
de $V$~: on a $B(e_{i},e_{j})=B(f_{i},f_{j})=0$ et
$B(e_{i},f_{j})=\delta_{ij}$, $1\leq i,j\leq r$. Soit $\tau$
l'endomorphisme linéaire de $V$ qui vérifie $\tau.e_{i}=f_{i}$ et
$\tau.f_{i}=e_{i}$, $1\leq i\leq r$. Alors $\tau\in\mathbold{S}$,
$\epsilon(\tau)=-1$ et
$\mathbold{G}=\mathbold{G}_{0}\sqcup\tau\mathbold{G}_{0}$. Le centre
$\mathbold{Z}$ de $\mathbold{G}_{0}$ est connexe et égal au centre de
$\mathbold{H}_{V}$.

Les orbites de $\mathbold{G}$ dans $^{u}\!\mathfrak{z}^{*}_{\g,red}$
sont les $\{\nu z^{*},-\nu z^{*}\}$, $\nu\in\mathbb{K}^×$. Si
$\nu\in\mathbb{K}^×$, l'orbite de type réductif et
unipotent qui rencontre $\g^{*}_{\nu z^{*}}$ est celle de la forme
$\nu g_{z^{*}}$ et l'on a $R_{\nu g_{z^{*}}}=\mathbold{S}_{0}$.

 Soit $\mathfrak{j}$ la sous-algèbre de Cartan de $\mathfrak{sp}(V)$
 admettant $(e_{1},\ldots,e_{r},f_{1},\ldots,f_{r})$ comme base de
 diagonalisation. Soit $\mathbold{J}$ (resp. $\mathbold{J}'$) le
 centralisateur (resp. normalisateur) de $\mathfrak{j}$ dans
 $\mathbold{S}_{0}$. On a alors
 $\mathbold{H}=\mathbold{J}\mathbold{Z}$ et $\mathbold{H}'=
 \mathbold{J}'\mathbold{Z}\sqcup\tau\mathbold{J}'\mathbold{Z}$.

 Soit $\nu\in\mathbb{K}^×$. On voit donc que
 $\mathbold{N}_{\mathbold{R}_{\nu
     g_{z^{*}}}}(\mathfrak{j})\mathbold{H}
 =\mathbold{J}'\mathbold{Z}\subsetneqq \mathbold{H}'$.

\subsection{Algèbres de Lie résolubles}\label{1l}

Dans ce numéro, $\g$ est une algèbre de Lie algébrique, $\mathfrak{z}$
son centre, $\mathbold{G}$ un groupe algébrique connexe d'algèbre de
Lie $\mathfrak{g}$ et $\mathbold{Z}$ la composante neutre du centre de
$\mathbold{G}$.

\begin{lem}\label{lem1l1}
Soit $\a$ un idéal de $\g$, $g$ un élément de $\g^*$ et $\nu $ la
restriction de $g$ à $\a\cap \z$.
%et
%soit $\a^*_\nu$ le sous-espace affine de $\a^*$ formé des éléments
%dont la restriction à $\z\cap \a$ est égale à $\nu$
Alors la projection de l'orbite $\mathbold{G}.g$ sur $\a^*$ est
ouverte dans $\a^*_\nu$ si et seulement si l'on a $\g(g)\cap \a
\subset\z$.
\end{lem}
\begin{dem}
En effet, chacune de ces conditions est équivalente à ce que
l'application de restriction est une surjection de l'espace
tangent $\g.g \subset \g^*$ à $\mathbold{G}.g$ en $g$ sur l'orthogonal  de
$\a\cap \z$ dans $\a^*$.
\end{dem}
\begin{co}\label{co1l1}
Soit $\mathfrak{a}$ un idéal de $\g$ contenu dans $^{u}\!\mathfrak{g}$
et contenant $^{u}\!\mathfrak{z}$ et $g\in\mathfrak{g}^{*}$ une forme
de type réductif. Alors, si $\nu=g_{\mid^{u}\!\mathfrak{z}}$ et
$a=g_{\mid\mathfrak{a}}$, l'orbite de $a$ sous $\mathbold{G}$ est
ouverte dans $\mathfrak{a}^{*}_{\nu}$.
\end{co}
\begin{dem}
Ceci résulte de ce que $\g(g)\cap \a$ est contenu dans
${}^u\!(\g(g))$, et donc dans $\a\cap \z$ par définition d'une forme
de type réductif.
\end{dem}
\begin{co}\label{co1l2}
 Soit $\nu \in
{}^u\!\z^*$ une forme quasi-réductive. Alors, pour tout idéal
$\mathfrak{a}$ de $\g$ contenu dans $^{u}\!\mathfrak{g}$ et
contenant $^{u}\!\mathfrak{z}$, l'espace affine $\a^{*}_{\nu}$ est
préhomogène sous l'action de $\mathbold{G}/\mathbold{Z}$.
\end{co}

Dans la suite de ce numéro, on suppose que $\g$ est résoluble.
\begin{theo}\label{theo1l1}
On suppose que $\g$ est résoluble. Soit $g\in\g^{*}$ et $\nu=g_{\vert
  ^{u}\!\mathfrak{z}}$. Alors, les assertions suivantes sont
équivalentes~:

(i) $g$ est de type réductif,

(ii) $\g(g)\cap ^{u}\!\g\subset\mathfrak{z}$,

(iii) la projection de $\mathbold{G}.g$ sur $^{u}\!\g^{*}$ est un
ouvert dans $(^{u}\!\g^{*})_{\nu}$.
\end{theo}
\begin{dem}
(i) $\Rightarrow$ (ii), par définition des formes de type
  réductif. Comme $\g$ est une algèbre de Lie résoluble, on a
  $^{u}\!(\g(g))=( ^{u}\!\g)(g)$~; d'où la réciproque. L'équivalence
  (ii) $\Leftrightarrow$ (iii) résulte du lemme \ref{lem1l1} appliqué
  à l'idéal $\a= ^{u}\!\g$.
\end{dem}

Le corollaire suivant est évident~:
\begin{co}\label{co1l3}
On suppose que $\g$ est résoluble. Soit $\nu \in
{}^u\!\z^*$. Alors $\nu$ est une forme quasi-réductive si et
seulement si l'espace affine $( ^{u}\!\g^{*})_{\nu}$ est
préhomogène sous l'action de $\mathbold{G}/\mathbold{Z}$.
\end{co}
\begin{theo}\label{theo1l2}
On suppose que $\g$ est résoluble. Soit $\nu \in {}^u\!\z^*$ une forme
quasi-réductive et soit $\Omega\subset ( ^{u}\!\g^{*})_{\nu}$ la
$\mathbold{G}$-orbite ouverte. Si $g\in \g^{*}_{\nu}$, les assertions
suivantes sont équivalentes~:

(i) $g$ est de type réductif,

(ii) $\g(g)\cap ^{u}\!\g= ^{u}\!\mathfrak{z}$,

(iii) la restriction de $g$ à $^{u}\!\g$ est contenue dans $\Omega$,

(iv) $g$ est fortement régulière.
\end{theo}
\begin{dem}
L'équivalence des assertions (i), (ii) et (iii) résulte du théorème
précédent. L'équivalence de (i) et (iv) résulte du théorème \ref{theo1d1} d).
\end{dem}
\subsection{Exemples.}\label{1m}
\subsubsection{}\label{1m1} Soit $\s$ une algèbre de Lie semi-simple,
et soit $\mathfrak{b}\subset \s$ une sous-algèbre de Borel. C'est
un résultat bien connu de Kostant que $\mathfrak{b}$ est
quasi-réductive (voir \cite{kostant} et \cite{joseph-1977}).

\subsubsection{}\label{1m2} En considérant  la sous-algèbre de
Borel $\mathfrak{g}$ de $\mathfrak{sl}(3,\mathbb{K})$ constituée des
matrices triangulaires supérieures, on obtient une algèbre de Lie
quasi-réductive de centre trivial, de base
$h_{1},h_{2},x_{1},x_{2},z$ pour laquelle les crochets non nuls sont
$[h_{i},x_{j}]=2\delta_{ij}x_{j}$, $[h_{i},z]=2z$ et
$[x_{1},x_{2}]=z$. On note $\mathbold{G}$ le sous-groupe de Borel
correspondant. Le radical unipotent $^{u}\!\mathfrak{g}$ est
l'algèbre de Heisenberg de base $x_{1},x_{2},z$.

On a donc $\mathfrak{z}^{*}=\{0\}$, $\g^{*}_{0}=\g^{*}$ et
  $^{u}\!\mathfrak{g}^{*}_{0}=\,^{u}\!\mathfrak{g}^{*}$. On vérifie que
  l'ensemble des formes linéaires $u\in\,^{u}\!\mathfrak{g}^{*}$ qui
  vérifient $u(z)\neq0$ est une $\mathbold{G}-$orbite ouverte. Ainsi,
  toute  forme linéaire $g\in\g^{*}$ telle $g(z)\neq 0$ est de type
  réductif. Ces formes sont fortement régulières et leur orbite est de
  dimension 4.

On rappelle l'algèbre de Frobenius $\mathfrak{b}$ introduite au numéro
\ref{1i4}. L'algèbre de Lie quotient $\g/\mathbb{K}z$ est isomorphe au
produit $\mathfrak{b}\times\mathfrak{b}$. On en déduit facilement
qu'une forme linéaire $g\in\g^{*}$ telle que $g(z)=0$ est régulière si
et seulement si $g(x_{1})g(x_{2})\neq0$. Une telle forme $g$ vérifie
$g(g)=\mathbb{K}z$~; elle n'est donc pas de type réductif. On voit
donc que le théorème \ref{theo1l2} n'est plus vrai si l'on y remplace
dans l'assertion (iv) \flqq fortement régulière\frqq~ par \flqq
régulière\frqq.

\section{Algèbres de Lie quasi-réductives unimodulaires}\label{2}

Les algèbres de Lie quasi-réductives et unimodulaires ont une
structure particulièrement simple que nous allons décrire. Nos
résultats généralisent ceux de Anh qui concernent les groupes
algébriques sur $\mathbb{R}$ unimodulaires qui admettent une série
discrète modulo le centre (voir \cite{anh-1978} et \cite{anh-1980}).
Ils ont été énoncés dans \cite{duflo-1982} sans démonstration. Comme
ceci est important pour cet article, il nous semble opportun d'en
écrire une.

%\subsection{}\label{2a}
%On commence par le résultat suivant~:
%\begin{lem}\label{lem2a1}
%Soit $\mathfrak{g}$ une algèbre de Lie quasi-réductive, $\mathfrak{z}$ son
%centre, $\mathbold{G}$ un groupe algébrique connexe d'algèbre de Lie
%$\mathfrak{g}$, $\mathfrak{a}\subset\,^{u}\!\mathfrak{g}$ un idéal
%$\mathbold{G}$-invariant contenant $^{u}\!\mathfrak{z}$ et
%$g\in\mathfrak{g}^{*}$ une forme réductive. Alors, si
%$\mu=g_{\mid^{u}\!\mathfrak{z}}$ et $a=g_{\mid\mathfrak{a}}$, l'orbite de $a$
%sous $\mathbold{G}$ est ouverte dans $\mathfrak{a}^{*}_{\mu}$.
%\end{lem}
%\begin{dem}
%La forme bilinéaire $(X,Y)\mapsto\langle a,[X,Y]\rangle$ induit une
%dualité entre $\mathfrak{g}/\mathfrak{g}(a)$ et
%$\mathfrak{a}/\mathfrak{a}(g)$. Or,
%$\mathfrak{a}(g)=\mathfrak{a}\cap\mathfrak{g}(g)=\,^{u}\!\mathfrak{z}$
%puisque la forme $g$ est réductive, de sorte que
%$\dim\mathbold{G}.a=
%\dim\mathfrak{a}/^{u}\!\mathfrak{z}=\dim\mathfrak{a}^{*}_{\mu}$.
%Comme $\mathbold{G}.a$ est contenu dans $\mathfrak{a}^{*}_{\mu}$, le
%résultat est clair.
%\end{dem}

\subsection{}\label{2b}
Voici le résultat principal de ce numéro~:
\begin{theo}\label{theo2b1}
Soit $\mathfrak{g}$ une algèbre de Lie algébrique unimodulaire et
quasi-réductive. On désigne par $\mathfrak{z}$ le centre de
$\mathfrak{g}$ et par $\mathfrak{z}_{u}$ celui de
$^{u}\!\mathfrak{g}$. Soit $g\in\mathfrak{g}^{*}$ une forme linéaire
de type unipotent et réductif dont on désigne par $u$ la restriction
à $^{u}\!\mathfrak{g}$. Alors, on a

(i)
$\mathfrak{z}_{u}=\,^{u}\!(\mathfrak{g}(g))=
\,^{u}\!\mathfrak{g}(u)=\,^{u}\!\mathfrak{z}$,

(ii) $\mathfrak{g}(g)=\mathfrak{g}(u)$,

(iii) $\mathfrak{g}=\mathfrak{g}(u)+\,^{u}\!\mathfrak{g}$.
\end{theo}
\begin{dem}
On va raisonner par récurrence sur la dimension de $\mathfrak{g}$. Le
résultat étant évident si $\mathfrak{g}=0$, on suppose que
$\dim\mathfrak{g}>0$ et le résultat vrai en dimension strictement
inférieure.

\subsection{}\label{2c}
Dans ce numéro nous décrivons une situation dans laquelle
l'hypothèse de récurrence s'applique. Soit $\mathbold{G}$ un groupe
algébrique connexe d'algèbre de Lie $\mathfrak{g}$ et soit
$\mathfrak{a}\subset\,^{u}\!\mathfrak{g}$ un idéal
$\mathbold{G}$-invariant contenant strictement
$^{u}\!\mathfrak{z}$. On pose $a=g_{\mid\mathfrak{a}}$,
$\mu=g_{\mid^{u}\!\mathfrak{z}}$, $\mathfrak{h}=\mathfrak{g}(a)$,
$\mathbold{H}=\mathbold{G}(a)$ et $\mathfrak{q}=\ker
a\cap\mathfrak{a}(a)$. Alors, $\mathfrak{q}$ est un idéal
$\mathbold{H}$-invariant de $^{u}\!\mathfrak{h}$. On désigne par
$\mathbold{Q}$ le sous-groupe unipotent de $\mathbold{G}$ d'algèbre de
Lie $\mathfrak{q}$, lequel est contenu dans $\mathbold{H}$. On pose
$\mathbold{G}_{1}=\mathbold{H}/\mathbold{Q}$, qui est un groupe
algébrique sur $\mathbb{K}$ d'algèbre de Lie
$\mathfrak{g}_{1}=\mathfrak{h}/\mathfrak{q}$, et on désigne par
$g_{1}$ la forme linéaire sur $\mathfrak{g}_{1}$ induite par
$h=g_{\mid\mathfrak{h}}$. Enfin, on désigne par $\mathfrak{z}_{1}$
(resp. $\mathfrak{z}_{1,u}$) le centre de $\mathfrak{g}_{1}$
(resp. $^{u}\!\mathfrak{g}_{1}$).

D'après le lemme \ref{co1l1}, on a
$\dim\mathbold{G}.a=\dim(\mathfrak{a}/^{u}\!\mathfrak{z})$, de sorte que
$\dim\mathfrak{g}_{1}<\dim\mathfrak{g}$. D'autre part, d'après \cite[I.16 et
I.18]{duflo-1982}, $g_{1}$ est de type unipotent et l'on a
\begin{equation*}
\mathfrak{g}_{1}(g_{1})=
\mathfrak{h}(h)/\mathfrak{q}=(\mathfrak{g}(g)+\mathfrak{a}(a))/\mathfrak{q}.
\end{equation*}
De plus, $\mathfrak{a}(a)/\mathfrak{q}$ est un idéal central de
$\mathfrak{g}_{1}$ contenu dans $^{u}\!\mathfrak{g}_{1}$. Par suite,
on a~:
\begin{equation*}
^{u}\!\mathfrak{z}_{1}\subset\,^{u}\!(\mathfrak{g}_{1}(g_{1}))=
(^{u}\!(\mathfrak{g}(g))+\mathfrak{a}(a))/\mathfrak{q}=\mathfrak{a}(a)/\mathfrak{q}
\subset\,^{u}\!\mathfrak{z}_{1}.
\end{equation*}
On en déduit que $\mathfrak{g}_{1}$ est quasi-réductive et que $g_{1}$ est
une forme de type unipotent et réductif.

La forme bilinéaire $(X,Y)\mapsto\langle a,[X,Y]\rangle$ induit une
dualité $\mathfrak{h}$-invariante entre $\mathfrak{g}/\mathfrak{h}$ et
$\mathfrak{a}/^{u}\!\mathfrak{z}$ ainsi qu'une forme symplectique
$\mathfrak{h}$-invariante sur $\mathfrak{a}/\mathfrak{a}(a)$. On
déduit de ceci que l'action de $\mathfrak{h}$ dans
$\mathfrak{h}/\mathfrak{a}(a)$ est unimodulaire et, puisque
$\mathfrak{a}(a)/\mathfrak{q}$ est central dans $\mathfrak{g}_{1}$,
que $\mathfrak{g}_{1}$ est unimodulaire. On voit donc que l'hypothèse
de récurrence s'applique à $\mathfrak{g}_{1}$ et $g_{1}$, de sorte
que, posant $u_{1}=(g_{1})_{\mid^{u}\!\mathfrak{g}_{1}}$, on a

(i) $\mathfrak{z}_{1,u}=\,^{u}\!(\mathfrak{g}_{1}(g_{1}))
=\,^{u}\!\mathfrak{g}_{1}(u_{1})=\,^{u}\!\mathfrak{z}_{1}=
\mathfrak{a}(a)/\mathfrak{q}$,

(ii) $\mathfrak{g}_{1}(g_{1})=\mathfrak{g}_{1}(u_{1})$,

(iii) $\mathfrak{g}_{1}=\mathfrak{g}_{1}(u_{1})+\,^{u}\!\mathfrak{g}_{1}$.

\subsection{}\label{2d}
Dans ce numéro, on montre que $^{u}\!\mathfrak{z}=\mathfrak{z}_{u}$~;
supposons au contraire que
$^{u}\!\mathfrak{z}\subsetneqq\mathfrak{z}_{u}$ et soit
$\mathfrak{a}\subset\mathfrak{z}_{u}$ un idéal
$\mathbold{G}$-invariant, contenant strictement $^{u}\!\mathfrak{z}$
et minimal pour cette propriété. On applique à $\mathfrak{a}$ les
résultats du numéro \ref{2c} dont on reprend les notations.

On a $\mathfrak{a}=\mathfrak{a}(a)\subset\mathfrak{h}$. Soit
$\mathbold{R}$ un facteur réductif de $\mathbold{G}$ et $\mathfrak{r}$
son algèbre de Lie. Alors, on a
\begin{eqnarray*}
\mathfrak{h} & = &\mathfrak{r}(a)\oplus\,^{u}\!\mathfrak{g}\\
^{u}\!\mathfrak{h} & = & \,^{u}\!(\mathfrak{r}(a))\oplus\,^{u}\!\mathfrak{g}
\end{eqnarray*}
de sorte que $\mathfrak{r}(a)$ s'identifie via la projection naturelle à
une sous-algèbre de $\mathfrak{g}_{1}$ et que l'on peut alors écrire~:
\begin{eqnarray*}
\mathfrak{g}_{1} & = &
\mathfrak{r}(a)\oplus\,^{u}\!\mathfrak{g}/\mathfrak{q}\\
^{u}\!\mathfrak{g}_{1} & = &
\,^{u}\!(\mathfrak{r}(a))\oplus\,^{u}\!\mathfrak{g}/\mathfrak{q}.
\end{eqnarray*}

Soit $\mathfrak{m}\subset\,^{u}\!\mathfrak{g}$ un supplémentaire
$\mathfrak{r}$-invariant de $\mathfrak{a}$, que l'on identifie via la
projection naturelle à un sous-espace de $\mathfrak{g}_{1}$. Alors, on
a
\begin{equation*}
^{u}\!\mathfrak{g}_{1}=
\,^{u}\!\mathfrak{g}/\mathfrak{q}=\mathfrak{m}\oplus\mathfrak{a}/\mathfrak{q},
\end{equation*}
$\mathfrak{m}$ étant un supplémentaire $\mathfrak{r}(a)$-invariant de
$\mathfrak{a}/\mathfrak{q}=\mathfrak{z}_{1,u}$ dans $^{u}\!\mathfrak{g}_{1}$.

Comme, d'après les résultats du numéro \ref{2c}, on a
$^{u}\!\mathfrak{g}_{1}(u_{1})=\mathfrak{z}_{1,u}$, il vient~:
\begin{equation*}
^{u}\!\mathbold{G}_{1}.u_{1}=u_{1}+\mathfrak{z}_{1,u}^{\perp},
\end{equation*}
si bien, que quitte à translater $g$ par un élément de
$^{u}\!\mathbold{H}$, on peut supposer que $u_{1}$ s'annule sur
$^{u}\!(\mathfrak{r}(a))\oplus\mathfrak{m}$. Mais alors, on a
$^{u}\!(\mathfrak{r}(a))\subset\,^{u}\!\mathfrak{g}_{1}(u_{1})=
\mathfrak{z}_{1,u}$
et donc $^{u}\!(\mathfrak{r}(a))=0$. Autrement dit, $\mathfrak{r}(a)$
est une algèbre de Lie réductive et
\begin{equation*}
^{u}\!\mathfrak{g}_{1}=
\mathfrak{m}\oplus\mathfrak{z}_{u}/\mathfrak{q}=
\,^{u}\!\mathfrak{g}/\mathfrak{q}
\end{equation*}
d'où suit que
\begin{equation*}
^{u}\!\mathfrak{g}  =  \,^{u}\!\mathfrak{h},
\end{equation*}
et aussi
\begin{equation*}
  ^{u}\!\mathfrak{g}(u)/\mathfrak{q} =\!^{u}\!\mathfrak{g}_{1}(u_{1})=
\mathfrak{z}_{1,u}.
\end{equation*}
On voit donc d'une part que, quitte à translater par un élément de
$^{u}\!\mathbold{H}=\,^{u}\!\mathbold{G}$, on peut supposer que
$u_{\mid\mathfrak{m}=0}$ et d'autre part que
$\mathfrak{z}_{u}/\mathfrak{q}$ est contenu dans
$\mathfrak{z}_{1,u}=\mathfrak{a}/\mathfrak{q}$ de sorte que
$\mathfrak{z}_{u}=\mathfrak{a}$.  On déduit de cette dernière
assertion que l'action de $\mathbold{R}$ dans
$\mathfrak{z}_{u}/^{u}\!\mathfrak{z}$ est irréductible. Il suit
également de ce qui précède que
\begin{equation}\label{eq2c1}
^{u}\!\mathfrak{g}(u)=\mathfrak{z}_{u}.
\end{equation}

Par ailleurs, on a le résultat suivant~:
\begin{lem}\label{lem2c1}
Soit $\mathbold{V}$ un espace vectoriel de dimension finie sur
$\mathbb{K}$ et $\mathbold{R}$ un sous-groupe algébrique réductif du
groupe des automorphismes affines de $\mathbold{V}$. On suppose que
$\mathbold{R}$ a une orbite ouverte dans $\mathbold{V}$ et que le
centre de $\mathbold{R}$ est fini. Alors le groupe
d'isotropie d'un point de l'orbite ouverte n'est pas réductif.
\end{lem}
\begin{dem}
On peut supposer que $\mathbold{R}$ est connexe. Comme $\mathbold{R}$
est réductif et agit par automorphismes affines, il possède un point
fixe dans $\mathbold{V}$~: on peut donc supposer que
l'action de $\mathbold{R}$ dans $\mathbold{V}$ est linéaire, cas qui a
été traité par Anh dans le cas réel (voir \cite[Lemma
1]{anh-1978}). Pour la commodité du lecteur nous donnons la
démonstration dans le cas général.

Supposons que le stabilisateur d'un point de l'orbite
ouverte $\mathbold{O}$ soit réductif. Alors, il suit de
\cite{borel-harish-chandra-1962} que $\mathbold{O}$ est une variété
affine. Par suite son complémentaire est l'ensemble des zéros d'une
fonction polynôme $P$ sur $\mathbold{V}$ et il est clair que $P$ est un
semi-invariant pour $\mathbold{R}$ dont le poids est un caractère non
trivial : par suite le centre de $\mathbold{R}$ n'est pas
fini.
\end{dem}

Reprenons la démonstration de notre assertion. Posons
$\mu=g_{\mid^{u}\!\mathfrak{z}}$. Comme $\mathfrak{a}=\mathfrak{z}_{u}$, on a
$\mathbold{G}.a=\mathbold{R}.a$ de sorte que, compte tenu du lemme
\ref{co1l1}, $\mathbold{R}$ possède une orbite ouverte dans l'espace
affine $\mathfrak{a}^{*}_{\mu}$.  Comme le $\mathbold{R}$-module
$\mathfrak{z}_{u}/^{u}\!\mathfrak{z}$ est irréductible, le lemme \ref{lem2c1}
montre alors que le centre de $\mathbold{R}$ agit par un caractère non
trivial $\alpha$ dans $\mathfrak{z}_{u}/^{u}\!\mathfrak{z}$. Il existe donc un
sous-groupe à un paramètre $\gamma$ du centre de $\mathbold{R}$ et un
entier naturel non nul $n_{0}$ tels que
$\alpha\circ\gamma(t)=t^{n_{0}}$, $t\in\mathbb{K}^×$.

Soit $e_{1},\ldots,e_{s}$ une base de $\mathfrak{m}$ constituée de vecteurs
propres pour le sous-groupe à un paramètre $\gamma$~: si $1\leq j\leq
s$, il existe un entier relatif $n_{j}$ tel que
$\Ad\gamma(t).e_{j}=t^{n_{j}}e_{j}$. D'autre part, soit
$\mathfrak{z}_{u,n_{0}}=\{X\in\mathfrak{z}_{u} : \Ad\gamma(t).X=t^{n_{0}}X\}$.
Alors, on a $\mathfrak{z}_{u}=
\mathfrak{z}_{u,n_{0}}\oplus\,^{u}\!\mathfrak{z}$ et
donc
$^{u}\!\mathfrak{g}=
\mathfrak{m}\oplus\mathfrak{z}_{u,n_{0}}\oplus\,^{u}\!\mathfrak{z}$.
Cette décomposition induit la décomposition duale
$^{u}\!\mathfrak{g}^{*}
=\mathfrak{m}^{*}\oplus\mathfrak{z}^{*}_{u,n_{0}}
\oplus\,^{u}\!\mathfrak{z}^{*}$.
Cela dit, on a vu que l'on peut supposer que $u_{\mid\mathfrak{m}}=0$
de sorte que l'on peut écrire $u=v+w$ avec
$v\in\mathfrak{z}^{*}_{u,n_{0}}$ et $w\in\,^{u}\!\mathfrak{z}^{*}$. De
plus, on a vu que
$^{u}\!\mathfrak{g}_{1}(u_{1})=\mathfrak{z}_{1,u}=
\mathfrak{z}_{u}/\mathfrak{q}$,
de sorte que la forme
$(\beta_{u})_{\mid\mathfrak{m}}=(\beta_{u_{1}})_{\mid\mathfrak{m}}$
est une forme symplectique. Si $1\leq i,j\leq s$, on a alors pour
$t\in\mathbb{K}^×$,
\begin{equation*}
\begin{split}
\langle u,[e_{i},e_{j}]\rangle & =
\langle\sideset{}{^{*}}\Ad\gamma(t).u,\Ad\gamma(t)[e_{i},e_{j}]\rangle\\
& = t^{n_{i}+n_{j}-n_{0}}\langle v,[e_{i},e_{j}]\rangle +
t^{n_{i}+n_{j}}\langle w,[e_{i},e_{j}]\rangle,
\end{split}
\end{equation*}
d'où résulte que $n_{i}+n_{j}\in\{0,n_{0}\}$ si
$\beta_{u}(e_{i},e_{j})\neq0$. Or, on a
\begin{equation*}
(\beta_{u})_{\mid\mathfrak{m}}=\sum_{\substack{1\leq i<j\leq
s\\ \beta_{u}(e_{i},e_{j})\neq
0}}\beta_{u}(e_{i},e_{j})e_{i}^{*}\wedge e_{j}^{*}
\end{equation*}
et, par suite
\begin{equation*}
\sideset{}{^{*}}\Ad\gamma(t^{-1})((\beta_{u})_{\mid\mathfrak{m}})=
t^{n_{0}}\sum_{\substack{1\leq i<j\leq
s\\ n_{i}+n_{j}=n_{0}}}\beta_{u}(e_{i},e_{j})e_{i}^{*}\wedge e_{j}^{*}
+\sum_{\substack{1\leq i<j\leq
s\\ n_{i}+n_{j}=0}}\beta_{u}(e_{i},e_{j})e_{i}^{*}\wedge e_{j}^{*}.
\end{equation*}
On en déduit que
\begin{equation*}
\sideset{}{_{\mathfrak{m}}}\det\Ad\gamma(t)=
\sideset{}{_{(e_{1},\ldots,e_{s})}}\Pfaff
\sideset{}{^{*}}\Ad\gamma(t^{-1})((\beta_{u})_{\mid\mathfrak{m}})/
\sideset{}{_{(e_{1},\ldots,e_{s})}}\Pfaff((\beta_{u})_{\mid\mathfrak{m}})
\end{equation*}
est un polynôme en $t$ et donc qu'il existe un entier naturel $n$ tel
que, pour $t\in\mathbb{K}^×$,
\begin{equation*}
\sideset{}{_{\mathfrak{m}}}\det\Ad\gamma(t)=t^{n}.
\end{equation*}
Si $r=\dim\mathfrak{z}_{u,n_{0}}$, on a donc, pour $t\in\mathbb{K}^×$,
\begin{equation*}
\sideset{}{_{^{u}\!\mathfrak{g}}}\det\Ad\gamma(t)=t^{n+rn_{0}}
\end{equation*}
avec $n+rn_{0}>0$ puisque $n_{0}>0$ et
$r=\dim\mathfrak{z}_{u}/^{u}\!\mathfrak{z}>0$. Ceci contredit le fait que
$\mathfrak{g}$ est unimodulaire. Compte tenu de la relation \ref{eq2c1}, on
a montré que
$\mathfrak{z}_{u}=
\,^{u}\!\mathfrak{z}=\,^{u}\!(\mathfrak{g}(g))=\,^{u}\!\mathfrak{g}(g)
=\,^{u}\!\mathfrak{g}(u)$.

\subsection{}\label{2e}
Pour achever la démonstration du théorème, on distingue deux
cas.

Commençons par remarquer que si $\mu=0$, on a
$^{u}\!\mathfrak{g}=\,^{u}\!\mathfrak{z}$. En effet, dans cette
situation, l'algèbre de Lie
$\mathfrak{g}'=\mathfrak{g}/^{u}\!\mathfrak{z}$ et la forme linéaire
$g'$ induite sur $\mathfrak{g}'$ par $g$ satisfont les hypothèses du
théorème~; si l'on désigne par $\mathfrak{z}'$
(resp. $\mathfrak{z}'_{u}$) le centre de $\mathfrak{g}'$
(resp. $^{u}\!\mathfrak{g}'$), il résulte du numéro précédent que l'on
a
$\mathfrak{z}'_{u}=\,^{u}\!\mathfrak{z}'=\,^{u}\!\mathfrak{g}'(g')=0$,
si bien que
$^{u}\!\mathfrak{g}/^{u}\!\mathfrak{z}=\,^{u}\!\mathfrak{g}'=0$.

Supposons donc que $^{u}\!\mathfrak{g}=\,^{u}\!\mathfrak{z}$. Alors,
$\mathfrak{g}$ est le produit direct de son unique facteur réductif
$\mathfrak{r}$ par son radical unipotent $^{u}\!\mathfrak{z}$. Comme
$g$ est de type unipotent, $g_{\mid\mathfrak{r}}=0$ et donc $g=0$. Les
assertions (i), (ii) et (iii) du théorème sont alors claires.

Supposons maintenant
$^{u}\!\mathfrak{z}=\mathfrak{z}_{u}\subsetneqq\,^{u}\!\mathfrak{g}$. Alors
$\mu\neq0$ tandis que l'algèbre de Lie $\mathfrak{g}'=\mathfrak{g}/\ker\mu$ et
la forme linéaire $g'$ induite sur $\mathfrak{g}'$ par $g$ satisfont les
hypothèses du théorème. On peut donc supposer que
$\dim\,^{u}\!\mathfrak{z}=1$.

Remarquons que pour achever de démontrer le théorème, il suffit de
montrer que $\mathfrak{g}(g)$ contient un facteur réductif
$\mathfrak{r}$ de $\mathfrak{g}$. En effet, dans ce cas on a
$\mathfrak{r}\oplus\,^{u}\!\mathfrak{z}=\mathfrak{g}(g)\subset\mathfrak{g}(u)$
et l'assertion (iii) est vérifiée. Cependant, si on applique le lemme
\ref{co1l1} à l'idéal $\mathfrak{a}=\,^{u}\!\mathfrak{g}$, on voit que
$\dim\mathfrak{g}/\mathfrak{g}(u)=\dim\,^{u}\!\mathfrak{g}/^{u}\!\mathfrak{z}$,
d'où résulte que
$\dim\mathfrak{g}(u)=\dim\mathfrak{r}+\dim\,^{u}\!\mathfrak{z}$ et
donc que $\mathfrak{g}(u)=\mathfrak{g}(g)$ et
$^{u}\!(\mathfrak{g}(g))=\,^{u}\!\mathfrak{g}(u)=\,^{u}\mathfrak{z}$,
ce qui montre (i) et (ii). Réciproquement, si les conditions (i) et
(ii) sont satisfaites, tout facteur réductif de $\mathfrak{g}(g)$ est
un facteur réductif de $\mathfrak{g}$.

Montrons donc que $\mathfrak{g}(g)$ contient un facteur réductif de
$\mathfrak{g}$. On procède par récurrence sur la dimension de
$\mathfrak{g}$. On suppose donc $\dim\mathfrak{g}>0$ et le résultat
établi en dimension strictement inférieure. Comme on a vu, on peut
supposer que $\dim\,^{u}\!\mathfrak{z}=1$ et $\mu\neq0$. Soit
$\mathfrak{a}\subset\,^{u}\!\mathfrak{g}$ un idéal
$\mathbold{G}$-invariant contenant $^{u}\mathfrak{z}$ et tel que
l'idéal $\mathfrak{a}/^{u}\!\mathfrak{z}$ soit non trivial et central
dans $\mathfrak{g}/^{u}\!\mathfrak{z}$. Alors, $\mathfrak{a}$ contient
strictement $^{u}\!\mathfrak{z}=\mathfrak{z}_{u}$ et, avec les
notations du numéro \ref{2c}, le couple $(\mathfrak{g}_{1},g_{1})$
vérifie les assertions (i) à (iii) du théorème.

Commençons maintenant par montrer que
$^{u}\mathbold{G}.a=a+(\mathfrak{a}/^{u}\mathfrak{z})^{*}$, où
$(\mathfrak{a}/^{u}\mathfrak{z})^{*}$ est identifié à l'orthogonal de
$^{u}\mathfrak{z}$ dans $\mathfrak{a}^{*}$. Comme
$[\,^{u}\mathfrak{g},\mathfrak{a}]\subset\,^{u}\mathfrak{z}$, il est
clair que pour $X\in\,^{u}\!\mathfrak{g}$ et $Y\in\mathfrak{a}$, on a
\begin{equation*}
\begin{split}
\langle\sideset{}{^{*}}\Ad(\exp X).a,Y\rangle & = \langle
a,\Ad(\exp-X).Y\rangle\\
& = \langle a,Y-[X,Y]\rangle\\
& = \langle a+\sideset{}{^{*}}\ad X.a,Y\rangle
\end{split}
\end{equation*}
si bien que $^{u}\mathbold{G}.a=a+\,^{u}\!\mathfrak{g}.a$ est un
sous-espace affine de $\mathfrak{a}^{*}$. De plus, comme
$\mathfrak{z}_{u}$ est le centre de $^{u}\!\mathfrak{g}$, il est clair
que $^{u}\mathbold{G}.a\subset
a+(\mathfrak{a}/^{u}\mathfrak{z})^{*}$. Supposons que l'inclusion soit
stricte. Alors, il existe $Y\in\mathfrak{a}\backslash\mathfrak{z}_{u}$ tel que,
pour tout $X\in\,^{u}\!\mathfrak{g}$, $\langle
a,[X,Y]\rangle=0$. Mais, si $X\in\,^{u}\!\mathfrak{g}$,
$[X,Y]\in\mathfrak{z}_{u}=\,^{u}\!\mathfrak{z}$, tandis que
$\dim\mathfrak{z}_{u}=1$ et $\mu=a_{\mid\mathfrak{z}_{u}}\neq0$~: on
en déduit que pour tout $X\in\,^{u}\!\mathfrak{g}$, $[X,Y]=0$ et donc
que $Y\in\mathfrak{z}_{u}=\,^{u}\!\mathfrak{z}$, ce qui est
absurde. On a donc bien
$^{u}\mathbold{G}.a=a+(\mathfrak{a}/^{u}\mathfrak{z})^{*}$.

Soit $\mathfrak{r}$ un facteur réductif de $\mathfrak{g}$ et soit
$\mathfrak{m}\subset\mathfrak{a}$ un supplémentaire
$\mathfrak{r}$-invariant de $\mathfrak{z}_{u}$. D'après ce qui
précède, on peut supposer que $a_{\mid\mathfrak{m}}=0$. Mais alors, on
a $\mathfrak{r}\subset\mathfrak{g}(a)=\mathfrak{h}$ et $\mathfrak{r}$,
identifié à son image par la projection naturelle de $\mathfrak{h}$
sur $\mathfrak{g}_{1}$, est un facteur réductif de
$\mathfrak{g}_{1}$. Comme $(\mathfrak{g}_{1},g_{1})$ vérifie les
hypothèses du théorème, il résulte de remarques antérieures que
$\mathfrak{g}_{1}(g_{1})$ contient un facteur réductif de
$\mathfrak{g}_{1}$~: il existe donc
$x\in\mathbold{G}_{1}=\mathbold{H}/\mathbold{Q}$ tel que $\Ad
x.\mathfrak{r}\subset\mathfrak{g}_{1}(g_{1})$. Maintenant, on a
$\mathfrak{g}_{1}(g_{1})=\mathfrak{h}(h)/\mathfrak{q}$ et aussi
$\mathfrak{h}(h)=\mathfrak{g}(g)+\mathfrak{a}(a)$. Il en résulte que,
si $y\in\mathbold{H}$ a pour image $x$ par la projection naturelle de
$\mathbold{H}$ sur $\mathbold{G}_{1}$, alors $\Ad
y.\mathfrak{r}\subset\mathfrak{h}(h)$ et donc que $\mathfrak{g}(g)$
contient un facteur réductif de $\mathfrak{g}$.
\end{dem}

\begin{co}\label{co2e1}
Soit $\mathfrak{g}$ une algèbre de Lie algébrique, unimodulaire et
préhomogène. Alors, $\mathfrak{g}$ est le produit direct d'un tore et d'une
algèbre de Lie unipotente et préhomogène.
\end{co}
\begin{dem}
Soit $\mathfrak{j}$ le tore facteur réductif du centre de
$\mathfrak{g}$. D'après le théorème \ref{theo2b1} (ii) et (iii), on a
$\mathfrak{g}=\mathfrak{j}\oplus\,^{u}\!\mathfrak{g}$.
\end{dem}

\subsection{}\label{2f}
Soit $\g$ une algèbre de Lie unimodulaire et quasi-réductive. Il résulte du
théorème \ref{theo2b1} que son radical unipotent est une algèbre de
Lie préhomogène. Dans ce cas la description des formes
quasi-réductives pour $\g$ se ramène à la description de celles pour
son radical unipotent et, dans le cas connexe au moins, les sous-groupes
réductifs canoniques sont tous conjugués~:

\begin{theo}\label{theo2f1}
Soit $\g$ une algèbre de Lie unimodulaire et réductif.

(i) On a
$^{u}\!\mathfrak{z}^{*}_{\g,red}=\,^{u}\!\mathfrak{z}^{*}_{^{u}\!\g,red}$. Plus
précisément, soit $\mu\in^{u}\!\mathfrak{z}^{*}_{^{u}\!\g,red}$,
$\mathfrak{r}$ un facteur réductif de $\g$ et $\mathfrak{m}$ un
supplémentaire $\mathfrak{r}$-invariant de $^{u}\!\mathfrak{z}$ dans
$^{u}\!\mathfrak{g}$. Alors, la forme linéaire $g_{\mu}$ sur $\g$ qui
prolonge $\mu$ et s'annule sur $\mathfrak{r}\oplus\mathfrak{m}$ est de
type réductif et unipotent.

(ii) Soit $\mathbold{G}$ un groupe algébrique d'algèbre de Lie $\g$
tel que $^{u}\!\mathbold{Z}$ soit un sous-groupe central. Alors, les
sous-groupes réductifs canoniques de $\mathbold{G}$ sont ses facteurs
réductifs. En particulier, ils sont tous conjugués sous l'action de
$^{u}\!\mathbold{G}$.
\end{theo}
\begin{dem}
(i) L'inclusion
  $^{u}\!\mathfrak{z}^{*}_{\g,red}\subset\,^{u}\!\mathfrak{z}^{*}_{^{u}\!\g,red}$
  résulte du théorème \ref{theo2b1} (i).  Réciproquement, soit
  $\mu\in\,^{u}\!\mathfrak{z}^{*}_{^{u}\!\g,red}$, $\mathfrak{r}$ un
  facteur réductif de $\g$, $\mathfrak{m}$ un supplémentaire
  $\mathfrak{r}$-invariant de $^{u}\!\mathfrak{z}$ dans $^{u}\!\mathfrak{g}$ et
  $g_{\mu}\in\g^{*}_{\mu}$ telle que
  $g_{\vert\mathfrak{r}\oplus\mathfrak{m}}=0$. Soit $u$ la restriction
  de $g_{\mu}$ à $^{u}\!\mathfrak{g}$. Comme $\mu$ est quasi-réductive
  relativement à $^{u}\!\mathfrak{g}$, on a
  $^{u}\!\mathfrak{g}(u)=\,^{u}\!\mathfrak{z}$. Par construction il est
  clair que $g_{\mu}$ est stabilisé par $\mathfrak{r}$. On en déduit
  immédiatement que
  $\g(g_{\mu})=\mathfrak{r}\oplus\,^{u}\!\mathfrak{z}$. D'où
  l'assertion.

(ii) Soit $\mathbold{R}$ un facteur réductif de $\mathbold{G}$ et
  $\mu\in\,^{u}\!\mathfrak{z}^{*}_{\g,red}$. Si l'on choisit le
  supplémentaire $\mathfrak{m}$ de $^{u}\!\mathfrak{z}$ dans
  $^{u}\!\mathfrak{g}$ du (i), $\mathbold{R}$-invariant, ce qui est
  loisible, il est immédiat que
  $\mathbold{G}(g_{\mu})=\mathbold{R}\,^{u}\!\mathbold{Z}$. D'où
  l'assertion.
\end{dem}

L'exemple suivant montre que l'hypothèse selon laquelle
$^{u}\!\mathbold{Z}$ est un sous-groupe central de $\mathbold{G}$ est
nécessaire, à la fois pour l'assertion du théorème précédent et pour
l'assertion (iii) du théorème \ref{theo1f1}. Le groupe
$\mathbold{G}_{0}$ est égal à $\mathbold{K}\times\mathbold{K}$ où
$\mathbold{K}$ est le produit semi-direct
$\mathrm{Sp}(V)\mathbold{H}_{V}$ et $\mathbold{G}$ est le produit
semi-direct $\{1,\sigma\}\mathbold{G}_{0}$, où $\sigma$ est
l'automorphisme de permutation. L'algèbre de Lie de $\mathbold{G}$ est
$\g=(\mathfrak{sp}(V)\oplus\mathfrak{h}_{V})
\times(\mathfrak{sp}(V)\oplus\mathfrak{h}_{V})$ dont le centre a pour
base les vecteurs $z_{1}=(z,0)$ et $z_{2}=(0,z)$.  On vérifie qu'il y a
deux classes de conjugaison des facteurs réductifs canoniques de
$\mathbold{G}$~: celle de $\mathrm{Sp}(V)\times\mathrm{Sp}(V)$, qui
correspond aux formes de type réductif et unipotent $g$ vérifiant
$g(z_{1})\neq g(z_{2})$ et celle de
$\{1,\sigma\}(\mathrm{Sp}(V)\times\mathrm{Sp}(V))$, qui correspond aux
formes de type réductif et unipotent $g$ vérifiant $g(z_{1})=
g(z_{2})$.

\subsection{}\label{2g}
Le résultat suivant montre que, dans le cas unimodulaire, les orbites
de type réductif et unipotent sont entièrement déterminées, parmi
celles de type unipotent, par leur restriction au radical unipotent du
centre.

\begin{pr}\label{pr2g1}
Soit $\mathfrak{g}$ une algèbre de Lie unimodulaire quasi-réductive de
centre $\mathfrak{z}$ et $g\in\g^{*}$ une forme de type
unipotent. Alors les assertions suivantes sont équivalentes~:

(i) $g$ est de type réductif,

(ii) la restriction de $g$ à $^{u}\!\mathfrak{z}$ est quasi-réductive.
\end{pr}
\begin{dem}
Il suffit de montrer que (ii) $\Rightarrow$ (i). Supposons donc que la
restriction $\mu$ de $g$ à $^{u}\!\mathfrak{z}$ est
quasi-réductive. Par définition, il existe $g'\in\g^{*}_{\mu}$ une
forme linéaire de type réductif et unipotent. Soit $u$ (resp. $u'$) la
restriction de $g$ (resp. $g'$) à $^{u}\!\mathfrak{g}$ et soit
$\mathbold{G}$ un groupe algébrique connexe d'algèbre de Lie $\g$. Il
suit du théorème \ref{theo2b1} que
$^{u}\!\mathbold{G}.u=u+\,^{u}\!\mathfrak{z}^{\perp}$ et donc que
$u'\in\,^{u}\!\mathbold{G}.u$. Ainsi quitte à remplacer $g'$ par un
élément de sa $^{u}\!\mathbold{G}$-orbite, on peut supposer que
$u'=u$. Mais alors, on a
$\mathfrak{h}=\mathfrak{g}(u)=\mathfrak{r}\oplus\,^{u}\!\mathfrak{z}$
avec $\mathfrak{r}$ un facteur réductif de $\mathfrak{g}$ contenu dans
$\mathfrak{g}(g')$, de sorte que $g'_{\mid\mathfrak{r}}=0$. De plus,
$h=g_{\mid\mathfrak{h}}$ est de type unipotent et donc
$h_{\mid\mathfrak{r}}=0$ (voir \cite[I.19]{duflo-1982}). On en déduit
que $g'_{\mid\mathfrak{r}}=0=g_{\mid\mathfrak{r}}$ et donc que $g'=g$.
\end{dem}

\subsection{}\label{2h}
Il est connu depuis longtemps que les orbites coadjointes des algèbres
de Lie unimodulaires sont génériquement fermées (voir
\cite{charbonnel-1982}). Pour les algèbres de Lie quasi-réductives, on
a le résultat plus précis suivant~:
\begin{theo}\label{theo2h1}
Soit $\g$ une algèbre de Lie quasi-réductive et $\mathbold{G}$ un
groupe algébrique d'algèbre de Lie $\g$. Alors, les assertions
suivantes sont équivalentes~:

(i) $\g$ est unimodulaire,

(ii) il existe une orbite sous l'action du groupe adjoint connexe de
type réductif et fermée,

(iii) toute $\mathbold{G}$-orbite de type réductif est fermée.
\end{theo}
\begin{dem}
(i) $\Rightarrow$ (iii) Selon la démonstration de
  \cite[théorème 3.2]{charbonnel-dixmier-1981}, on peut supposer que
  $\mathbb{K}=\mathbb{C}$. Soit $g\in\g^{*}_{red}$ et soit
  $\mathbold{R}$ un facteur réductif de $\mathbold{G}$ contenant
  $\mathbold{R}_{g}$. Soit $\mathbold{J}$ un tore maximal de
  $(\mathbold{R}_{g})_{0}$. C'est un sous-groupe de Cartan de
  $\mathbold{R}_{0}$. Il existe donc une décomposition d'Iwasawa
  $\mathbold{R}=KA\mathbold{N}$ où $K$ (resp. $\mathbold{N}$) est un
  sous-groupe compact (resp. unipotent) maximal de $\mathbold{R}$ et
  $A$ est un sous-groupe déployé maximal de $\mathbold{J}$. Alors,
  $\mathbold{N}\,^{u}\!\mathbold{G}$ est un sous-groupe unipotent de
  $\mathbold{G}$, normalisé par $\mathbold{J}$, et l'on a
  $\mathbold{G}.g=K\mathbold{N}\,^{u}\!\mathbold{G}.g$. Il est alors
  clair que $\mathbold{G}.g$ est fermé.

(iii) $\Rightarrow$ (ii) est évident.

(ii) $\Rightarrow$ (i) Supposons que $\g$ admette une orbite $O$ sous
  l'action du groupe adjoint connexe de type réductif et fermée. Soit
  $\mathfrak{j}$ une sous-algèbre de Cartan-Duflo de $\g$, $\h$ (resp.
  $\mathbold{H}$) son centralisateur dans $\g$
  (resp. $\mathbold{G}$). Alors, d'après le théorème \ref{theo1j1},
  $\h$ est une algèbre de Lie préhomogène de centre
  $\mathfrak{z}_{\h}=\mathfrak{j}\oplus\,^{u}\!\mathfrak{z}$ et
  $O\cap\h^{*}$ est une réunion finie d'orbites de type réductif sous
  l'action du groupe adjoint connexe de $\h$ ou de
  $\mathbold{H}_{0}$~: ces orbites sont donc fermées. Soit $\omega$
  l'une d'entre elles. Il suit du théorème \ref{theo1h1} (ii) qu'il
  existe $\mu\in\mathfrak{z}^{*}_{\h}$ tel que
  $\omega=\h^{*}_{\mu}$. Soit alors $\mathbold{R}$ un facteur réductif
  de $\mathbold{H}_{0}$ d'algèbre de Lie $\mathfrak{r}$. Il agit donc
  par affinités dans l'espace affine $\h^{*}_{\mu}$ et il y admet donc
  un point fixe, disons $g$. On a donc
  $\h(g)=\mathfrak{r}\oplus\,^{u}\!\mathfrak{z}$. Il est donc clair
  que l'action adjointe de $\mathbold{R}$ dans $\h(g)$ est
  unimodulaire. Par ailleurs $\mathbold{R}$ laisse invariant la forme
  symplectique induite par $\beta_{g}$ sur le quotient $\h/\h(g)$. On
  voit donc que son action adjointe dans $\h$ est unimodulaire et donc
  que $\h$ est unimodulaire. Mais alors, il suit du corollaire
  \ref{co2e1} que $\h=\mathfrak{j}\oplus\,^{u}\!\mathfrak{h}$. Comme
  $\h$ est le centralisateur de $\mathfrak{j}$ dans $\g$, on voit que
  $\mathfrak{j}$ est une sous-algèbre de Cartan d'un facteur réductif
  de $\g$. D'autre part, il existe $g\in\h^{*}$ qui soit une forme
  fortement régulière relativement à $\g$. Alors la restriction de
  $\beta_{g}$ à $[\mathfrak{j},\mathfrak{g}]$ est une forme
  symplectique $\mathfrak{j}$-invariante, de sorte que l'action
  adjointe de $\mathfrak{j}$ dans $[\mathfrak{j},\mathfrak{g}]$ est
  unimodulaire. Comme $\g=\h\oplus[\mathfrak{j},\mathfrak{g}]$ et
  $\mathfrak{j}$ est central dans $\h$, on voit que l'action adjointe
  de $\mathfrak{j}$ dans $\g$ est unimodulaire. Par suite, $\g$ est
  unimodulaire.
\end{dem}

\section{Exemple des sous-algèbres paraboliques de $\mathfrak{so}(n,\mathbb{K})$
  et de certaines sous-algè\-bres de Lie de
  $\mathfrak{sp}(2n,\mathbb{K})$}\label{3}

Une question qui se pose assez vite dans la théorie, est la
détermination parmi les sous-algèbres paraboliques des algèbres de Lie
simples, de celles qui sont quasi-réductives. Dans
\cite{panyushev-2005}, Panyushev étudie l'indice des sous-algèbres de
Lie \flqq seaweed\frqq, i.e. qui sont intersection de deux
sous-algèbres paraboliques, des algèbres de Lie classiques. Ce
faisant, il démontre en particulier que les sous-algèbres seaweed de
$\mathfrak{gl}_{n}$ et paraboliques de $\mathfrak{sp}(2n,\mathbb{K})$ sont
quasi-réductives ({\it loc. cit.} théorèmes 4.3 et 5.2), de même
qu'une classe particulière de sous-algèbres paraboliques de
$\mathfrak{so}(n,\mathbb{K})$ ({\it loc. cit.}  théorème 5.4). Dans
cette section, nous allons donner une expression pour l'indice des
sous-algèbres paraboliques de $\mathfrak{so}(n,\mathbb{K})$ et
caractériser celles parmi elles qui sont quasi-réductives. Nous
verrons en particulier que, contrairement à celles des autres algèbres
de Lie classiques, elles ne sont pas toutes quasi-réductives~: ce fait
était déjà une conséquence de \cite[numéro 3.2]{tauvel-yu-2004-a}.

Comme on le verra, cette dernière question se ramène à déterminer
parmi les sous-algèbres de $\mathfrak{sp}(2n,\mathbb{K})$ qui stabilisent un
drapeau en position générique, celles qui sont quasi-réductives. Dans
\cite{dvorsky-2003} Dvorsky introduit ces algèbres de
Lie particulières et montre que le calcul de l'indice de certaines
sous-algèbres seaweed de $\mathfrak{so}(n,\mathbb{K})$ se ramène au calcul de
leur indice.

\subsection{}\label{3a}
Nous commençons par un résultat général qui nous sera utile pour la
suite.
\begin{defi}
Soit $\mathfrak{g}$ une algèbre de Lie algébrique sur $\mathbb{K}$. On
appelle \emph{rang} de $\mathfrak{g}$ sur $\mathbb{K}$ et on note
$\rang\mathfrak{g}$ la dimension commune de ses sous-algèbres de
Cartan-Duflo.
\end{defi}

Soit $\mathfrak{g}$ une algèbre de Lie algébrique sur $\mathbb{K}$,
$\mathfrak{n}$ un idéal de $\mathfrak{g}$ contenu dans
$^{u}\!\mathfrak{g}$ ou, de manière équivalente, un idéal unipotent de
$\mathfrak{g}$, $g\in\mathfrak{g}^{*}$ une forme linéaire fortement
régulière, $n$ la restriction de $g$ à $\mathfrak{n}$, $\mathfrak{h}$
le stabilisateur de $n$ dans $\mathfrak{g}$, $h$ la restriction de $g$
à $\mathfrak{h}$. Alors le noyau $\mathfrak{q}$ de
$n_{\mid\mathfrak{n}(n)}$, est un idéal unipotent de $\mathfrak{h}$.
On pose $\mathfrak{g}_{1}=\mathfrak{h}/\mathfrak{q}$ et on désigne par
$g_{1}$ la forme linéaire sur $\mathfrak{g}_{1}$ induite par $h$.

\begin{lem}\label{lem3a1}
On garde les notations précédentes. Alors, on a~:

(i) $\ind\mathfrak{g}=\ind\mathfrak{g}_{1}+
\dim\mathfrak{q}-\dim\mathfrak{g}/(\mathfrak{h}+\mathfrak{n})$.

(ii) Toute sous-algèbre de Cartan-Duflo de $\mathfrak{g}$ est
conjuguée à un tore maximal de $\mathfrak{h}(h)$ et a même dimension
qu'une sous-algèbre de Cartan-Duflo de $\mathfrak{g}_{1}$. En
particulier, on a $\rang\mathfrak{g}=\rang\mathfrak{g}_{1}$.
\end{lem}
\begin{dem}
D'après \cite[I.16]{duflo-1982}, on a
$\mathfrak{h}(h)=\mathfrak{g}(g)+\mathfrak{n}(n)$ et
$\mathfrak{n}(n).g=(\mathfrak{h}+\mathfrak{n})^{\perp}$. Or $g_{1}$ est
une forme fortement régulière sur $\mathfrak{g}_{1}$ et
$\mathfrak{g}_{1}(g_{1})=\mathfrak{h}(h)/\mathfrak{q}$. Le lemme est
alors clair.
\end{dem}

\subsection{}\label{3a'} Soit $V$ un espace vectoriel non nul de dimension
finie sur $\mathbb{K}$ et soit $\cur{V}=\{\{0\}=V_{0}\subsetneq
V_{1}\subsetneq\cdots\subsetneq V_{t-1}\subsetneq V_{t}=V\}$ un
drapeau de $V$. On désigne par $\mathfrak{q}_{\cur{V}}$ la
sous-algèbre parabolique de $\mathfrak{gl}(V)$ constituée des
endomorphismes qui laissent invariant $\cur{V}$ et par
$\mathbold{Q}_{\cur{V}}$ le sous-groupe algébrique connexe de
$\mathrm{GL}(V)$ d'algèbre de Lie $\mathfrak{q}_{\cur{V}}$.

\begin{defi}\label{def3a'1}
Soit $\cur{V}=\{\{0\}=V_{0}\subsetneq V_{1}\subsetneq\cdots\subsetneq
V_{t-1}\subsetneq V_{t}=V\}$ un drapeau. On note $h=h(\cur{V})$ le
nombre d'indices $i$ tels que $1\leq i\leq t-1$ et $\dim V_{i}$ et
$\dim V_{i+1}$ soient tous deux impairs.  On dit que le drapeau
$\cur{V}$ vérifie la propriété $\cur{P}$, si deux espaces consécutifs
de la suite $\cur{V}$ ne peuvent être tous deux de dimension impaire,
c'est à dire si $h(\cur{V})=0$.
\end{defi}

\begin{defi}\label{def3a'2}
Soit $\xi$ une forme bilinéaire alternée sur $V$. On dit que le
drapeau $\cur{V}$ est \emph{générique relativement à $\xi$} ou que la
forme $\xi$ est\emph{ générique relativement à $\cur{V}$}, si pour
$1\leq i\leq t$, la restriction de $\xi$ à $V_{i}$ est de rang
maximum, savoir $2[\frac{1}{2}\dim V_{i}]$, et
$V_{i}^{\perp_{\xi}}\cap V_{i-1}=0$.
\end{defi}

\begin{lem}\label{lem3a'1}
Soit $V$ un espace vectoriel de dimension finie sur $\mathbb{K}$ et
$\cur{V}=\{\{0\}=V_{0}\subsetneq V_{1}\subsetneq\cdots\subsetneq
V_{t-1}\subsetneq V_{t}=V\}$ un drapeau de $V$. Alors

(i) une forme $\xi\in\wedge^{2}V^{*}$ est générique relativement à
$\cur{V}$, si et seulement s'il existe une base $e_{1},\ldots,e_{r}$
de $V$ adaptée au drapeau $\cur{V}$ telle que, notant
$e_{1}^{*},\ldots,e_{r}^{*}$ la base duale, on ait $\xi=\sum_{1\leq
2i\leq r}e_{2i-1}^{*}\wedge e_{2i}^{*}$,

(ii) l'ensemble des formes bilinéaires alternées $\xi$ sur $V$ qui
sont génériques relativement à $\cur{V}$ est une orbite sous l'action
de $\mathbold{Q}_{\cur{V}}$, ouverte dans $\wedge^{2}V^{*}$.
\end{lem}
\begin{dem}
(i) Seule l'implication directe mérite démonstration. Pour ce faire,
on procède par récurrence sur l'entier $t$. Lorsque $t=1$ le résultat
est évident. Supposons donc $t>1$ et le résultat vrai à l'ordre
$t-1$. Soit $\xi\in\wedge^{2}V^{*}$ une forme générique relativement à
$\cur{V}$. Comme $\xi_{\mid V_{t-1}}$ est évidemment générique
relativement au drapeau $\cur{V}'=\{\{0\}=V_{0}\subsetneq
V_{1}\subsetneq\cdots\subsetneq V_{t-1}\}$ de $V_{t-1}$, on dispose
d'une base $e_{1},\ldots,e_{r_{t-1}}$ de $V_{t-1}$ adaptée à ce
drapeau telle que $\xi_{\mid V_{t-1}}=\sum_{1\leq 2i\leq
r_{t-1}}e_{2i-1}^{*}\wedge e_{2i}^{*}$.

Si $r_{t-1}$ est pair, $\xi_{\mid V_{t-1}}$ est symplectique de sorte
que $V=V_{t-1}\oplus V_{t-1}^{\perp_{\xi}}$ et il suffit de compléter
$e_{1},\ldots,e_{r_{t-1}}$ par une base $e_{r_{t-1}+1},\ldots,e_{r}$
de $V_{t-1}^{\perp_{\xi}}$ telle que $\xi_{\mid
V_{t-1}^{\perp_{\xi}}}=\sum_{r_{t-1}+1\leq 2i\leq r}e_{2i-1}^{*}\wedge
e_{2i}^{*}$.

Si $r_{t-1}$ est impair, $e_{r_{t-1}}$ est contenu dans le noyau de
$\xi_{\mid V_{t-1}}$. Il existe alors un vecteur $e_{r_{t-1}+1}$ de
$V$ tel que $\xi(e_{r_{t-1}+1},e_{i})=0$, $1\leq i\leq r_{t-1}-1$ et
$\xi(e_{r_{t-1}+1},e_{r_{t-1}})=1$. Alors remplaçant dans le drapeau
$\cur{V}$ le sous-espace $V_{t-1}$ par
$V_{t-1}\oplus\mathbb{K}e_{r_{t-1}+1}$, on se ramène au cas précèdent.

(ii) Désignons par $\cur{O}$ l'ensemble des formes bilinéaires
alternées sur $V$ qui sont génériques relativement à $\cur{V}$. Comme
l'ensemble des bases de $V$ adaptées au drapeau $\cur{V}$ sont dans la
même $\mathbold{Q}_{\cur{V}}$-orbite, il résulte de (i) que $\cur{O}$
est lui-même une $\mathbold{Q}_{\cur{V}}$-orbite. Il reste à montrer
que c'est un ouvert de Zariski de $\wedge^{2}V^{*}$.

Posons $r_{i}=\dim V_{i}$, $0\leq i\leq t$. Soit $e_{1},\ldots,e_{r}$
une base de $V$ adaptée au drapeau $\cur{V}$ et
$e_{1}^{*},\ldots,e_{r}^{*}$ la base duale.  Alors, le fait que la
forme bilinéaire alternée $\xi$ soit un élément de $\cur{O}$ se
traduit par les conditions ouvertes suivantes~: pour $1\leq i\leq t$,
les matrices $(\xi(e_{k},e_{l}))_{1\leq k,l\leq r_{i}}$ et
$(\xi(e_{k},e_{l}))_{1\leq k\leq r_{i},1\leq l\leq r_{i-1}}$ sont
respectivement de rang $2[\frac{r_{i}}{2}]$ et $r_{i-1}$. D'où le
lemme.
\end{dem}

\subsection{Indice et quasi-réductivité des algèbres de Lie
$\mathfrak{r}_{\cur{V}}$}\label{3b}

On suppose que $V$ est un espace vectoriel de dimension finie sur
$\mathbb{K}$ muni d'une forme bilinéaire alternée $\xi$ de rang
maximal et on désigne par $\mathfrak{gl}(V)(\xi)$ son stabilisateur dans
$\mathfrak{gl}(V)$. Lorsque $\xi$ est symplectique,
$\mathfrak{gl}(V)(\xi)$ est notée $\mathfrak{sp}(V,\xi)$ ou plus
simplement $\mathfrak{sp}(V)$ et c'est l'algèbre de Lie du groupe
symplectique correspondant. On se donne un drapeau $\cur{V}$ de $V$ et
on désigne par $\mathfrak{r}_{\cur{V}}$ la sous-algèbre de Lie de
$\mathfrak{gl}(V)(\xi)$ constituée des endomorphismes stabilisant
$\cur{V}$ qui est aussi la sous-algèbre de Lie
$\mathfrak{q}_{\cur{V}}(\xi)$ stabilisant $\xi$ de l'algèbre
parabolique $\mathfrak{q}_{\cur{V}}$.

Dans \cite{dvorsky-2003}, Dvorsky a déterminé l'indice des algèbres de
Lie $\mathfrak{r}_{\cur{V}}$ avec $\xi$ symplectique et $\cur{V}$ un
drapeau générique. Nous retrouvons le résultat de
Dvorsky et caractérisons celles de ces algèbres de Lie qui sont
quasi-réductives.

\begin{pr}\label{pr3b1}
Soit $V$ un espace de dimension finie sur $\mathbb{K}$ muni d'une
forme bilinéaire alternée de rang maximal $\xi$ et
$\cur{V}=\{\{0\}=V_{0}\subsetneq V_{1}\subsetneq\cdots\subsetneq
V_{t-1}\subsetneq V_{t}=V\}$ un drapeau générique relativement à
$\xi$.

a) Supposons que $\dim V$ est pair, de sorte que $\xi$ est
symplectique. Alors,

(i) on a $\ind\mathfrak{r}_{\cur{V}}=\sum_{i=1}^{t}[\frac{1}{2}(\dim
V_{i}-\dim V_{i-1})]$,

(ii) la dimension du radical unipotent d'un stabilisateur générique de
la représentation coadjointe de $\mathfrak{r}_{\cur{V}}$ est $h(\cur{V})$,

(iii) l'algèbre de Lie $\mathfrak{r}_{\cur{V}}$ est quasi-réductive si
et seulement si le drapeau $\cur{V}$ vérifie la propriété $\cur{P}$.

b) Supposons que $\dim V$ est impair et désignons par $\cur{V}'$ le
drapeau du sous-espace $V_{t-1}$ obtenu en supprimant l'espace $V$ du
drapeau $\cur{V}$. Alors,

(iv) on a \begin{equation*}
\ind\mathfrak{r}_{\cur{V}}=\begin{cases}
\sum_{i=1}^{t}[\frac{1}{2}(\dim
V_{i}-\dim V_{i-1})]-1 & \text{si $\dim V_{t-1}<\dim V-1$}\\
\sum_{i=1}^{t}[\frac{1}{2}(\dim
V_{i}-\dim V_{i-1})]+1 & \text{si $\dim V_{t-1}=\dim V-1$},
\end{cases}
\end{equation*}

(v) la dimension du radical unipotent d'un stabilisateur générique de
la représentation coadjointe de $\mathfrak{r}_{\cur{V}}$ est $h(\cur{V}')$,

(vi) l'algèbre de Lie $\mathfrak{r}_{\cur{V}}$ est quasi-réductive si
et seulement si le drapeau $\cur{V}'$ vérifie la propriété $\cur{P}$.
\end{pr}

\begin{dem}
Supposons qu'il existe $1\leq i\leq t$ tel que $\dim V_{i}$ soit
pair. Alors $\xi$ induit une forme symplectique sur le sous-espace
$U=V_{i}$ de $V$ de sorte que $V$ est la somme directe de $U$ et de
son orthogonal $W$ relativement à $\xi$. De plus si $\cur{U}$
(resp. $\cur{W}$) désigne le drapeau $\{\{0\}=V_{0}\subsetneq
V_{1}\subsetneq\cdots\subsetneq V_{i-1}\subsetneq V_{i}\}$
(resp. $\{\{0\}=W_{0}\subsetneq W_{1}=V_{i+1}\cap
W\subsetneq\cdots\subsetneq W_{j}=V_{i+j}\cap
W\subsetneq\cdots\subsetneq W_{t-i}=W\}$), $\cur{U}$ (resp. $\cur{W}$)
est générique relativement à $\xi_{\mid U}$ (resp. $\xi_{\mid
W}$). Par suite, $\mathfrak{r}_{\cur{V}}$ s'identifie au produit
direct $\mathfrak{r}_{\cur{U}}\times\mathfrak{r}_{\cur{W}}$. Il est
alors clair qu'il suffit de vérifier les assertions de la proposition
pour $\mathfrak{r}_{\cur{U}}$ et $\mathfrak{r}_{\cur{W}}$. Par
ailleurs lorsque $t=1$ et $\dim V$ est pair, $\mathfrak{r}_{\cur{V}}$
est l'algèbre de Lie $\mathfrak{sp}(V)$, tandis que lorsque $\dim V=1$,
$t$ est égal à $1$ et $\mathfrak{r}_{\cur{V}}=\mathfrak{gl}(V)$ est un
tore~: dans ces deux cas, le résultat est évident.

On se ramène donc à démontrer la proposition dans l'un des cas
suivants~:

- $\dim V$ est pair, $t\geq2$ et pour $1\leq i\leq t-1$, $\dim V_{i}$
est un nombre impair,

- $\dim V>1$ est impair et $t=1$,

- $\dim V$ et $\dim V_{t-1}$ sont impairs.

\subsection{}\label{3b'}
Dans ce numéro, on se place dans le premier cas. On pose $\dim
V_{i}=2p_{i}+1$, $1\leq i\leq t-1$, $\dim
V=2p_{t}=2p$ et $a_{i}=\dim V_{i}-\dim V_{i-1}$, $1\leq i\leq
t$. Alors $a_{1}=2p_{1}+1=2q_{1}+1$ et
$a_{t}=2(p_{t}-p_{t-1})-1=2q_{t}+1$ sont impairs, tandis que
$a_{i}=2(p_{i}-p_{i-1})=2(q_{i}+1)$, $2\leq i\leq t-1$, est pair. On a
donc $h(\cur{V})=t-2$.

D'après le lemme \ref{lem3a'1}, on peut trouver une base $e_{1},\ldots,
e_{2p}$ de $V$ adaptée au drapeau $\cur{V}$ dans laquelle la matrice
de $\xi$ soit diagonale par blocs de la forme
\begin{equation*}
J_{2p}=
\left[\begin{array}{ccccccccc}
J_{2q_{1}} & & & & & & & &\\
&J& & & & & & & \\
&& J_{2q_{2}}& & & & & & \\
&&&J& & & & &\\
&&&&\ddots& & & &\\
&&&&&J&&&\\
&&&&&&J_{2q_{t-1}}&&\\
&&&&&&&J&\\
&&&&&&&&J_{2q_{t}}
\end{array}
\right]
\end{equation*}
où
$J_{2}=J=
\begin{bmatrix}
0 & 1\\
-1 & 0
\end{bmatrix}$ et, pour $k\in\mathbb{N}^{*}$, $J_{2q}$ est la matrice
diagonale par blocs d'ordre $2q$
\begin{equation*}
J_{2q}=\left[\begin{array}{cccc}
J&&&\\
&J&&\\
&&\ddots&\\
&&&J
\end{array}
\right],
\end{equation*}
$J_{0}$ étant la matrice vide. On désigne par $\mathfrak{sp}_{2q}$
l'algèbre de Lie des matrices $X$ d'ordre $2q$ sur $\mathbb{K}$ telles
que $^{t}XJ_{2q}+J_{2q}X=0$.

Posant $p_{0}=-1$, on introduit les sous-espaces
$W_{i}=\mathbb{K}e_{2p_{i-1}+3}\oplus\mathbb{K}e_{2p_{i-1}+4}\oplus
\cdots\oplus\mathbb{K}e_{2p_{i}}$, $1\leq i \leq t$. Pour $1\leq
i\leq t$, $W_{i}$ est un sous-espace symplectique de dimension $2q_{i}$ de
$V$ et on identifie $\mathfrak{sp}(W_{i})$ à la sous-algèbre de Lie de
$\mathfrak{sp}(V)$ constituée des éléments agissant trivialement sur
l'orthogonal de $W_{i}$.

Si $u,v\in V$, on définit un élément $u\odot v$ de $\mathfrak{sp}(V)$ en
posant
\begin{equation*}
u\odot v(x)=\xi(u,x)v+\xi(v,x)u\mbox{, }x\in V.
\end{equation*}
Alors, l'application $(u,v)\mapsto u\odot v$ induit un isomorphisme
d'espaces vectoriels de $S^{2}(V)$ sur $\mathfrak{sp}(V)$. De plus, si
$u_{1},v_{1},u_{2},v_{2}\in V$, on a
\begin{equation*}
[u_{1}\odot v_{1},u_{2}\odot v_{2}]=\xi(u_{1},u_{2})v_{1}\odot v_{2}+
\xi(u_{1},v_{2})v_{1}\odot u_{2}+\xi(u_{2},v_{1})u_{1}\odot v_{2}+
\xi(v_{1},v_{2})u_{1}\odot u_{2}
\end{equation*}
et si $X\in\mathfrak{sp}(V)$, $u,v\in V$
\begin{equation*}
[X,u\odot v]=(X.u)\odot v+u\odot(X.v).
\end{equation*}

Soit $\mathfrak{l}_{\cur{V}}\subset\mathfrak{sp}(V)$ la sous-algèbre
constituée des éléments dont la matrice dans la base
$e_{1},\ldots,e_{2p}$ est diagonale par blocs de la forme
\begin{equation*}
\left[
\begin{array}{ccccccccccc}
A_{1}&&&&&&&&&&\\
&a_{1}&&&&&&&&&\\
&&-a_{1}&&&&&&&&\\
&&&A_{2}&&&&&&&\\
&&&&a_{2}&&&&&&\\
&&&&&-a_{2}&&&&&\\
&&&&&&\ddots&&&&\\
&&&&&&&A_{t-1}&&&\\
&&&&&&&&a_{t-1}&\\
&&&&&&&&&-a_{t-1}&\\
&&&&&&&&&&A_{t}
\end{array}\right]
\end{equation*}
avec $A_{i}\in\mathfrak{sp}_{2q_{i}}$, $1\leq i\leq t$ et
$a_{i}\in\mathbb{K}$, $1\leq i\leq t-1$. En fait
$\mathfrak{l}_{\cur{V}}$ est l'ensemble des éléments $A$ de
$\mathfrak{sp}(V)$ qui s'écrivent
\begin{equation*}
A=\sum_{i=1}^{t}A_{i}+\sum_{i=1}^{t-1}a_{i}H_{i},
\end{equation*}
avec $H_{i}=-e_{2p_{i}+1}\odot e_{2p_{i}+2}$, $a_{i}\in\mathbb{K}$,
$1\leq i\leq t-1$ et $A_{i}\in\mathfrak{sp}(W_{i})$, $1\leq i\leq t$.

Pour $\alpha=(\alpha_{2},\ldots,\alpha_{t})\in W_{2}\times\cdots\times
W_{t}$, $\beta=(\beta_{1},\ldots,\beta_{t-1})\in
W_{1}\times\cdots\times W_{t-1}$,
$x=(x_{1},\ldots,x_{t-1})\in\mathbb{K}^{t-1}$ et
$y=(y_{2},\ldots,y_{t-1})\in\mathbb{K}^{t-2}$, on pose
\begin{equation*}
X_{(\alpha,\beta,x,y)}=\sum_{i=1}^{t-1}(\beta_{i}+\alpha_{i+1})\odot
e_{2p_{i}+1}+\frac{1}{2}\sum_{i=1}^{t-1}x_{i}e_{2p_{i}+1}\odot
e_{2p_{i}+1} + \sum_{i=2}^{t-1}y_{i}e_{2p_{i-1}+1}\odot
e_{2p_{i}+1}.
\end{equation*}
Un calcul immédiat montre que, pour $(\alpha,\beta,x,y)$ et
$(\alpha',\beta',x',y')\in (W_{2}\times\cdots\times
W_{t})\times(W_{1}\times\cdots\times W_{t-1})\times\mathbb{K}^{t-1}
\times\mathbb{K}^{t-2}$, on a
\begin{equation*}
\begin{split}
[X_{(\alpha,\beta,x,y)},X_{(\alpha',\beta',x',y')}]&=
\sum_{i=1}^{t-1}(\xi(\beta_{i},\beta'_{i})+\xi(\alpha_{i+1},\alpha'_{i+1}))
e_{2p_{i}+1}\odot e_{2p_{i}+1}\\
&\quad+
\sum_{i=1}^{t-2}(\xi(\alpha_{i+1},\beta'_{i+1})+
\xi(\beta_{i+1},\alpha'_{i+1}))
e_{2p_{i}+1}\odot e_{2p_{i+1}+1}.
\end{split}
\end{equation*}

On vérifie aisément que $\mathfrak{l}_{\cur{V}}$ est un facteur
réductif de $\mathfrak{r}_{\cur{V}}$ et que l'application
$(\alpha,\beta,x,y)\mapsto X_{(\alpha,\beta,x,y)}$
(resp. $(x,y)\mapsto X_{(0,0,x,y)}$) induit un isomorphisme d'espaces
vectoriels de $(W_{2}\times\cdots\times
W_{t})\times(W_{1}\times\cdots\times
W_{t-1})\times\mathbb{K}^{t-1}\times\mathbb{K}^{t-2}$
(resp. $\mathbb{K}^{t-1}\times\mathbb{K}^{t-2}$) sur le radical
unipotent $\mathfrak{n}_{\cur{V}}$ de $\mathfrak{r}_{\cur{V}}$
(resp. le centre $\mathfrak{z}_{\cur{V}}$ de
$\mathfrak{n}_{\cur{V}}$). En particulier, $\mathfrak{n}_{\cur{V}}$
est une algèbre de Lie d'indice de nilpotence 2 et dont le centre est
de dimension $2t-3$. On désigne par $\mathbold{R}_{\cur{V}}$
(resp. $\mathbold{L}_{\cur{V}}$, $\mathbold{N}_{\cur{V}}$) le
sous-groupe algébrique connexe de $\mathrm{GL}(V)$ d'algèbre de Lie
$\mathfrak{r}_{\cur{V}}$ (resp. $\mathfrak{l}_{\cur{V}}$,
$\mathfrak{n}_{\cur{V}}$), de sorte que $\mathbold{N}_{\cur{V}}$ est
le radical unipotent et $\mathbold{L}_{\cur{V}}$ un facteur réductif
de $\mathbold{R}_{\cur{V}}$.

On identifie $(W_{2}\times\cdots\times
W_{t})\times(W_{1}\times\cdots\times W_{t-1})$ à un sous-espace, noté
$\mathfrak{m}_{\cur{V}}$, de $\mathfrak{n}_{\cur{V}}$ au moyen de
l'application $(\alpha,\beta)\mapsto X_{(\alpha,\beta,0,0)}$. On pose
$Z_{i}=e_{2p_{i}+1}\odot e_{2p_{i}+1}$, $1\leq i\leq t-1$,
$T_{i}=e_{2p_{i}+1}\odot e_{2p_{i+1}+1}$, $1\leq i\leq t-2$. Alors, on
a $\mathfrak{n}_{\cur{V}}=
\mathfrak{m}_{\cur{V}}\oplus\mathfrak{z}_{\cur{V}}$, tandis que
$Z_{1},\ldots,Z_{t-1},T_{1},\ldots,T_{t-2}$ est une base de
$\mathfrak{z}_{\cur{V}}$ dont on désigne par
$Z_{1}^{*},\ldots,Z_{t-1}^{*},T_{1}^{*},\ldots,T_{t-2}^{*}$ la base
duale. De plus la décomposition en somme directe précédente permet de
voir $\mathfrak{z}_{\cur{V}}^{*}$ comme un sous-espace de
$\mathfrak{n}_{\cur{V}}^{*}$.

Si $w_{i}\in W_{i}$, on désigne par $\xi(w_{i},.)_{\alpha}$
(resp. $\xi(w_{i},.)_{\beta}$) la forme linéaire sur
$\mathfrak{m}_{\cur{V}}$ définie par
$(\alpha,\beta)\mapsto\xi(w_{i},\alpha_{i})$, $2\leq i \leq t$
(resp. $(\alpha,\beta)\mapsto\xi(w_{i},\beta_{i})$, $1\leq i\leq t-1$).

Soit $n=\sum_{i=1}^{t-1}\zeta_{i}Z_{i}^{*}+
\sum_{i=1}^{t-2}\tau_{i}T_{i}^{*}$ un élément de
$\mathfrak{z}_{\cur{V}}^{*}\subset\mathfrak{n}_{\cur{V}}^{*}$. Alors,
on a
\begin{equation}\label{eq3b'1}
\begin{split}
\exp X_{(\alpha,\beta,x,y)}.n &=n-\zeta_{1}\xi(\beta_{1},.)_{\beta}
-\sum_{i=2}^{t-1}\xi(\zeta_{i-1}\alpha_{i}+\tau_{i-1}\beta_{i},.)_{\alpha}\\
&\quad -\sum_{i=2}^{t-1}\xi(\zeta_{i}\beta_{i}+\tau_{i-1}\alpha_{i},.)_{\beta}
-\zeta_{t-1}\xi(\alpha_{t},.)_{\alpha}
\end{split}
\end{equation}
Soit ${\mathfrak{z}^{*}_{\cur{V}}}'$
(resp. ${\mathfrak{n}_{\cur{V}}^{*}}'$) l'ouvert de Zariski de
$\mathfrak{z}_{\cur{V}}^{*}$ (resp. $\mathfrak{n}_{\cur{V}}^{*}$)
complémentaire des zéros de l'élément
$Z_{1}Z_{t-1}\prod_{i=2}^{t-1}(Z_{i}Z_{i-1}-T_{i-1}^{2})$ de
l'algèbre symétrique de $\mathfrak{z}_{\cur{V}}$
(resp. $\mathfrak{n}_{\cur{V}}$).

On déduit de \ref{eq3b'1} que si $n\in\mathfrak{z}{_{\cur{V}}^{*}}'$, on
a $\mathbold{N}_{\cur{V}}.n=n+\mathfrak{z}_{\cur{V}}^{\perp}$. Par
suite $\mathfrak{z}{_{\cur{V}}^{*}}'$ est un ensemble de représentants
des $\mathbold{N}_{\cur{V}}$-orbites dans
$\mathfrak{n}{_{\cur{V}}^{*}}'$.

Soit $g$ une forme fortement régulière sur
$\mathfrak{g}=\mathfrak{r}_{\cur{V}}$. On peut supposer que la
restriction $n$ de $g$ à $\mathfrak{n}=\mathfrak{n}_{\cur{V}}$ est un
élément de $\mathfrak{n}{_{\cur{V}}^{*}}'$ et, quitte à translater $g$
par un élément de $\mathbold{N}_{\cur{V}}$, un élément de
$\mathfrak{z}{_{\cur{V}}^{*}}'$. On peut alors écrire
$n=\sum_{i=1}^{t-1}\zeta_{i}Z_{i}^{*}+
\sum_{i=1}^{t-2}\tau_{i}T_{i}^{*}$.

Mais alors, si $A\in\mathfrak{l}_{\cur{V}}$ s'écrit
$A=\sum_{i=1}^{t}A_{i}+\sum_{i=1}^{t-1}a_{i}H_{i}$ avec
$A_{i}\in\mathfrak{sp}(W_{i})$, $1\leq i\leq t$, on a
\begin{equation*}
A.n=-2\sum_{i=1}^{t-1}a_{i}\zeta_{i}Z_{i}^{*}-
\sum_{i=1}^{t-1}(a_{i-1}+a_{i})\tau_{i-1}T_{i-1}^{*}.
\end{equation*}
On déduit de ceci que
\begin{equation*}
\mathfrak{h}=\mathfrak{g}(n)=
\prod_{i=1}^{t}\mathfrak{sp}(W_{i})\times\mathfrak{z}_{\cur{V}}.
\end{equation*}

Soit $\mathfrak{q}$ le noyau de la restriction de $n$ à
$\mathfrak{z}_{\cur{V}}$ et
$\mathfrak{g}_{1}=\mathfrak{h}/\mathfrak{q}$. Alors on a
$\ind\mathfrak{g}_{1}+\dim\mathfrak{q}=\ind\mathfrak{h}$ tandis que
\begin{equation*}
  \ind\mathfrak{h}=
\sum_{i=1}^{t}\ind\mathfrak{sp}(W_{i})+\dim\mathfrak{z}_{\cur{V}}
=\sum_{i=1}^{t}\ind\mathfrak{sp}(W_{i})+2t-3.
\end{equation*}
Comme
$\dim\mathfrak{r}_{\cur{V}}-\dim(\mathfrak{h}+\mathfrak{n}_{\cur{V}})
=t-1$ et $\ind\mathfrak{sp}(W_{i})=q_{i}$, $1\leq i\leq t$, il résulte
alors du lemme \ref{lem3a1} que
\begin{equation*}
\ind\mathfrak{r}_{\cur{V}}=\sum_{i=1}^{t}q_{i}+t-2
\end{equation*}
et que les sous-algèbres de Cartan-Duflo de $\mathfrak{r}_{\cur{V}}$
sont conjuguées aux sous-algèbres de Cartan-Duflo de
$\prod_{i=1}^{t}\mathfrak{sp}(W_{i})\times\mathfrak{z}_{\cur{V}}$ de
sorte que
\begin{equation*}
\rang\mathfrak{r}_{\cur{V}}=\sum_{i=1}^{t}q_{i}.
\end{equation*}
Il est alors clair que
\begin{equation*}
\begin{split}
\ind\mathfrak{r}_{\cur{V}}&=
\sum_{i=1}^{t}[\frac{1}{2}(\dim V_{i}-\dim V_{i-1})]\\
&=\rang\mathfrak{r}_{\cur{V}}+h(\cur{V}).
\end{split}
\end{equation*}
D'où les assertions (i) et (ii).  Par ailleurs, comme
$\mathfrak{r}_{\cur{V}}$ contient une sous-algèbre de Cartan de
$\mathfrak{sp}(V)$, aucun élément non nul de l'idéal unipotent
$\mathfrak{z}_{\cur{V}}$ n'est central dans
$\mathfrak{r}_{\cur{V}}$. On voit donc que $\mathfrak{r}_{\cur{V}}$
est quasi-réductive si et seulement si $h(\cur{V})=0$. D'où la
proposition dans le cas considéré.

\subsection{}\label{3b''}
Dans ce numéro, on suppose que $\dim V=2p+1$ est impair, $p\geq1$ et
$t=1$. Alors, on a $\mathfrak{r}_{\cur{V}}=\mathfrak{gl}(V)(\xi)$ et
$h(\cur{V})=0$.  Soit $V'\subset V$ un hyperplan supplémentaire du
noyau de $\xi$ si bien que $\xi_{\mid V'}$ est une forme
symplectique. On identifie $\mathfrak{sp}(V')$ (resp. $\mathbb{K}$) à
une sous-algèbre de Lie de $\mathfrak{gl}(V)$ notée $\mathfrak{sp}(V')$
(resp. $\mathfrak{t}_{V'}$) de la manière suivante~:

- si $A\in\mathfrak{sp}(V')$, on l'étend en un élément encore noté $A$
  de $\mathfrak{gl}(V)$ en décidant que $A_{\mid\ker\xi}=0$,

- si $t\in\mathbb{K}$, on le considère comme l'élément semi-simple
  $H_{t}$ de $\mathfrak{gl}(V)$ agissant trivialement dans $V'$ et par
  multiplication par $t$ dans $\ker\xi$.

Alors les algèbres de Lie $\mathfrak{sp}(V')$ et $\mathfrak{t}_{V'}$
sont le commutant l'une de l'autre dans $\mathfrak{gl}(V)$ et
$\mathfrak{l}=\mathfrak{sp}(V')\times\mathfrak{t}_{V'}$ est un
facteur réductif de $\mathfrak{gl}(V)(\xi)$.

Soit $e\in\ker\xi$ un vecteur non nul. Si $w\in V'$, on désigne par
$X_{w}$ l'élément de $gl(V)$ tel que $X_{w}(v)=\xi(w,v)e$. Alors,
l'application $w\mapsto X_{w}$ induit une bijection de $V'$ sur le
radical unipotent $\mathfrak{n}$ de
$\mathfrak{g}=\mathfrak{gl}(V)(\xi)$, lequel est abélien. Une simple
vérification montre alors que l'on a
\begin{equation*}
\mathfrak{gl}(V)(\xi)=(\mathfrak{sp}(V')\times\mathbb{K})
\oplus V',
\end{equation*}
$\mathfrak{sp}(V')$ agissant de manière naturelle dans $V'$ et
$\mathbb{K}$ par multiplication des scalaires.

Si l'on identifie $V'$ avec son dual au moyen de $\xi$, l'action
coadjointe de $\mathfrak{sp}(V')\times\mathfrak{t}_{V'}$ dans $V'$ est
donnée par $\ad^{*}(A,t).v=(A-t)v$,
$(A,t)\in\mathfrak{sp}(V')\times\mathfrak{t}_{V'}$, $v\in V'$. On en
déduit que le groupe adjoint connexe de $\mathfrak{g}$ a une unique
orbite ouverte dans le dual $\mathfrak{n}^{*}=V'$, savoir
$V'-\{0\}$. Ainsi, si $g\in\mathfrak{g}^{*}$ est une forme fortement
régulière, sa restriction $n$ à $\mathfrak{n}$ est un vecteur non nul
de $V'$. Si $\mathfrak{q}$ est le noyau de $n$, on a donc
\begin{align*}
\mathfrak{h}&= \mathfrak{g}(n)=\mathfrak{l}(n)\oplus\mathfrak{n},\\
\mathfrak{g}_{1}&=
\mathfrak{h}/\mathfrak{q}=\mathfrak{l}(n)\times\mathbb{K},
\end{align*}
tandis que $\dim\mathfrak{g}/(\mathfrak{h}+\mathfrak{n})=\dim
V'=\dim\mathfrak{q}+1$.  Alors, appliquant le lemme \ref{lem3a1}, on
voit que $\ind\mathfrak{gl}(V)(\xi)=\ind\mathfrak{l}(n)$ et
$\rang\mathfrak{gl}(V)(\xi)=\rang\mathfrak{l}(n)$.

Cependant, le calcul montre que $\mathfrak{l}(n)$ est isomorphe à
l'algèbre de Lie $\mathfrak{g}$ du numéro \ref{1g4} construite sur un
sous-espace symplectique $V''$ de l'orthogonal de $n$ dans $V'$. Il
suit alors de ce numéro, que $\mathfrak{l}(n)$ est quasi-réductive,
d'indice et de rang tous deux égaux à $\frac{1}{2}\dim V''=p-1$. Ceci
achève la démonstration de la proposition dans le cas considéré.

\subsection{}\label{3b'''}
Dans ce numéro, on suppose que $\dim V_{t}=2p+1$ et $\dim
V_{t-1}=2(p-m)-1$ sont impairs, de sorte que $m$ est un entier
naturel. Alors, il existe une base $e_{1},\ldots,e_{2p+1}$ de $V$ qui
soit adaptée au drapeau $\cur{V}$ et dans laquelle la matrice de $\xi$
soit diagonale par blocs de la forme
\begin{equation*}
\begin{bmatrix}
J_{2(p-m-1)} & 0 & 0 & 0\\
0 & J & 0 & 0\\
0 & 0 & J_{2m} & 0\\
0 & 0 & 0 & 0
\end{bmatrix}
\end{equation*}
Considérons les sous-espaces $V''\subset V'$ et $V'''$ de $V$, avec
$V'=\mathbb{K}e_{1}\oplus\cdots\oplus\mathbb{K}e_{2p}$,
$V''=V_{t-1}\oplus\mathbb{K}e_{2(p-m)}$ et
$V'''=\mathbb{K}e_{2(p-m)+1}\oplus\cdots\oplus\mathbb{K}e_{2p}$. Alors
, $V''$ et $V'$ sont des sous-espaces symplectiques de $V$, de sorte
que l'on a les sommes directes orthogonales $V'=V''\oplus V'''$ et
$V=V''\oplus V''^{\perp_{\xi}}$, la première étant symplectique. Dans
la suite, on identifie les éléments de $\mathfrak{gl}(V)$
(resp. $\mathfrak{gl}(V')$, $\mathfrak{gl}(V'')$, $\mathfrak{gl}(V''')$,
$\mathfrak{gl}(V''^{\perp_{\xi}})$) avec leur matrice dans la base
$e_{1},\ldots,e_{2p+1}$ (resp. $e_{1},\ldots,e_{2p}$,\,
$e_{1},\ldots,e_{2(p-m)}$,\, $e_{2(p-m)+1},\ldots,e_{2p}$,
$e_{2(p-m)+1},\ldots,e_{2p+1}$).
Introduisons le drapeau $\cur{V}''=\{\{0\}=V_{0}\subsetneq
V_{1}\subsetneq\cdots\subsetneq V_{t-1}\subsetneq V''\}$ de
$V''$. Alors les éléments de $\mathfrak{r}_{\cur{V}}$ sont de la forme
\begin{equation*}
X_{A,a,\alpha,y,\beta,z,D,d,\gamma} =
\begin{bmatrix}
A & 0 & J_{2(p-m-1)}\,^{t}\alpha & 0 & 0\\
\alpha & a & y & \beta & 0\\
0 & 0 & -a & 0 & 0\\
0 & 0 & J_{2m}\,^{t}\beta & D & 0\\
0 & 0 & z & \gamma & d
\end{bmatrix}
\end{equation*}
avec  $\alpha\in\mathbb{K}^{2(p-m-1)}$, $\beta,
\gamma\in\mathbb{K}^{2m}$ des vecteurs lignes, $D\in\mathfrak{sp}(V''')$,
$a,y,z,d\in\mathbb{K}$ et
\begin{align*}
Y_{A,a,\alpha,y} & =\begin{bmatrix}
A & 0 & J_{2(p-m-1)}\,^{t}\alpha\\
\alpha & a & y \\
0 & 0 & -a
\end{bmatrix}\in\mathfrak{r}_{\cur{V}''}\\
T_{D,d,\gamma} &=
\begin{bmatrix}
D & 0\\
\gamma & d
\end{bmatrix}\in\mathfrak{gl}(V''^{\perp_{\xi}})(\xi).
\end{align*}
Alors, on identifie
$\mathfrak{r}_{\cur{V}''}\times\mathfrak{sp}(V''')\times\mathbb{K}$ à
une sous-algèbre de Lie de $\mathfrak{r}_{\cur{V}}$ au moyen de
l'application $(Y_{A,a,\alpha,y},D,d)\mapsto
X_{A,a,\alpha,y,0,0,D,d,0}$.

On désigne par $\mathfrak{n}$ le sous-espace de
$\mathfrak{g}=\mathfrak{r}_{\cur{V}}$ constitué des matrices
$Z_{\beta,\gamma,y,z}=X_{0,0,0,y,\beta,z,0,0,\gamma}$, lorsque
$\beta,\gamma$ parcourent $\mathbb{K}^{2m}$ et $y,z$ parcourent
$\mathbb{K}$. On pose $\mathfrak{n}_{1}=\{Z_{\beta,0,0,0} :
\beta\in\mathbb{K}^{2m}\}$, $\mathfrak{n}_{2}=\{Z_{0,\gamma,0,0} :
\gamma\in\mathbb{K}^{2m}\}$, $Z_{1}=Z_{0,0,1,0}$, $Z_{2}=Z_{0,0,0,1}$
et $\mathfrak{z}=\mathbb{K}Z_{1}\oplus\mathbb{K}Z_{2}$. Alors,
$\mathfrak{n}=\mathfrak{n}_{1}\oplus\mathfrak{n}_{2}\oplus\mathfrak{z}$
est un idéal unipotent de $\mathfrak{r}_{\cur{V}}$ dont le centre est
$\mathfrak{z}$. De plus, on a $\mathfrak{r}_{\cur{V}}=
(\mathfrak{r}_{\cur{V}''}\times\mathfrak{sp}(V''')\times\mathbb{K})\oplus
\mathfrak{n}$.

En fait, pour
$\beta,\beta',\gamma,\gamma'\in\mathbb{K}^{2m}$,
$y,y',z,z',d\in\mathbb{K}$, $D\in\mathfrak{sp}(V''')$ et
$Y_{A,a,\alpha,t}\in\mathfrak{r}_{\cur{V}''}$, on a
\begin{align*}
[Z_{\beta,\gamma,y,z},Z_{\beta',\gamma',y',z'}] & =
2\beta J_{2m}\,^{t}\beta' Z_{1}+
(\gamma
J_{2m}\,^{t}\beta'-\gamma'J_{2m}\,^{t}\beta)Z_{2}\\
[(Y_{A,a,\alpha,t},D,d),Z_{\beta,\gamma,y,z}] &=
Z_{a\beta-\beta D,d\gamma-\gamma D,2ay,(a+d)z}.
\end{align*}

Soit, pour $i=1,2$, $Z_{i}^{*}$ la forme linéaire sur $\mathfrak{n}$,
nulle sur $\mathfrak{n}_{1}\oplus\mathfrak{n}_{2}$ et telle que
$\langle Z_{i}^{*},Z_{j}\rangle=\delta_{i,j}$, $j=1,2$. Enfin, pour
$\alpha$ un vecteur ligne de $\mathbb{K}^{2m}$, on désigne par
$\langle\alpha,.\rangle_{1}$ et $\langle\alpha,.\rangle_{2}$ les
formes linéaires sur $\mathfrak{n}$ définies par
\begin{align*}
\langle\alpha,Z_{\beta,\gamma,y,z}\rangle_{1} &=
\alpha\,^{t}\beta\\
\langle\alpha,Z_{\beta,\gamma,y,z}\rangle_{2} &=
\alpha\,^{t}\gamma.
\end{align*}
 Alors, si
$\lambda_{1},\lambda_{2}\in\mathbb{K}$, on a pour
$\beta,\gamma\in\mathbb{K}^{2m}$ et $y,z\in\mathbb{K}$,
\begin{equation}
\sideset{}{^{*}}\ad
Z_{\beta,\gamma,y,z}.(\lambda_{1}Z_{1}^{*}+\lambda_{2}Z_{2}^{*})
=-\langle2\lambda_{1}\beta J_{2m}+\lambda_{2}\gamma J_{2m},\rangle_{1}
-\langle\lambda_{2}\beta J_{2m},\rangle_{2}.\label{eq3b'''1}
\end{equation}
Il s'ensuit que $\{\lambda_{1}Z_{1}^{*}+\lambda_{2}Z_{2}^{*} :
\lambda_{2}\neq0\}$ est un ensemble de représentants des orbites
régulières dans $\mathfrak{n}^{*}$. Plus précisément, si
$n=\lambda_{1}Z_{1}^{*}+\lambda_{2}Z_{2}^{*}$ avec $\lambda_{2}\neq0$,
on a $\mathfrak{n}(n)=\mathfrak{z}$.

D'autre part, pour $(Y_{A,a,\alpha,y},D,d)\in\mathfrak{r}_{\cur{V}''}
\times\mathfrak{sp}(V''')\times\mathbb{K}$ on a
\begin{equation}
\sideset{}{^{*}}\ad(Y_{A,a,\alpha,y},D,d).
(\lambda_{1}Z_{1}^{*}+\lambda_{2}Z_{2}^{*})=
-2a\lambda_{1}Z_{1}^{*}-(a+d)\lambda_{2}Z_{2}^{*}.\label{eq3b'''2}
\end{equation}
L'application $\chi:Y_{A,a,\alpha,y}\mapsto a$ est un caractère de
$\mathfrak{r}_{\cur{V}''}$ dont on note $\mathfrak{r}^{0}_{\cur{V}''}$
le noyau.

Soit $g$ une forme fortement régulière sur
$\mathfrak{g}=\mathfrak{r}_{\cur{V}}$. Par généricité et compte tenu
de la relation \ref{eq3b'''1}, on peut supposer que sa restriction à
$\mathfrak{n}$ est de la forme
$n=\lambda_{1}Z_{1}^{*}+\lambda_{2}Z_{2}^{*}$ avec
$\lambda_{1}\lambda_{2}\neq 0$. Il suit alors des relations
\ref{eq3b'''1} et \ref{eq3b'''2} que, reprenant les notations du
lemme \ref{lem3a1}, on a
\begin{align*}
\mathfrak{h}=\mathfrak{g}(n) & =
(\mathfrak{r}^{0}_{\cur{V}''}\times\mathfrak{sp}(V'''))+\mathfrak{z},\\
\mathfrak{g}_{1}=\mathfrak{h}/\mathfrak{q}
&=\mathfrak{r}^{0}_{\cur{V}''}\times\mathfrak{sp}(V''').
\end{align*}
On voit alors que
$\dim\mathfrak{g}/(\mathfrak{h}+\mathfrak{n})=2$. Comme
$\dim\mathfrak{q}=1$ et
$\rang\mathfrak{sp}(V''')=\ind\mathfrak{sp}(V''')=m$, on déduit du lemme
cité que
\begin{align}
\ind\mathfrak{r}_{\cur{V}}& =\ind\mathfrak{r}^{0}_{\cur{V}''}
+m-1,\label{eq3b'''3}\\
\rang\mathfrak{r}_{\cur{V}}& =\rang\mathfrak{r}^{0}_{\cur{V}''}
+m.\label{eq3b'''4}
\end{align}

\subsection{}\label{3b''''}
Dans ce numéro, avant de conclure, nous allons comparer l'indice et le
rang de $\mathfrak{r}^{0}_{\cur{V}''}$ à ceux de
$\mathfrak{r}_{\cur{V}''}$. On pose
$\mathfrak{g}=\mathfrak{r}_{\cur{V}''}$ et
$\mathfrak{n}=\mathbb{K}Z_{1}$ où $Z_{1}=Z_{0,0,1,0}=Y_{0,0,0,1}$ est
l'élément de $\mathfrak{r}_{\cur{V}''}$ introduit au numéro
précédent. Alors, on a $[Y_{A,a,\alpha,y},Z_{1}]=2aZ_{1}$, de sorte
que $\mathfrak{n}$ est un idéal unipotent de
$\mathfrak{r}_{\cur{V}''}$.

Soit $g$ une forme fortement régulière sur $\mathfrak{g}$. On peut
supposer que $g(Z_{1})\neq 0$. Désignant par $n$ la restriction de $g$ à
$\mathfrak{n}$, on voit donc que $\mathfrak{h}=\mathfrak{g}(n)$ est le
centralisateur de $Z_{1}$ dans $\mathfrak{g}$, soit
$\mathfrak{h}=\mathfrak{r}^{0}_{\cur{V}''}$. Appliquant le lemme
\ref{lem3a1} on voit que l'on a
\begin{align}
\ind\mathfrak{r}^{0}_{\cur{V}''}&
=1+\ind\mathfrak{r}_{\cur{V}''},\label{eq3b''''1}\\
\rang\mathfrak{r}^{0}_{\cur{V}''}
&=\rang\mathfrak{r}_{\cur{V}''}.\label{eq3b''''2}
\end{align}
Rassemblant ces formules avec les formules \ref{eq3b'''3} et \ref{eq3b'''4},
on obtient
\begin{align}
\ind\mathfrak{r}_{\cur{V}}& =\ind\mathfrak{r}_{\cur{V}''}
+m,\label{eq3b''''3}\\
\rang\mathfrak{r}_{\cur{V}}& =\rang\mathfrak{r}_{\cur{V}''}
+m.\label{eq3b''''4}
\end{align}
Par ailleurs, l'algèbre $\mathfrak{r}_{\cur{V}}$ contient le tore
constitué des matrices de la forme $x_{A,a,0,0,D,d,0,0}$ où $A$
(resp. $D$) est une matrice diagonale dans $\mathfrak{sp}(V'')$
(resp. $\mathfrak{sp}(V''')$) et $a,d\in\mathbb{K}$. Or le
centralisateur de ce tore dans $\mathfrak{gl}(V)$ est le tore des
matrices diagonales. Il en résulte que le radical unipotent du centre
de $\mathfrak{r}_{\cur{V}}$ est trivial, si bien que cette dernière
est quasi-réductive si et seulement si son rang est égal à son indice.

Comme $V''$ est un espace symplectique, l'assertion a) de la
proposition s'applique à l'algèbre $\mathfrak{r}_{\cur{V}''}$. Dans le
cas considéré, l'assertion b) de la proposition résulte alors des
formules \ref{eq3b''''3} et \ref{eq3b''''4} et du fait que
$h(\cur{V}')=h(\cur{V}'')$.
\end{dem}

\subsection{Indice et quasi-réductivité des sous algèbres paraboliques
de $\mathfrak{so}(E)$.}\label{3c}

Soit $E$ un espace vectoriel de dimension finie supérieure ou égale à
$3$ sur $\mathbb{K}$ muni d'une forme bilinéaire symétrique non
dégénérée $B$. On désigne par $\mathfrak{so}(E)$ l'algèbre de Lie du
groupe orthogonal correspondant. Soit $\cur{V}=\{\{0\}=V_{0}\subsetneq
V_{1}\subsetneq\cdots\subsetneq V_{t}=V\}$ un drapeau de sous-espaces
isotropes de $E$ avec $\dim V\geq 1$. On désigne par
$\mathfrak{p}_{\cur{V}}$ la sous-algèbre de Lie de $\mathfrak{so}(E)$
constituée des endomorphismes de $E$ stabilisant le drapeau
$\cur{V}$~: c'est une sous-algèbre parabolique de $\mathfrak{so}(E)$
et toutes les sous-algèbres paraboliques de $\mathfrak{so}(E)$ sont
obtenues ainsi. Lorsque le drapeau $\cur{V}=\{\{0\}=V_{0}\subsetneq
V_{1}=V\}$ ne contient qu'un seul sous-espace isotrope non nul $V$, on
note la sous-algèbre parabolique correspondante plus simplement
$\mathfrak{p}_{V}$ au lieu de $\mathfrak{p}_{\cur{V}}$.

Si $u,v\in E$, on définit l'élément $u\wedge_{B}v$ de $\mathfrak{so}(E)$
en posant, pour $x\in E$~:
\begin{equation*}
u\wedge_{B}v(x)=B(x,u)v-B(x,v)u.
\end{equation*}
Il existe une unique application linéaire $\wedge_{B}:E\otimes
E\mapsto\mathfrak{so}(E)$ telle que $\wedge_{B}(u\otimes
v)=u\wedge_{B}v$, laquelle passe au quotient en un isomorphisme
d'espaces vectoriels de $\wedge^{2}E$ sur $\mathfrak{so}(E)$. Plus
généralement, si $V,W$ sont des sous-espaces de $E$ tels que $V\cap
W=0$, $\wedge_{B}$ induit un isomorphisme d'espaces vectoriels de
$\wedge^{2}V$ (resp. $V\otimes W$) sur un sous-espace de
$\mathfrak{so}(E)$ noté $\wedge_{B}^{2}V$ (resp. $V\wedge_{B}W$).

Si $u_{1},u_{2},v_{1},v_{2}\in E$, on a
\begin{equation*}
\begin{split}
[u_{1}\wedge_{B}v_{1},u_{2}\wedge_{B}v_{2}]&=
B(u_{1},u_{2})v_{1}\wedge_{B}v_{2}+B(u_{1},v_{2})u_{2}\wedge_{B}v_{1}\\
&\quad +B(v_{1},u_{2})u_{1}\wedge_{B}v_{2}+B(v_{1},v_{2})u_{1}\wedge_{B}u_{2}
\end{split}
\end{equation*}

On munit $\mathfrak{so}(E)$ de la forme bilinéaire symétrique non
dégénérée et invariante $\kappa(X,Y)=-\frac{1}{2}\Tr XY$. Si
$u_{1},u_{2},v_{1},v_{2}\in E$, on a
\begin{equation*}
\kappa(u_{1}\wedge_{B}v_{1},u_{2}\wedge_{B}v_{2})=\det
\begin{bmatrix}
B(u_{1},u_{2}) & B(u_{1},v_{2})\\
B(v_{1},u_{2}) & B(v_{1},v_{2})
\end{bmatrix}.
\end{equation*}

\subsection{}\label{3d}
Soit $V$ un sous-espace isotrope non nul de $E$. On désigne par
$\mathfrak{n}_{V}$ le radical unipotent de la sous-algèbre parabolique
$\mathfrak{p}_{V}$ et par $\mathbold{N}_{V}$ le sous-groupe unipotent
de $\mathrm{GL}(V)$ d'algèbre de Lie $\mathfrak{n}_{V}$.

Soit $W$ un sous-espace isotrope supplémentaire dans $E$ de
l'orthogonal $V^{\perp_{B}}$ de $V$ pour $B$. Alors $B$ induit une
dualité entre $V$ et $W$, la restriction de $B$ à $(V\oplus
W)^{\perp_{B}}$ est non dégénérée et on a $E=V\oplus(V\oplus
W)^{\perp_{B}}\oplus W$. Si de plus $F$ est un sous-espace de $E$ tel que
$B_{\mid F}$ soit non dégénérée, on identifie les algèbres de Lie
$\mathfrak{gl}(V)$ et $\mathfrak{so}(F)$ à des sous-algèbres de Lie de
$\mathfrak{so}(E)$ de la manière suivante~:

- si $X\in\mathfrak{gl}(V)$, on l'étend en un élément encore noté $X$ de
$\mathfrak{so}(E)$ en décidant que $X_{\mid (V\oplus W)^{\perp_{B}}}=0$
et, si $y\in W$, $X.y$ est l'élément de $W$ tel que
\begin{equation*}
B(X.y,x)+B(y,X.x)=0\mbox{, }x\in V,
\end{equation*}

- si $X\in\mathfrak{so}(F)$, on l'étend en un élément encore noté $X$ de
$\mathfrak{so}(E)$ en décidant que $X_{\mid F^{\perp_{B}}}=0$.

Alors, si on prend $F=(V\oplus W)^{\perp_{B}}$, l'application
$(X,Y)\mapsto X+Y$ permet d'identifier
$\mathfrak{gl}(V)\times\mathfrak{so}(F)$ à une sous-algèbre de Lie
$\mathfrak{l}_{V}$ de $\mathfrak{gl}(E)$ qui
%\begin{equation*}
%\mathfrak{l}_{V}=\mathfrak{gl}(V)\oplus\mathfrak{so}(F)
%\end{equation*}
est un facteur réductif de $\mathfrak{p}_{V}$.  D'autre part, on a
$V^{\perp_{B}}=F\oplus V$. L'application $\wedge_{B}$ envoie
$V^{\perp_{B}}\otimes V$ sur $\mathfrak{n}_{V}$, induisant, lorsque
$\dim V>1$, un isomorphisme de $\wedge^{2}V$ sur le centre
$\mathfrak{z}_{V}$ de $\mathfrak{n}_{V}$ de sorte qu'alors
$\mathfrak{z}_{V}=\wedge_{B}^{2}V$. Par ailleurs, on a
$\mathfrak{n}_{V}=F\wedge_{B} V\oplus\mathfrak{z}_{V}$, cette
décomposition étant $\mathfrak{l}_{V}$-stable, si $\dim V>1$, tandis
que $\mathfrak{n}_{V}=F\wedge_{B} V=\mathfrak{z}_{V}$ est isomorphe à
$F$, sinon. L'algèbre de Lie $\mathfrak{n}_{V}$ est d'indice de
nilpotence 2, si $\dim V>1$ et $\dim F>0$, commutative sinon. De plus,
si $X\in\mathfrak{gl}(V)$, $Y\in\mathfrak{so}(F)$, $u\in F$, $v,v'\in V$,
on a
\begin{align*}
[X,u\wedge_{B} v] &= u\wedge_{B}X.v\\
[Y,u\wedge_{B} v] &= Y.u\wedge_{B}v\\
[Y,v\wedge_{B} v'] &= 0\\
[X,v\wedge_{B} v'] &= X.v\wedge_{B}v'+v\wedge_{B}X.v'.
\end{align*}

La forme $\kappa$ permet d'identifier $\wedge_{B}^{2}V^{*}$
(resp. $(F\wedge_{B} V)^{*}$) de manière
$\mathfrak{l}_{V}$-équivariante avec $\wedge_{B}^{2}W$
(resp. $F\wedge_{B} W$). On a alors les décompositions duales
\begin{align*}
\mathfrak{n}_{V}&=F\wedge_{B} V\oplus\wedge_{B}^{2}V\\
\mathfrak{n}_{V}^{*}&=F\wedge_{B} W\oplus\wedge_{B}^{2}W.
\end{align*}

Dans les numéros qui suivent, nous allons déterminer les
stabilisateurs dans $\mathfrak{p}_{V}$ des éléments génériques $n$ de
$\mathfrak{n}_{V}^{*}$ ainsi que leur quotient par le noyau
$\mathfrak{q}$ de $n_{\mid\mathfrak{n}_{V}(n)}$.

\subsection{}\label{3d'}
Dans ce numéro, on suppose que $\dim V=1$, de sorte que
$\mathfrak{gl}(V)=\mathbb{K}$. Le choix de bases, $e_{1}$ de $V$ et
$f_{1}$ de $W$, telles que $B(e_{1},f_{1})=1$ permet d'identifier $F$ avec
$\mathfrak{n}_{V}$ (resp. $\mathfrak{n}_{V}^{*}=F\wedge_{B}^{2}W$)
au moyen de l'application $u\mapsto u\wedge_{B}e_{1}$
(resp. $u\mapsto u\wedge_{B}f_{1}$). Alors,
$\mathfrak{l}_{V}=\mathbb{K}\times\mathfrak{so}(F)$ et
$\mathfrak{p}_{V}=(\mathbb{K}\times\mathfrak{so}(F))\oplus F$,
$\mathbb{K}$ (resp. $\mathfrak{so}(F)$) agissant par multiplication dans
$F$ (resp.  de manière naturelle). Les éléments génériques de
$\mathfrak{n}_{V}^{*}$ sont, modulo l'identification précédente, les
vecteurs non isotropes de $F$. Soit donc $n\in F=\mathfrak{n}_{V}^{*}$
un vecteur non isotrope et $F_{0}$ son orthogonal dans $F$
relativement à $B$, qui n'est autre que $\mathfrak{q}$. Alors, on a
\begin{align}
\mathfrak{p}_{V}(n) & =\mathfrak{so}(F_{0})\oplus F\label{eq3e'1}\\
\mathfrak{p}_{V}(n)/\mathfrak{q} &=
\mathfrak{so}(F_{0})\oplus F/F_{0},\label{eq3e'2}
\end{align}
$F/F_{0}$ étant un idéal unipotent central.

\subsection{}\label{3e}
Dans ce numéro et les trois suivants, on suppose que $\dim V>1$.  Soit
$\xi\in\wedge_{B}^{2}W$ un élément de rang maximal. Nous allons
déterminer un système de représentants des $\mathbold{N}_{V}$-orbites
dans $\mathfrak{n}_{V}^{*}$ dont la restriction à $\mathfrak{z}_{V}$
est $\xi$.

On choisit un système libre $e_{1},\ldots,e_{s},f_{1},\ldots,f_{s}$ de
$E$ tel que $s=[\frac{1}{2}\dim E]$, que
$B(e_{i},e_{j})=B(f_{i},f_{j})=0$ et $B(e_{i},f_{j})=\delta_{ij}$,
$1\leq i,j\leq s$, que $e_{1},\ldots,e_{r}$ soit une base de $V$,
$f_{1},\ldots,f_{r}$ une base de $W$ et que
$\xi=\sum_{i=1}^{[\frac{r}{2}]}f_{2i-1}\wedge_{B}f_{2i}$. Alors
$e_{r+1},\ldots,e_{s}$, $f_{r+1},\ldots,f_{s}$ est soit une base de
$F$, si $E$ est de dimension paire, soit se complète en une base
$e_{r+1},\ldots,e_{s}, f_{r+1},\ldots,f_{s},\tilde{u}$ de $F$, où $\tilde{u}$
est un vecteur de $F$ tel que $B(\tilde{u},\tilde{u})=1$ et
$B(\tilde{u},e_{i})=B(\tilde{u},f_{i})=0$, $1\leq i\leq s$, sinon.

On considère $\xi$ comme un élément de $\mathfrak{n}_{V}^{*}$ noté
alors $n_{\xi}$. Soit $T\in F\wedge_{B}V$ que l'on écrit
$T=\sum_{i=1}^{r}u_{i}\wedge_{B}e_{i}$, avec $u_{i}\in F$, $1\leq
i\leq r$. Alors, on a
\begin{equation*}
T.n_{\xi}=\sum_{i=1}^{[\frac{r}{2}]}(u_{2i}\wedge_{B}f_{2i-1}-
u_{2i-1}\wedge_{B}f_{2i}).
\end{equation*}
On en déduit que

- si $r$ est pair, $\mathfrak{n}_{V}(n_{\xi})=\mathfrak{z}_{V}$ et
$\mathbold{N}_{V}.n_{\xi}=n_{\xi}+F\wedge_{B}W=
n_{\xi}+\mathfrak{z}_{V}^{\perp}$
est l'unique $\mathbold{N}_{V}$-orbite dont la restriction à
$\mathfrak{z}_{V}$ est $\xi$,

- si $r$ est impair, $\mathfrak{n}_{V}(n_{\xi})=
F\wedge_{B}e_{r}\oplus\mathfrak{z}_{V}$ et
$\mathbold{N}_{V}.n_{\xi}=n_{\xi}+\sum_{i=1}^{r-1}F\wedge_{B}f_{i}$. Comme
les éléments de $F\wedge_{B}W$ sont fixés par $\mathbold{N}_{V}$, on
voit que $n_{\xi}+F\wedge_{B}f_{r}$ est un ensemble de représentants
des $\mathbold{N}_{V}$-orbites dans $\mathfrak{n}_{V}^{*}$ dont la
restriction à $\mathfrak{z}_{V}$ est $\xi$. De plus, pour tout $n\in
n_{\xi}+F\wedge_{B}f_{r}$, on a $\mathfrak{n}_{V}(n)=
F\wedge_{B}e_{r}\oplus\mathfrak{z}_{V}$.

\subsection{}\label{3f}
Dans ce numéro, on suppose que $r=2p$ est pair avec $p\geq1$. On
désigne alors par $\mathfrak{n}{_{V}^{*}}'$ l'ouvert de Zariski de
$\mathfrak{n}_{V}^{*}$ constitué des formes linéaires $n$ telles que
$n_{\mid\mathfrak{z}_{V}}$ soit une forme symplectique $\xi$ sur
$V$. Soit $n\in\mathfrak{n}{_{V}^{*}}'$. Compte tenu de ce qui précède
on a $\mathfrak{n}_{V}(n)=\mathfrak{z}_{V}$ et, quitte à translater
par un élément de $\mathbold{N}_{V}$, on peut supposer que
$n=n_{\xi}$. On désigne par $\mathfrak{sp}(V,\xi)$ ou plus simplement
par $\mathfrak{sp}(V)$ s'il n'y a pas risque de confusion, la
sous-algèbre de Lie de $\mathfrak{gl}(V)$ qui est l'algèbre de Lie du
groupe symplectique correspondant. Alors, on voit immédiatement que
\begin{equation*}
\mathfrak{p}_{V}(\xi)=
(\mathfrak{sp}(V)\times\mathfrak{so}(F))\oplus\mathfrak{n}_{V}.
\end{equation*}
On en déduit facilement que
\begin{align}%\label{eq3f1}
\mathfrak{p}_{V}(n) &=
(\mathfrak{sp}(V)\times\mathfrak{so}(F))\oplus\mathfrak{z}_{V}\label{eq3f1}\\
\mathfrak{p}_{V}(n)/\mathfrak{q} &=
(\mathfrak{sp}(V)\times\mathfrak{so}(F))\oplus\mathfrak{z}_{V}/\mathfrak{q},
\label{eq3f2}
\end{align}
l'idéal $\mathfrak{z}_{V}$ (resp. $\mathfrak{z}_{V}/\mathfrak{q}$)
étant unipotent (resp. unipotent et central).

\subsection{}\label{3f'}
Dans ce numéro, on suppose que $r=2p+1$ est impair avec $p\geq1$ et
que $F$ est nul, de sorte que $\mathfrak{n}_{V}=\wedge^{2}_{B}V$ et
$\mathfrak{n}_{V}^{*}=\wedge^{2}_{B}W$. On désigne par
$\mathfrak{n}{_{V}^{*}
}'$ l'ouvert de Zariski de $\mathfrak{n}_{V}^{*}$
constitué des $\xi$ de rang maximal. On prend
$n=\xi\in\mathfrak{n}{_{V}^{*}
}'$. Alors, on a
\begin{align}%\label{eq3f'1}
\mathfrak{p}_{V}(n) &=
\mathfrak{gl}(V)(\xi)\oplus\mathfrak{n}_{V}\label{eq3f'1}\\
\mathfrak{p}_{V}(n)/\mathfrak{q} &=
\mathfrak{gl}(V)(\xi)\oplus\mathfrak{n}_{V}/\mathfrak{q},\label{eq3f'2}
\end{align}
l'idéal $\mathfrak{n}_{V}$ (resp. $\mathfrak{n}_{V}/\mathfrak{q}$)
étant unipotent (resp. unipotent et central).

\subsection{}\label{3g}
Dans ce numéro, on suppose que $r=2p+1$ est impair avec $p\geq1$ et
que $F$ est non nul. Un élément $n$ de $\mathfrak{n}_{V}^{*}$ s'écrit
de manière unique sous la forme $n=\xi+n_{1}$, avec
$\xi\in\wedge_{B}^{2}W$ et $n_{1}\in F\wedge_{B}W$. Alors, $n_{1}$, vu
comme un élément de $\mathfrak{so}(E)$ envoie $V$ dans $F$. On désigne
par $\mathfrak{n}{_{V}^{*}
}'$ le sous-ensemble de $\mathfrak{n}_{V}^{*}$
constitué des $n$ pour lesquels $\xi$ est une forme bilinéaire
alternée de rang maximal $2p$ sur $V$ et, si $e$ est un élément non
nul du noyau de $\xi$ dans $V$, $n_{1}(e)$ est un vecteur non isotrope
de $F$. Alors $\mathfrak{n}{_{V}^{*}
}'$ est un ouvert de Zariski non
vide $\mathbold{N}_{V}$-invariant de $\mathfrak{n}_{V}^{*}$.

Soit $n\in\mathfrak{n}{_{V}^{*}
}'$ et $\xi\in\wedge_{B}^{2}W$ la
restriction de $n$ à $\mathfrak{z}_{V}$. On reprend les notations du
numéro \ref{3e}. Compte tenu des résultats de celui-ci, quitte à
translater par un élément de $\mathbold{N}_{V}$, on peut supposer que
$n=n_{\xi}+v_{r}\wedge_{B}f_{r}$ avec $v_{r}\in F$ tel que
$B(v_{r},v_{r})\neq0$.

Commençons par étudier $\mathfrak{p}_{V}(\xi)$. On a clairement
\begin{equation*}
\mathfrak{p}_{V}(\xi)=
(\mathfrak{gl}(V)(\xi)\times\mathfrak{so}(F))\oplus\mathfrak{n}_{V}.
\end{equation*}
\'Ecrivons $V=V'\oplus\mathbb{K}e_{r}$ avec
$V'=\mathbb{K}e_{1}\oplus\cdots\oplus\mathbb{K}e_{2p}$. Alors, la
matrice de $\xi$ dans la base $e_{1},\ldots,e_{r}$ est de la forme
\begin{equation*}
\begin{bmatrix}
J_{2p} & 0\\
0 & 0
\end{bmatrix}
\end{equation*}
de sorte que $\xi_{\mid V'}$ est une forme symplectique et que, si
on identifie les endomorphismes de $V$ avec leur matrice dans cette
base, les éléments de $\mathfrak{gl}(V)(\xi)$ s'écrivent
\begin{equation*}
X_{A,a,c}=
\begin{bmatrix}
A & 0\\
c & a
\end{bmatrix}
\end{equation*}
avec $A\in\mathfrak{sp}_{2p}$, $c\in\mathbb{K}^{2p}$ un vecteur ligne et
$a\in\mathbb{K}$. Soit, $X_{A,a,c}\in\mathfrak{gl}(V)(\xi)$,
$Y\in\mathfrak{so}(F)$ et $T=\sum_{i=1}^{r}u_{i}\wedge_{B}e_{i}\in
F\wedge_{B}V$. Alors, on a
\begin{equation*}
\begin{split}
(X_{A,a,c}+Y+T).n
&=\sum_{i=1}^{p}[(u_{2i}-c_{2i-1}v_{r})\wedge_{B}f_{2i-1}
  -(u_{2i-1}+c_{2i}v_{r})\wedge_{B}f_{2i}]\\ &\quad
+(Y-a)v_{r}\wedge_{B}f_{r}.
\end{split}
\end{equation*}
On en déduit que $X_{A,a,c}+Y+T\in\mathfrak{p}_{V}(n)$ si et seulement
si $a=0$, $Y.v_{r}=0$ et $u_{2i}=c_{2i-1}v_{r}$,
$u_{2i-1}=-c_{2i}v_{r}$, $1\leq i\leq p$.

On désigne par $F_{0}$ l'orthogonal de $v_{r}$ dans $F$. Alors
$B_{\mid F_{0}}$ est non dégénérée de sorte que l'on peut identifier
$\mathfrak{so}(F_{0})$ à une sous-algèbre de Lie de $\mathfrak{so}(F)$. On
identifie $\mathfrak{sp}(V')$ avec la sous-algèbre de Lie de
$\mathfrak{gl}(V)$ constituée des matrices de la forme $X_{A,0,0}$,
$A\in\mathfrak{sp}_{2p}$. Si $c\in\mathbb{K}^{2p}$, on pose
$T(c)=v_{r}\wedge_{B}\sum_{i=1}^{p}(c_{2i-1}e_{2i}-c_{2i}e_{2i-1})$,
$\tilde{T}(c)=X_{0,c,0}+T(c)$ et on désigne par $\mathfrak{t}_{V'}$ le
sous-espace vectoriel de $\mathfrak{so}(E)$ constitué des
$\tilde{T}(c)$, $c\in\mathbb{K}^{2p}$. Alors, on a
\begin{align}
\mathfrak{p}_{V}(n)&=(\mathfrak{sp}(V')\times\mathfrak{so}(F_{0}))\oplus
\mathfrak{t}_{V'}\oplus\mathfrak{n}_{V}(n)\label{eq3g1}\\
%& =\mathfrak{so}(F_{0})\times(\mathfrak{sp}(V')\oplus\mathfrak{n}_{V}(n)),
%\notag\\
\mathfrak{n}_{V}(n)&=F\wedge_{B}e_{r}\oplus\mathfrak{z}_{V},\label{eq3g2}
\end{align}
$\mathfrak{t}_{V'}\oplus\mathfrak{n}_{V}(n)$ étant le radical
unipotent de $\mathfrak{p}_{V}(n)$ et
$\mathfrak{sp}(V')\times\mathfrak{so}(F_{0})$ en étant un facteur
réductif. %De plus $\mathfrak{p}_{V}(n)$ est le produit direct de ses
%sous-algèbres de Lie $\mathfrak{so}(F_{0})$ et
%$\mathfrak{sp}(V')\oplus\mathfrak{t}_{V'}\oplus\mathfrak{n}_{V}(n)$~:
%\begin{equation}
%\mathfrak{p}_{V}(n)=\mathfrak{so}(F_{0})\times(\mathfrak{sp}(V')
%\oplus\mathfrak{t}_{V'}\oplus\mathfrak{n}_{V}(n)).\label{eq3g2}
%\end{equation}

Soit $\mathfrak{q}$ le noyau de $n_{\mid \mathfrak{n}_{V}(n)}$. Alors
$\mathfrak{p}_{V}(n)/\mathfrak{q}$ est le produit direct de
$\mathfrak{so}(F_{0})$ et de l'algèbre de Lie quotient $\mathfrak{s}_{V'}=
\mathfrak{sp}(V')\oplus\mathfrak{t}_{V'}
\oplus\mathfrak{n}_{V}(n)/\mathfrak{q}$ dont un facteur réductif est
$\mathfrak{sp}(V')$ et dont le radical unipotent est $\mathfrak{u}_{V'}=
\mathfrak{t}_{V'}\oplus\mathfrak{n}_{V}(n)/\mathfrak{q}$.

Soit $c\in\mathbb{K}^{2p}$ et $A\in\mathfrak{sp}_{2p}$. Alors, on a
\begin{equation*}
[X_{A,0,0},\tilde{T}(c)]=\tilde{T}(-cA).
\end{equation*}
On voit donc que le sous-espace $\mathfrak{t}_{V'}$ est invariant sous
l'action du facteur réductif $\mathfrak{sp}(V')$ et que, si
$e_{1}^{*},\ldots,e_{2p}^{*}$ est la base duale de la base
$e_{1},\ldots,e_{2p}$ de $V'$, l'application $\tilde{T}(c)\mapsto
c^{*}=\sum_{i=1}^{2p}c_{i}e_{i}^{*}$ est un isomorphisme de
$\mathfrak{sp}(V')$-modules de $\mathfrak{t}_{V'}$ sur $V'^{*}$.

D'autre part, si $c,c'\in\mathbb{K}^{2p}$, on a
\begin{equation*}
\begin{split}
[\tilde{T}(c),\tilde{T}(c')] & =[T(c),T(c')]\\
& =
\sum_{1\leq i<j\leq p}
\{(c_{2i-1}c'_{2j-1}-c'_{2i-1}c_{2j-1})e_{2i}\wedge_{B}e_{2j}+
(c_{2i}c'_{2j}-c'_{2i}c_{2j})e_{2i-1}\wedge_{B}e_{2j-1}\} \\
&\quad +\sum_{1\leq
i,j\leq p}(c_{2i}c'_{2j-1}-c'_{2i}c_{2j-1})e_{2j}\wedge_{B}e_{2i-1}
\end{split}
\end{equation*}
de sorte que
\begin{equation*}
\langle n,[\tilde{T}(c),\tilde{T}(c')]\rangle=\sum_{i=1}^{p}
(c_{2i-1}c'_{2i}-c'_{2i-1}c_{2i}).
\end{equation*}

La forme symplectique $\xi$ induit un isomorphisme $\tilde{\xi}$
d'espaces vectoriels de $V'$ sur $V'^{*}$ lequel envoie l'élément $x$
de $V'$ sur la forme linéaire $y\mapsto\xi(x,y)$. Cet isomorphisme
transporte $\xi$ en une forme symplectique sur $V'^{*}$ encore notée
$\xi$ telle que $\xi=\sum_{i=1}^{p}e_{2i-1}\wedge e_{2i}$. Soit alors
$\mathfrak{h}_{V'^{*}}=V'^{*}\oplus\mathbb{K}Z$ l'algèbre de
Heisenberg construite sur $V'^{*}$comme au numéro \ref{1g2} et munie
de l'action par dérivation de $\mathfrak{sp}(V')$ introduite au numéro
\ref{1g3}.

On déduit des remarques  précédentes que
l'application
\begin{equation}\label{eq3g3}
\tilde{T}(c)+Y+\mathfrak{q}\mapsto c^{*}+\langle n,Y\rangle Z\mbox{,
}c\in\mathbb{K}^{2p}\mbox{, }Y\in\mathfrak{n}_{V}(n),
\end{equation}
est un isomorphisme $\mathfrak{sp}(V')$-équivariant d'algèbres de Lie de
$\mathfrak{u}_{V'}=\mathfrak{t}_{V'}\oplus\mathfrak{n}_{V}(n)/\mathfrak{q}$
sur $\mathfrak{h}_{V'^{*}}$. En particulier, l'algèbre de Lie
$\mathfrak{s}_{V'}$ est isomorphe au produit semi-direct de
$\mathfrak{sp}(V')$ par $\mathfrak{h}_{V'^{*}}$.

\subsection{}\label{3h}
Dans ce numéro, nous énonçons le résultat principal de cette section.
\begin{theo}\label{theo3h1}
Soit $E$ un espace vectoriel de dimension $q\geq3$ sur $\mathbb{K}$
muni d'une forme bilinéaire symétrique non dégénérée $B$,
$\cur{V}=\{\{0\}=V_{0}\subsetneq V_{1}\subsetneq\cdots\subsetneq
V_{t}=V\}$ un drapeau de sous-espaces isotropes de $E$ avec $r=\dim
V\geq1$. On note $\cur{V}'$ le drapeau tel que
$\cur{V'}=\cur{V}-\{V\}$, si $r$ est impair égal à $\frac{q}{2}$, et
$\cur{V'}=\cur{V}$, sinon.

(i) On a les formules suivantes pour l'indice de la sous-algèbre
parabolique $\mathfrak{p}_{\cur{V}}$~:

\begin{equation}\label{eq3h1}
\ind\mathfrak{p}_{\cur{V}}=\begin{cases}[\frac{q}{2}]-r+
\sum_{i=1}^{t}[\frac{1}{2}(\dim V_{i}-\dim V_{i-1})] & \text{si $r$
est pair,}\\
[\frac{q-1}{2}]-r+
\sum_{i=1}^{t}[\frac{1}{2}(\dim V_{i}-\dim V_{i-1})] & \text{si $r$
est impair et $r<\frac{q}{2}$,}\\
\sum_{i=1}^{t}[\frac{1}{2}(\dim V_{i}-\dim V_{i-1})]-1 & \text{si $r$
est impair, $r=\frac{q}{2}$ et}\\ &\text{$\dim V_{t-1}<r-1$,}\\
\sum_{i=1}^{t}[\frac{1}{2}(\dim V_{i}-\dim V_{i-1})]+1 & \text{si $r$
est impair, $r=\frac{q}{2}$ et}\\ &\text{$\dim V_{t-1}=r-1$,}
\end{cases}
\end{equation}

(ii) La dimension du radical unipotent d'un stabilisateur générique de
la représentation coadjointe de $\mathfrak{p}_{\cur{V}}$ est
$h(\cur{V}')$.

(iii) L'algèbre de Lie $\mathfrak{p}_{\cur{V}}$ est quasi-réductive si
et seulement si le drapeau $\cur{V'}$ vérifie la propriété $\cur{P}$.
\end{theo}

\begin{dem}
Commençons par remarquer que, le centre d'une sous-algèbre de Lie
parabolique d'une algèbre de Lie simple étant trivial,
$\mathfrak{p}_{\cur{V}}$ sera quasi-réductive si et seulement si
$\ind\mathfrak{p}_{\cur{V}}=\rang\mathfrak{p}_{\cur{V}}$.  Il est donc
clair que (iii) est conséquence de (ii).

Reprenons maintenant les notations des numéros précédents et plus
particulièrement celles des numéros \ref{3d} et suivants.  L'algèbre
de Lie $\mathfrak{p}_{\cur{V}}$ est une sous-algèbre de Lie de
$\mathfrak{p}_{V}$ et ce sont toutes deux des sous-algèbres
paraboliques de $\mathfrak{so}(E)$. Par suite $\mathfrak{n}_{V}$ est un
idéal unipotent de $\mathfrak{p}_{\cur{V}}$ et l'on a
\begin{equation*}
\mathfrak{p}_{\cur{V}}=\mathfrak{p}_{\cur{V}}\cap\mathfrak{l}_{V}
\oplus\mathfrak{n}_{V}=(\mathfrak{q}_{\cur{V}}\times\mathfrak{so}(F))
\oplus\mathfrak{n}_{V}
\end{equation*}
où $\mathfrak{q}_{\cur{V}}$ désigne la sous-algèbre parabolique de
$\mathfrak{gl}(V)$ stabilisant le drapeau $\cur{V}$.

La démonstration consiste à appliquer plusieurs fois le lemme
\ref{lem3a1}, dont une première avec
$\mathfrak{n}=\mathfrak{n}_{V}$.

Soit donc $g$ une forme linéaire fortement régulière sur
$\mathfrak{g}=\mathfrak{p}_{\cur{V}}$, $n$ (resp. $\xi$) sa restriction à
$\mathfrak{n}=\mathfrak{n}_{V}$ (resp. $\mathfrak{z}_{V}$). On reprend
les notations $\mathfrak{h}$ et $\mathfrak{g}_{1}$ du numéro \ref{3a}.

\subsection{}\label{3i}
Dans ce numéro, on suppose que $r$ est pair. On a donc
$\cur{V}=\cur{V'}$. Par généricité de $g$, on peut supposer que $\xi$
est une forme symplectique sur $V$ générique relativement au drapeau
$\cur{V}$. Mais alors, il suit du numéro \ref{3e} que, quitte à
translater $g$ par un élément de $\mathbold{N}_{V}$, on peut également
supposer que $n=n_{\xi}$ de sorte que, compte tenu de \ref{eq3f1} et
\ref{eq3f2}, on a
$\mathfrak{h}=(\mathfrak{r}_{\cur{V}}\times\mathfrak{so}(F))
\oplus\mathfrak{z}_{V}$ et $\mathfrak{g}_{1}=
(\mathfrak{r}_{\cur{V}}\times\mathfrak{so}(F))
\oplus\mathfrak{z}_{V}/\mathfrak{q}$, l'idéal
$\mathfrak{z}_{V}/\mathfrak{q}$ étant unipotent et central. Alors, on
a
$\ind\mathfrak{g}_{1}=\ind\mathfrak{r}_{\cur{V}}+\ind\mathfrak{so}(F)+1$,
tandis que $\mathfrak{p}_{\cur{V}}/(\mathfrak{h}+\mathfrak{n}_{V})=
\mathfrak{q}_{\cur{V}}/\mathfrak{r}_{\cur{V}}$. Cependant, il résulte
de l'assertion (ii) du lemme \ref{lem3a'1} que
$\dim\mathfrak{q}_{\cur{V}}/\mathfrak{r}_{\cur{V}}=
\dim\wedge^{2}V^{*}=\dim\mathfrak{z}_{V}$. Appliquant alors le lemme
\ref{lem3a1} (i), on trouve que $\ind\mathfrak{p}_{\cur{V}}=
\ind\mathfrak{r}_{\cur{V}}+\ind\mathfrak{so}(F)$.  Compte tenu de la
proposition \ref{pr3b1} a) (i), la formule \ref{eq3h1} pour l'indice
de $\mathfrak{p}_{\cur{V}}$ est claire dans ce cas.

D'autre part, il résulte du lemme \ref{lem3a1} (ii) que les
sous-algèbres de Cartan-Duflo de $\mathfrak{p}_{\cur{V}}$ et
$\mathfrak{g}_{1}$ ont même dimension, autrement dit que
$\rang\mathfrak{p}_{\cur{V}}=
\rang(\mathfrak{r}_{\cur{V}}\times\mathfrak{so}(F))$. Par suite,
les dimensions des radicaux unipotents d'un 
stabilisateur générique de la représentation coadjointe de
$\p_{\cur{V}}$ ou de $\mathfrak{r}_{\cur{V}}$ sont les mêmes et
égales à $h(\cur{V})$, d'après la proposition \ref{pr3b1} a) (ii).
Il est alors clair que $\mathfrak{p}_{\cur{V}}$ est
quasi-réductive si et seulement si le drapeau $\cur{V}$ vérifie la
propriété $\cur{P}$. D'où le théorème dans ce cas.

\subsection{}\label{3i'}
Dans ce numéro,on suppose que $\dim V=1$ et on reprend les notations
du numéro \ref{3d'}. On a $h(\cur{V}')=0$. Alors $\mathfrak{g}=
\mathfrak{p}_{V}=(\mathbb{K}\times\mathfrak{so}(F))\oplus
F$ et, comme on l'a vu, l'idéal $\mathfrak{n}$ est abélien et son dual
s'identifient à $F$. Par généricité de $g$, on peut supposer que $n$
est un vecteur non isotrope de $F$. Si $F_{0}$ désigne l'orthogonal de
$n$ dans $F$, il suit des relations \ref{eq3e'1} et
\ref{eq3e'2} que l'on a
$\mathfrak{h}+\mathfrak{n}=\mathfrak{h}=\mathfrak{so}(F_{0})\oplus F$ et
$\mathfrak{g}_{1}=\mathfrak{so}(F_{0})\oplus F/F_{0}$. On en déduit tout
d'abord que $\dim\mathfrak{g}/(\mathfrak{h}+\mathfrak{n})=\dim F$ et
ensuite, compte tenu du lemme \ref{lem3a1}, que d'une part
$\ind\mathfrak{g}=\ind\mathfrak{so}(F_{0})=[\frac{1}{2}\dim F_{0}]
=[\frac{q-1}{2}]$ et d'autre part
$\rang\mathfrak{g}=\ind\mathfrak{so}(F_{0})$. Le théorème est donc
démontré dans le cas considéré.

\subsection{}\label{3i''}
Dans ce numéro, on suppose que $r=2p+1$ est impair, avec $p>0$ et que
$F$ est nul, de sorte que $q=2r$. On reprend les notations du numéro
\ref{3f'}. Par généricité de $g$, on peut supposer que $n=\xi$ est générique
relativement au drapeau $\cur{V}$, de sorte qu'en particulier,
$n=\xi\in\mathfrak{n}{_{V}^{*}
}'$. Alors, compte tenu des relations
\ref{eq3f'1} et \ref{eq3f'2}, on a
$\mathfrak{h}=\mathfrak{r}_{\cur{V}}\oplus\mathfrak{n}_{V}$ et
$\mathfrak{g}_{1}=
\mathfrak{r}_{\cur{V}}\oplus\mathfrak{n}_{V}/\mathfrak{q}$, l'idéal
$\mathfrak{n}_{V}/\mathfrak{q}$ étant unipotent et central. Le même
raisonnement qu'au numéro \ref{3i'} montre alors que les algèbres de
Lie $\mathfrak{p}_{\cur{V}}$ et $\mathfrak{r}_{\cur{V}}$ ont même
indice et même rang. Par suite, le théorème dans le cas considéré est
conséquence de la partie b) de la proposition \ref{pr3b1}.

\subsection{}\label{3j}
Désormais, on suppose que $r=2p+1$ est impair et $F$ non nul. On
reprend les notations du numéro \ref{3g}. Par généricité de $g$, on peut
supposer que $n$ est un élément de $\mathfrak{n}{_{V}^{*}
}'$ et que la
forme $\xi$ est générique relativement au drapeau $\cur{V}$. On
choisit les vecteurs $e_{1},\ldots,e_{s},f_{1},\ldots,f_{s}$ du numéro
\ref{3e} de telle sorte que la base $e_{1},\ldots,e_{r}$ de $V$ soit adaptée
au drapeau $\cur{V}$ et que $\xi=\sum_{1\leq i\leq
p}f_{2i-1}\wedge_{B}f_{2i}$. Alors, quitte à translater par un élément
de $\mathbold{N}_{V}$, on peut supposer de plus que
$n=n_{\xi}+v_{r}\wedge_{B}f_{r}$, avec $v_{r}$ un élément non isotrope
de $F$.

On pose $r_{i}=\dim V_{i}$, $1\leq i\leq t$. Alors la suite de
sous-espaces $\{\{0\}=V_{0}\subsetneq V_{1}\subsetneq\cdots\subsetneq
V_{t-1}\}$, complétée par $V'$ si $V_{t-1}\subsetneq V'$, est un
drapeau $\cur{V}'$ de $V'$ générique relativement à la forme
symplectique $\xi_{\mid V'}$.  On désigne par $\mathfrak{t}'_{V}$ le
sous-espace de $\mathfrak{t}_{V'}$ constitué des éléments
$\tilde{T}(c)$, avec $c\in\mathbb{K}^{2p}$ tel que $c_{i}=0$, $1\leq
i\leq r_{t-1}$. Alors, compte tenu de \ref{eq3g1} et \ref{eq3g2}, on a
\begin{align*}
\mathfrak{h} &=(\mathfrak{so}(F_{0})\times
\mathfrak{r}_{\cur{V}'})\oplus(\mathfrak{t}'_{V}\oplus\mathfrak{n}_{V}(n)),\\
\mathfrak{g}_{1}&=\mathfrak{so}(F_{0})\times
(\mathfrak{r}_{\cur{V}'}\oplus\mathfrak{t}'_{V}
\oplus\mathfrak{n}_{V}(n)/\mathfrak{q}),
\end{align*}
avec
\begin{equation*}
\mathfrak{n}_{V}(n)=F\wedge_{B}e_{r}\oplus\mathfrak{z}_{V}.
\end{equation*}
Par ailleurs, $\mathfrak{t}'_{V}\oplus\mathfrak{n}_{V}(n)$ est un
idéal unipotent de l'algèbre de Lie $\mathfrak{h}$.
%$\mathfrak{r}_{\cur{V}'}\oplus\mathfrak{t}'_{V}\oplus\mathfrak{n}_{V}(n)$.

Remarquons que $\mathfrak{r}_{\cur{V}}=\mathfrak{q}_{\cur{V}}(\xi)$
est la sous-algèbre de Lie de $\mathfrak{gl}(V)$ constituée des
éléments de la forme $X_{A,a,c}$ avec $A$ la matrice dans la base
$e_{1},\ldots,e_{2p}$ d'un élément de $\mathfrak{r}_{\cur{V}'}$, $c$
un élément de $\mathbb{K}^{2p}$ tel que $c_{i}=0$, $1\leq i\leq
r_{t-1}$ et $a\in\mathbb{K}$. Par suite, si on désigne par
$\mathfrak{q}_{\cur{V}}(\xi)^{0}$ la sous-algèbre de Lie de
$\mathfrak{q}_{\cur{V}}(\xi)$ constituée des matrices $X_{A,0,c}$, on
a
\begin{equation*}
\mathfrak{p}_{\cur{V}}(n)+\mathfrak{n}_{V}=
\mathfrak{q}_{\cur{V}}(\xi)^{0}\oplus\mathfrak{so}(F_{0})
\oplus\mathfrak{n}_{V}
\end{equation*}
tandis que, d'après l'assertion (ii) du lemme \ref{lem3a'1}
$\dim\mathfrak{q}_{\cur{V}}/\mathfrak{r}_{\cur{V}}=
\dim\wedge^{2}V^{*}$,
si bien que
\begin{align*}
\dim\mathfrak{g}/(\mathfrak{h}+\mathfrak{n})&=
\dim\mathfrak{p}_{\cur{V}}/(\mathfrak{p}_{\cur{V}}(n)+\mathfrak{n}_{V})\\
&=
\dim\mathfrak{q}_{\cur{V}}/\mathfrak{q}_{\cur{V}}(\xi)^{0}
+\dim\mathfrak{so}(F)/\mathfrak{so}(F_{0})\\
&=
\dim\mathfrak{q}_{\cur{V}}/\mathfrak{q}_{\cur{V}}(\xi)+1
+\dim\mathfrak{so}(F)/\mathfrak{so}(F_{0})\\
&=\dim\wedge^{2}V+\dim F\\
&=\dim\mathfrak{n}_{V}(n).
\end{align*}

Alors, désignant par $\mathfrak{s}_{\cur{V}'}$ l'algèbre de Lie
$\mathfrak{r}_{\cur{V}'}\oplus\mathfrak{t}'_{V}
\oplus\mathfrak{n}_{V}(n)/\mathfrak{q}$ et
appliquant le lemme \ref{lem3a1}, on voit que
\begin{equation}\label{eq3j1}
\ind\mathfrak{p}_{\cur{V}}=\ind\mathfrak{so}(F_{0})+
\ind\mathfrak{s}_{\cur{V}'}-1
\end{equation}
et que
\begin{equation}\label{eq3j2}
\rang\mathfrak{p}_{\cur{V}}=\rang\mathfrak{so}(F_{0})+
\rang\mathfrak{s}_{\cur{V}'}.
%\text{$\mathfrak{p}_{\cur{V}}$ est quasi-réductive si et seulement si
%$\rang\mathfrak{s}_{\cur{V}'}=\ind\mathfrak{s}_{\cur{V}'}-1$.}
\end{equation}

Avant d'étudier l'indice et le rang de $\mathfrak{s}_{\cur{V}'}$,
commençons par décrire son idéal unipotent
$\mathfrak{u}_{\cur{V}'}=\mathfrak{t}'_{V}
\oplus\mathfrak{n}_{V}(n)/\mathfrak{q}$. Remarquons que l'application
$\tilde{T}(c)\mapsto c^{*}$ permet d'identifier $\mathfrak{t}'_{V}$
avec le sous-espace $V_{t-1}^{\perp}$ orthogonal de $V_{t-1}$ dans
$V'^{*}$, qui est laissé invariant par $\mathfrak{r}_{\cur{V}'}$. Par
suite, l'application \ref{eq3g3} permet d'identifier l'algèbre de Lie
$\mathfrak{u}_{\cur{V}'}$ à la sous-algèbre
$V_{t-1}^{\perp}\oplus\mathbb{K}Z$ de l'algèbre de Lie de Heisenberg
$\mathfrak{h}_{V'^{*}}$ et l'algèbre de Lie $\mathfrak{s}_{\cur{V}'}$
est une sous-algèbre de Lie de $\mathfrak{s}_{V'}$, produit semi-direct
de $\mathfrak{r}_{\cur{V}'}$ par
$\mathfrak{u}_{\cur{V}'}=V_{t-1}^{\perp}\oplus\mathbb{K}Z$.

\subsection{}\label{3k}
Dans ce numéro, on suppose que $\dim V_{t-1}$ est pair, de sorte que
$h(\cur{V}')=h(\cur{V})$. Alors on a $V'=V_{t-1}\oplus
V_{t-1}^{\perp_{\xi}}$, le sous-espace $V_{t-1}^{\perp{\xi}}$ de $V'$
étant symplectique, et l'isomorphisme $\tilde{\xi}$ de $V'$ sur
$V'^{*}$ du numéro \ref{3g} envoie $V_{t-1}^{\perp_{\xi}}$ sur
$V_{t-1}^{\perp}$. De plus, si $\cur{V}''$ désigne le drapeau
$\{\{0\}\subsetneq V_{1}\subsetneq\cdots\subsetneq V_{t-2}\subsetneq
V_{t-1}\}$ de $V_{t-1}$, on a $\mathfrak{r}_{\cur{V}'}=
\mathfrak{r}_{\cur{V}''}\times\mathfrak{sp}(V_{t-1}^{\perp})$. D'autre
part, l'algèbre de Lie $\mathfrak{u}_{\cur{V}'}$ n'est autre que
l'algèbre de Heisenberg $\mathfrak{h}_{V_{t-1}^{\perp}}$, vue comme
sous-algèbre de Lie de $\mathfrak{h}_{V'^{*}}$. On voit alors que l'on
a
\begin{equation*}
\mathfrak{s}_{\cur{V}'}=\mathfrak{r}_{\cur{V}''}
\times(\mathfrak{sp}(V_{t-1}^{\perp})\oplus\mathfrak{h}_{V_{t-1}^{\perp}}).
\end{equation*}
D'après les résultats du numéro \ref{1g3}, on a
$\ind\mathfrak{sp}(V_{t-1}^{\perp})\oplus\mathfrak{h}_{V_{t-1}^{\perp}}=
1+\frac{1}{2}\dim V_{t-1}^{\perp}$. Compte tenu de la relation
\ref{eq3j1}, il vient donc
\begin{equation*}
\ind\mathfrak{p}_{\cur{V}}=\ind\mathfrak{so}(F_{0})+
\ind\mathfrak{r}_{\cur{V}''}+\frac{1}{2}\dim V_{t-1}^{\perp}.
\end{equation*}
La formule \ref{eq3h1} suit alors de ce que, d'une part
$\ind\mathfrak{so}(F_{0})=[\frac{1}{2}\dim F_{0}]=[\frac{q-1}{2}]-r$ et,
d'autre part, d'après l'assertion (i) de la proposition \ref{pr3b1},
$\ind\mathfrak{r}_{\cur{V}''}+\frac{1}{2}\dim V_{t-1}^{\perp}=
\sum_{i=1}^{t}[\frac{1}{2}(\dim V_{i}-\dim V_{i-1})]$.

D'autre part, compte tenu des résultats du numéro \ref{1g3}, on a
$\rang\mathfrak{s}_{\cur{V}'}= \rang\mathfrak{r}_{\cur{V}''}+
\rang(\mathfrak{sp}(V_{t-1}^{\perp})\oplus\mathfrak{h}_{V_{t-1}^{\perp}})
=\rang\mathfrak{r}_{\cur{V}''}+\frac{1}{2}\dim V_{t-1}^{\perp}$.
Il suit alors de la relation \ref{eq3j2} que les dimensions des
radicaux unipotents d'un stabilisateur générique de
la représentation coadjointe de $\p_{\cur{V}}$ ou de
$\mathfrak{r}_{\cur{V}''}$ sont les mêmes et égales à
$h(\cur{V}'')$, d'après la proposition \ref{pr3b1} a) (ii). Or on
a vu que $h(\cur{V}')=h(\cur{V})$ et il est clair que
$h(\cur{V}'')=h(\cur{V}')$. Les assertions (ii) et (iii) du
théorème sont démontrées dans ce cas.

\subsection{}\label{3l}
Dans ce numéro, on suppose que $\dim V_{t-1}$ est impair. Posons
$V''=V_{t-1}\oplus\mathbb{K}e_{r_{t-1}+1}$ de sorte que $V''$ est un
sous-espace symplectique de $V'$ et que l'on a la décomposition
orthogonale symplectique $V'=V''\oplus
V''^{\perp_{\xi}}$. Alors $\mathbb{K}e^{*}_{r_{t-1}+1}$ est le noyau
de la restriction de $\xi$ à $V_{t-1}^{\perp}$ et, si $V''^{\perp}$
désigne l'orthogonal de $V''$ dans $V'^{*}$, on a
$\mathfrak{u}_{\cur{V}'}=
\mathfrak{h}_{V''^{\perp}}\times\mathbb{K}e^{*}_{r_{t-1}+1}$.

D'autre part, la décomposition en somme directe
$\mathfrak{u}_{\cur{V}'}=V_{t-1}^{\perp}\oplus\mathbb{K}Z$ de la fin
du numéro \ref{3j}, induit la décomposition de l'espace dual
$\mathfrak{u}^{*}_{\cur{V}'}=V_{t-1}^{\perp_{\xi}}\oplus\mathbb{K}Z^{*}$,
où $Z^{*}$ désigne la forme linéaire nulle sur $V_{t-1}^{\perp}$ et
telle que $\langle Z^{*},Z\rangle=1$.

On pose $\dim V''^{\perp}=2m$ et on désigne par $\cur{V}''$ le drapeau
$\{\{0\}\subsetneq V_{1}\subsetneq\cdots\subsetneq V_{t-1}\subsetneq
V''\}$ de $V''$. Alors, les éléments de $\mathfrak{r}_{\cur{V}'}$ sont
les endomorphismes de $V'$ qui, si on les identifie avec leur matrice
dans la base $e_{1},\ldots,e_{2p}$, sont de la forme
\begin{equation*}
X_{A,a,\alpha,y,\beta,D}=
\begin{bmatrix}
A & 0 & J_{r_{t-1}-1}\,^{t}\alpha&0\\
\alpha &a&y&\beta\\
0&0&-a&0\\
0&0&J_{2m}\,^{t}\beta &D
\end{bmatrix}
\end{equation*}
où $\alpha\in\mathbb{K}^{r_{t-1}-1}$ et $\beta\in\mathbb{K}^{2m}$ sont
des vecteurs lignes, $a,y\in\mathbb{K}$, $A$ est une matrice carrée
d'ordre $r_{t-1}-1$ et $D\in\mathfrak{sp}_{2m}$, $A$, $\alpha$, $y$ et
$a$ étant tels que la
sous-matrice
\begin{equation*}
\begin{bmatrix}
A & 0 & J_{r_{t-1}-1}\,^{t}\alpha\\
\alpha &a&y\\
0&0&-a
\end{bmatrix}
\end{equation*}
soit la matrice dans la base $e_{1},\ldots,e_{r_{t-1}+1}$ d'un élément
de $\mathfrak{r}_{\cur{V}''}$.

Par ailleurs, on identifie $\mathfrak{sp}(V''^{\perp_{\xi}})$ avec la
sous-algèbre de Lie de $\mathfrak{r}_{\cur{V}''}$ constituée des
$X_{0,0,0,0,0,D}$, $D\in\mathfrak{sp}_{2m}$, au moyen de l'application qui
à l'élément de $\mathfrak{sp}(V''^{\perp_{\xi}})$ de matrice $D$ dans la
base $e_{r_{t-1}+2},\ldots e_{r_{t-1}+2m}$ fait correspondre
$X_{0,0,0,0,0,D}$.

Nous devons calculer l'indice et le rang de
$\mathfrak{g}=\mathfrak{s}_{\cur{V}'}$. Pour ce faire, on choisit $g$
une forme linéaire fortement régulière sur
$\mathfrak{s}_{\cur{V}'}$. Soit $n$ la restriction de $g$ à l'idéal
unipotent $\mathfrak{n}=\mathfrak{u}_{\cur{V}'}$ de
$\mathfrak{g}$. Par généricité de $g$ et quitte à translater par un
élément de $\mathbold{U}_{\cur{V}'}$ le sous-groupe unipotent de
$\mathrm{GL}(V)$ d'algèbre de Lie $\mathfrak{u}_{\cur{V}'}$, on peut
supposer que $n=\tau e_{r_{t-1}+1}+Z^{*}$, avec $\tau\neq 0$.

On commence par déterminer $\mathfrak{h}=\mathfrak{g}(n)$. Tout
d'abord si $c=\sum_{i=1}^{2m}c_{i}e^{*}_{r_{t-1}+1+i}\in²V''^{\perp}$,
on a
\begin{equation*}
c.n=\sum_{i=1}^{m}(c_{2i}e_{r_{t-1}+2i}-c_{2i-1}e_{r_{t-1}+2i+1})
\end{equation*}
tandis que, si $X_{A,a,\alpha,y,\beta,D}\in\mathfrak{r}_{\cur{V}'}$,
\begin{equation*}
X_{A,a,\alpha,y,\beta,D}.n=-\tau a e_{r_{t-1}+1}+
\tau\sum_{i=1}^{m}(\beta_{2i}e_{r_{t-1}+2i}-\beta_{2i-1}e_{r_{t-1}+2i+1}).
\end{equation*}
Il s'ensuit que $X_{A,a,\alpha,y,\beta,D}+c$ est dans $\mathfrak{h}$
si et seulement si $a=0$ et $c_{i}=-\tau\beta_{i}$, $1\leq i\leq2m$.

Si $\beta\in\mathbb{K}^{2m}$, on pose
$V(\beta)=-\tau\sum_{i=1}^{2m}\beta_{i}e^{*}_{r_{t-1}+1+i}$ et
$\tilde{V}(\beta)=X_{0,0,0,0,\beta,0}+V(\beta)$. On désigne par
$\mathfrak{v}_{\cur{V}'}$ le sous-espace de $\mathfrak{g}$ constitué
des $\tilde{V}(\beta)$, $\beta\in\mathbb{K}^{2m}$.

Remarquons que, pour $\beta,\beta'\in\mathbb{K}^{2m}$, on a
\begin{equation*}
[\tilde{V}(\beta),\tilde{V}(\beta')]=
\beta J_{2m}\,{t}\beta' Z_{1},
\end{equation*}
avec $Z_{1}=Y+\tau^{2}Z$ et $Y=X_{0,0,0,2,0,0}$. Ainsi,
$\mathfrak{v}_{\cur{V}'}\oplus\mathbb{K}Z_{1}$ est une sous-algèbre de
Lie unipotente de $\mathfrak{h}$ qui s'identifie, au moyen de
l'application
\begin{equation}\label{eq3l1}
\phi : \tilde{V}(\beta)+zZ_{1}\mapsto
\sum_{i=1}^{2m}\beta_{i}e^{*}_{r_{t-1}+1+i}+z
Z
\end{equation}
à l'algèbre de Lie de Heisenberg $\mathfrak{h}_{V''^{\perp}}$
construite sur le sous-espace symplectique $V''^{\perp}$ de $V'^{*}$.

Comme on l'a vu plus haut, $\mathfrak{r}_{\cur{V}''}$ est la
sous-algèbre de Lie de $\mathfrak{r}_{\cur{V}'}$ constituée des
éléments de la forme $X_{A,a,\alpha,y,0,0}$. L'application
$\chi:X_{A,a,\alpha,y,0,0}\mapsto a$ est un caractère de
$\mathfrak{r}_{\cur{V}''}$ dont on note $\mathfrak{r}^{0}_{\cur{V}''}$
le noyau. Alors, on constate que $\mathfrak{r}^{0}_{\cur{V}''}$ et
$\mathfrak{sp}(V''^{\perp_{\xi}})$ sont des sous-algèbres de Lie qui
commutent dans $\mathfrak{r}_{\cur{V}'}$, tandis que si
$X_{A,0,\alpha,y,0,D}\in\mathfrak{r}^{0}_{\cur{V}''}
\times\mathfrak{sp}(V''^{\perp_{\xi}})$, %l'élément $D$ de
%$\mathfrak{sp}(V''^{\perp_{\xi}})$ étant identifié à sa matrice dans la
%base $e_{r_{t-1}+2},\ldots e_{r_{t-1}+2m}$ de $V''^{\perp_{\xi}}$ et
%$\beta\in\mathbb{K}^{2m}$ étant un vecteur ligne,
on a
\begin{equation*}
[X_{A,0,\alpha,y,0,D},\tilde{V}(\beta)]=\tilde{V}(-\beta D).
\end{equation*}
On voit alors que l'application $\phi$ définie par la formule
\ref{eq3l1} est un isomorphisme
$\mathfrak{r}^{0}_{\cur{V}''}
\times\mathfrak{sp}(V''^{\perp_{\xi}})$-équivariant d'algèbres de Lie,
l'action de $\mathfrak{r}^{0}_{\cur{V}''}$ étant des deux côtés
triviale. Il suit de ces considérations que l'on a
\begin{equation*}
\mathfrak{h}=\mathfrak{r}^{0}_{\cur{V}''}\times
(\mathfrak{sp}(V''^{\perp_{\xi}})\oplus\mathfrak{h}_{V''^{\perp}})
\times\mathbb{K}e_{r_{t-1}+1}.
\end{equation*}
D'autre part, on a $\mathfrak{h}+\mathfrak{n}=
\mathfrak{r}^{0}_{\cur{V}''}\oplus\mathfrak{n}$ si bien que
$\dim\mathfrak{g}/(\mathfrak{h}+\mathfrak{n})$=1. Appliquant le lemme
\ref{lem3a1} et les résultats du numéro \ref{1g3}, on voit alors que
\begin{align}
\ind\mathfrak{s}_{\cur{V}'} &=\ind\mathfrak{r}^{0}_{\cur{V}''}
+\frac{1}{2}\dim V''^{\perp_{\xi}}+1\label{eq3l2}\\
\rang\mathfrak{s}_{\cur{V}'} &=\rang\mathfrak{r}^{0}_{\cur{V}''}
+\frac{1}{2}\dim V''^{\perp_{\xi}}\label{eq3l3}
\end{align}

%\subsection{}\label{3m}
Rassemblant ces formules avec les formules \ref{eq3b''''1} et
\ref{eq3b''''2}, on obtient
\begin{align*}
\ind\mathfrak{s}_{\cur{V}'} &=\ind\mathfrak{r}_{\cur{V}''}
+\frac{1}{2}\dim V''^{\perp_{\xi}}+2\label{}\\
\rang\mathfrak{s}_{\cur{V}'} &=\rang\mathfrak{r}_{\cur{V}''}
+\frac{1}{2}\dim V''^{\perp_{\xi}}\notag.
\end{align*}

La comparaison de ces deux dernières formules montre d'une part,
compte tenu de l'assertion \ref{eq3j2}, que dans notre situation
l'algèbre de Lie $\mathfrak{p}_{\cur{V}}$ n'est pas quasi-réductive et
d'autre part que la dimension du radical unipotent des stabilisateurs
génériques de la représentation coadjointe de $\p_{\cur{V}}$ est égale
à $\ind\mathfrak{r}_{\cur{V}''}-\rang\mathfrak{r}_{\cur{V}''}+1$ soit
$h(\cur{V}'')+1$, d'après la proposition \ref{pr3b1}. Comme les
sous-espaces $V_{t-1}$ et $V_{t}$ sont tous deux de dimension impaire,
ceci achève la preuve des assertions (ii) et (iii) du théorème.

Compte tenu de la relation \ref{eq3j1} et du fait que
$\ind\mathfrak{so}(F_{0})=[\frac{1}{2}\dim F_{0}]=[\frac{q-1}{2}]-r$, la
première de ces formules donne
\begin{equation*}
\ind\mathfrak{p}_{\cur{V}}=[\frac{q-1}{2}]-r+\ind\mathfrak{r}_{\cur{V}''}
+\frac{1}{2}\dim V''^{\perp_{\xi}}+1,
\end{equation*}
relation qui, du fait de l'assertion (i) de la proposition \ref{pr3b1}
et de ce que $\dim V_{t}-\dim V_{t-1}=2+\dim V''^{\perp_{\xi}}$, n'est
autre dans le cas d'espèce que la formule \ref{eq3h1}.
\end{dem}

\subsection{Quasi-réductivité des sous-algèbres paraboliques de
  $\mathfrak{so}(q,\mathbb{K})$ et diagrammes de Dynkin.}\label{3n}

Dans ce numéro et les suivants, nous reformulons la
caractérisation des sous-algèbres paraboliques de
$\mathfrak{so}(q,\mathbb{K})$ en terme de système de racines et de
diagrammes de Dynkin. Comme pour $3\leq q\leq 6$,
les algèbres $\mathfrak{so}(q,\mathbb{K})$ sont de type $A_{1}$,
$A_{1}\times A_{1}$, $C_{2}$ ou $A_{3}$, leurs sous-algèbres
paraboliques sont toutes quasi-réductives d'après
\cite{panyushev-2005}. Dans la suite, nous supposerons donc que
$q\geq7$.

Nous rappelons quelques notations relatives aux algèbres de Lie
semi-simples. Soit $\s$ une algèbre de Lie semi-simple et
$\mathfrak{b}\subset\s$ une sous-algèbre de Borel de radical
nilpotent $\n$. Ces données déterminent un ensemble de racines
simples $\Pi$ dans $(\mathfrak{b}/\n)^{*}$ ainsi que le diagramme
de  Dynkin correspondant. Nous indexons les
sous-algèbres paraboliques standard (c'est à dire, celles
contenant $\mathfrak{b}$) par les sous-ensembles $I\subset\Pi$ en
décidant que $I$ est l'ensemble des racines simples \emph{qui ne
sont pas} contenues dans la partie réductive de la sous-algèbre
parabolique $\p_{I}$. Par exemple, les sous-algèbres paraboliques
propres maximales sont les $\p_{\{\alpha\}}$, $\alpha\in\Pi$. On a
$\p_{I}=\cap_{\alpha\in I}\p_{\{\alpha\}}$,
$\p_{\Pi}=\mathfrak{b}$ et $\p_{\emptyset}=\s$.

Nous représentons un sous-ensemble $I\subset\Pi$ et la sous-algèbre
parabolique correspondante $\p_{I}$ en coloriant les racines sur le
diagramme de Dynkin~: en noir, celles qui sont situées dans $I$, en
blanc les autres, c'est à dire celles du facteur de Levi de $\p_{I}$.

\subsection{}\label{3o}
On suppose que le $\mathbb{K}$-espace vectoriel $E$ est de dimension
$q=2n+1$, $n\geq 3$. On se donne un drapeau $\cur{F}=\{F_{0}\subset
F_{1}\subset\cdots\subset F_{n}\}$ de sous-espaces isotropes, avec
$\dim F_{i}=i$. Deux tels drapeaux sont conjugués sous l'action de
$\mathrm{SO}(E)$. Le stabilisateur de $\cur{F}$ est une sous-algèbre
de Borel $\mathfrak{b}\subset\mathfrak{so}(E)$.

Soit $I$ un sous-ensemble de $\{1,\ldots,n\}$. On note $\cur{F}_{I}$
le sous-drapeau de $\cur{F}$ constitué des sous-espaces $F_{i}$, $i\in
I$ et on désigne par $\p_{I}$ son stabilisateur. Les algèbres $\p_{I}$
sont les sous-algèbres paraboliques standard de $\mathfrak{so}(E)$. On
indexe les racines simples $\Pi=\{\alpha_{1},\ldots\alpha_{n}\}$ de
telle manière que $\p_{\alpha_{i}}$ soit égal à $\p_{\{i\}}$, le
stabilisateur de $F_{i}$. On note $\alpha_{0}$ la plus basse
racine. Pour $n\geq 3$, on obtient le diagramme de Dynkin étendu suivant~:
\begin{figure}[H] %\label{fig:Bn}
\includegraphics{SOn-quasired.6}
\caption{$B_n$}\label{fig:Bn}
\end{figure}

Par exemple, le diagramme colorié suivant représente la sous-algèbre
parabolique \flqq minimale\frqq~ $\p_{\{1,3\}}$ de
$\mathfrak{so}(7,\mathbb{K})$, stabilisateur du drapeau
$\{F_{1},F_{3}\}$, obtenue à partir de $\mathfrak{b}$ en ajoutant le
sous-espace radiciel correspondant à la racine $-\alpha_{2}$~:
\begin{figure}[H]
\includegraphics{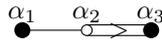}
\caption{$B_3$: stabilisateur de $F_{1,3}$}\label{fig:B3-13}
\end{figure}

\subsection{}\label{3p}
Nous utilisons des notations semblables à celles des numéros
précédents. Ici $E$ est de dimension $q=2n$, avec $n\geq 4$.  On se
donne un drapeau $\cur{F}=\{F_{0}\subset F_{1}\subset\cdots\subset
F_{n-1}\}$ de sous-espaces isotropes, avec $\dim F_{i}=i$. Deux tels
drapeaux sont conjugués sous l'action de $\mathrm{SO}(E)$. L'espace
$F_{n-1}$ est contenu dans exactement deux sous-espaces isotropes de
dimension $n$, que nous notons $F_{n}^{+}$ et $F_{n}^{-}$. Ils
vérifient $F_{n-1}=F_{n}^{+}\cap F_{n}^{-}$. On désigne par
$\cur{F}^{+}$ le drapeau $\cur{F}^{+}=\{F_{0}\subset
F_{1}\subset\cdots\subset F_{n-1}\subset F_{n}^{+}\}$ et de manière
similaire pour $\cur{F}^{-}$. Le stabilisateur de $\cur{F}$ est une
sous-algèbre de Borel $\mathfrak{b}\subset\mathfrak{so}(E)$. Le
stabilisateur $B$ de $\cur{F}$ dans $\mathrm{SO}(E)$ stabilise
également $\cur{F}^{+}$ et $\cur{F}^{-}$. Par contre, il existe des
éléments du stabilisateur dans $\mathrm{O}(E)$ de $\cur{F}$ qui
permutent les drapeaux $\cur{F}^{+}$ et $\cur{F}^{-}$.

Les sous-algèbres paraboliques maximales sont les stabilisateurs des
espaces $F_{i}$, $1\leq i\leq n-1$ et des espaces $F_{n}^{+}$ et
$F_{n}^{-}$. Nous indexons en conséquence les racines simples
$\Pi=\{\alpha_{1},\ldots\alpha_{n-1},\alpha_{n}^{+},\alpha_{n}^{-}\}$
et nous désignons par $\alpha_{0}$ la plus basse racine. Pour $n\geq
4$, on obtient le diagramme de Dynkin étendu suivant~:
\begin{figure}[H]
\includegraphics{SOn-quasired.9}
\caption{$D_n$} \label{fig:Dn}
\end{figure}

On introduit deux symboles $n^{+}$ et $n^{-}$. Soit $I$ un
sous-ensemble de $\{1,\ldots,n-1,n^{+},n^{-}\}$, que l'on identifie à
un sous-ensemble de $\Pi$ et soit $\p_{I}$ la sous-algèbre parabolique
standard correspondante de $\mathfrak{so}(E)$~: on a
$\p_{\{i\}}=\p_{\alpha_{i}}$ pour $1\leq i\leq n-1$,
$\p_{\{n^{+}\}}=\p_{\alpha_{n}^{+}}$ et $\p_{\{n^{-}\}}=\p_{\alpha_{n}^{-}}$.

Le stabilisateur de $F_{n-1}$ n'est pas une sous-algèbre parabolique
maximale~: c'est la sous-algèbre parabolique $\p_{\{n^{+},n^{-}\}}$,
intersection des stabilisateurs de $F_{n}^{+}$ et $F_{n}^{-}$. Plus
généralement, le stabilisateur d'un drapeau contenant $F_{n-1}$,
stabilise également chacun des deux drapeaux obtenus en ajoutant
$F_{n}^{+}$ et $F_{n}^{-}$.

\begin{figure}[H]%\label{fig:D4}
\includegraphics{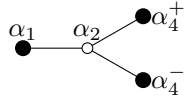}
\caption{$\p_{\{1,3\}}=\p_{\{1,4^+,4^-\}}$, stabilisateur du drapeau
  $\{F_1\subset F_3\}$. Ici $I^o=\{1,3\}$ } \label{fig:D4}
\end{figure}

\begin{figure}[H]%\label{fig:D5}
\includegraphics{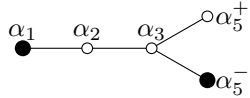}
\caption{$\p_{\{1,5^-\} }$, stabilisateur du drapeau $\{F_1\subset
  F_5^-\}$. Ici $I^o=\{1 \}$} \label{fig:D5}
\end{figure}

Soit $I\subset\{1,\ldots,n-2,n^{+},n^{-}\}$. On désigne par
$\widetilde{I}\subset\{1,\ldots,n-1,n^{+},n^{-}\}$ l'ensemble défini par~:
$$
\tilde I = (I\cap\{1,\dots,n-2 \})\cup \{n-1\}
$$
si $\{n^+, n^-\}\subset I$, et
$$
\tilde I = I
$$ si $\{n^{+},n^{-}\}$ n'est pas contenu dans $I$. L'application
$I\mapsto\widetilde{I}$ est injective. Son image est l'ensemble des
parties de $\{1,\ldots,n-1,n^{+},n^{-}\}$ qui ne contiennent pas
$\{n^{+},n^{-}\}$. Par suite, $\cur{F}_{\widetilde{I}}:=\{F_i, i\in
\widetilde I\}$ est un drapeau et la sous-algèbre parabolique $\p_{I}$
en est le stabilisateur. Nous illustrons ceci à l'aide des figures
\ref{fig:D4} et \ref{fig:D5} ci-après.

Nous aurons également besoin de l'ensemble
$I^{0}:=\widetilde{I}\cap\{1,\dots,n-2 ,n-1\}$. Ainsi,
$$
 I^o = (I\cap\{1,\dots,n-2 \})\cup \{n-1\}
$$
si $\{n^+, n^-\}\subset I$, et
$$
 I^o = I\cap\{1,\dots,n-2 \}
$$
si $\{n^+, n^-\}$ n'est pas contenu dans $I$.

\subsection{}\label{3q}
On considère $\mathfrak{so}(q,\mathbb{K})$ avec $q\geq7$. On utilise le
drapeau $\cur{F}=\{F_{0}\subset F_{1}\cdots\subset F_{n}\}$ (si
$q=2n+1$) ou $\cur{F}=\{F_{0}\subset F_{1}\cdots\subset F_{n-1}\}$ (si
$q=2n$). Le choix de $\cur{F}$ correspond au choix d'une sous-algèbre
de Borel $\mathfrak{b}\subset\mathfrak{so}(q,\mathbb{K})$ et permet de
définir la notion de sous-algèbres paraboliques standard, c'est à dire
de sous-algèbres de $\mathfrak{so}(q,\mathbb{K})$ contenant
$\mathfrak{b}$.

L'assertion (iii) du théorème \ref{theo3h1} peut s'énoncer ainsi~:
\begin{theo}\label{theo3q1}
Les sous-algèbres paraboliques standard non quasi-réductives de
$\mathfrak{so}(q,\mathbb{K})$ sont les sous-algèbres paraboliques
standard $\p$ pour lesquelles il existe deux éléments de dimension
impaire $F_{i}\subsetneq F_{j}$ de $\cur{F}$ tels que $F_{i}$ et
$F_{j}$ soient stables sous l'action de $\p$ et qu'aucun des
sous-espaces $F_{k}$ avec $i<k<j$ ne soit stable sous l'action de
$\p$.
\end{theo}

Nous reformulons ce théorème en utilisant la description des
sous-algèbres paraboliques en termes de sous-ensembles de l'ensemble
$\Pi$ des racines simples, indexées comme expliqué dans les précédents
numéros.
\begin{theo}\label{theo3q2}
Soit $n\geq3$ et $q=2n+1$. Soit $I\subset\{1,\ldots,n\}$. Alors, la
sous-algèbre parabolique $\p_{I}\subset\mathfrak{so}(q,\mathbb{K})$
n'est pas quasi-réductive si et seulement si il existe deux entiers
impairs $j<k$ tels que $j\in I$, $k\in I$ et $p\notin I$, pour
$j<p<k$.
\end{theo}

\begin{theo}\label{theo3q3}
Soit $n\geq4$ et $q=2n$. Soit
$I\subset\{1,\ldots,n-2,n^{+},n^{-}\}$. Alors, la sous-algèbre
  parabolique $\p_{I}\subset\mathfrak{so}(q,\mathbb{K})$ n'est pas
  quasi-réductive si et seulement si il existe deux entiers impairs
  $j<k$ tels que $j\in I^{0}$, $k\in I^{0}$ et $p\notin I^{0}$, pour
  $j<p<k$.
\end{theo}

\begin{rem}\label{rem3q1}
Soit $I\subset\Pi$ et $I'=\Pi\backslash I$. Ce dernier ensemble est
l'ensemble des racines simples d'un facteur de Levi de la sous-algèbre
parabolique $\p_{I}$. On désigne par $h$ le nombre de composantes
connexes de $I'$ de la forme $[i,j]$ avec $i$ et $j$ pairs et $2\leq
j\leq n-1$ (si $q=2n+1$) ou $2\leq j\leq n-2$ (si $q=2n$).

Alors, nous pouvons reformuler le théorème \ref{theo3q1} de la manière
suivante : la sous-algèbre parabolique $\p_{I}$ est quasi-réductive si
et seulement si $h=0$.

En fait, $h$ est égal à la dimension du radical unipotent d'un
stabilisateur générique de la représentation coadjointe de $\p_{I}$
(assertion (ii) du théorème \ref{theo3h1}).
\end{rem}

\begin{rem}\label{rem3q2}
Selon la remarque précédente, les deux exemples de plus petite
dimension de sous-algèbres paraboliques non quasi-réductives de
$\mathfrak{so}(q,\mathbb{K})$ son ceux des figures
\ref{fig:D4} et \ref{fig:D5}. L'exemple de la figure \ref{fig:D4} est
dû à Tauvel et Yu (\cite{tauvel-yu-2004-a}).
\end{rem}

\begin{rem}
Si $I$ ne contient que des entiers pairs (pour $q=2n+1$), ou si
$I^{0}$ ne contient que des entiers pairs (pour $q=2n$), alors les
théorèmes précédents montrent que l'algèbre $\p_{I}$ est
quasi-réductive. Ce résultat est dû à Panyushev
(\cite{panyushev-2005}).
\end{rem}

\begin{rem}
Nos résultats sont implicitement contenus dans les travaux de Dvorsky
(\cite{dvorsky-1995} et \cite{dvorsky-1996}). En effet, le nombre $h$
introduit dans la remarque \ref{rem3q1} est le même que l'entier $h$
défini dans \cite[formule 6.2]{dvorsky-1996} (pour
$\mathfrak{so}(2n+1,\mathbb{K})$) et dans \cite[formule 3.17 et lignes
  suivantes]{dvorsky-1995}.
\end{rem}

\section{Algèbres de Lie non algébriques}\label{sec:suppl}

Dans cette section, on suppose toujours que $\K$ est un corps de
caractéristique nulle. Par contre, $\g$ désigne une algèbre de Lie sur
$\K$ qui n'est plus donnée comme étant l'algèbre de Lie d'un groupe
algébrique. Telles qu'elles ont été formulées, les définitions
\ref{def:typered} de formes de type réductif et \ref{def:quasired}
d'algèbres quasi-réductives ont toujours un sens dans ce contexte plus
large. En fait, l'exemple suivant en est à l'origine~:

\subsection{Exemple.}\label{sec:suppl1}
On suppose que $\K=\R$. Soit $G$ un groupe de Lie connexe d'algèbre de
Lie $\g$ et soit $Z_{G}$ son centre. Soit $g\in\g^{*}$ telle que
$G(g)/Z_{G}$ soit compact. Alors, $g$ est de type réductif.

\subsection{}\label{sec:suppl2}
Dans la suite de cette section, on suppose que $\K$ est un corps
algébriquement clos de caractéristique nulle. On désigne par
$\hat{\mathbold{G}}$ le plus petit sous-groupe algébrique de
$\mathrm{GL}(\g)$ dont l'algèbre de Lie $\hat{\g}$ contient
$\g/\mathfrak{z}_{\g}$. En fait, $\hat{\g}$ est la plus petite
sous-algèbre de Lie algébrique contenant $\g/\mathfrak{z}_{\g}$.

Une autre propriété naturelle que l'on peut considérer pour une forme
linéaire $g$ est que son stabilisateur $\hat{\mathbold{G}}(g)$ dans
$\hat{\mathbold{G}}$ soit réductif. L'exemple du numéro suivant montre
que cette propriété est en général strictement plus forte que celle
d'être de type réductif. Le même exemple montre également qu'une
algèbre de Lie quasi-réductive $\g$ (en fait, même une algèbre de
Frobenius) peut admettre une enveloppe algébrique qui ne soit pas
quasi-réductive.

\subsection{Exemple.}\label{sec:suppl3}
On considère l'algèbre de Lie $\g$ admettant la base $h,x,y,z$ avec
les crochets non nuls $[x,y]=z$, $[h,x]=x+y$, $[h,y]=y$ et
$[h,z]=2z$. Le sous-espace $\h$ dont une base est $x,y,z$, est un
idéal, de centre $\K z$, isomorphe à l'algèbre de Lie de Heisenberg.
Soit $g=z^{*}\in\g^{*}$. On a $\g(g)=0$, de sorte que $\g$ est une
algèbre de Frobenius. En particulier, $\g$ est quasi-réductive.

Maintenant, soit $\hat{\g}$ l'algèbre de Lie admettant la base
$s,u,x,y,z$, telle que $\h$ soit un idéal de $\hat{\g}$, $[s,u]=0$,
$[s,x]=x$, $[s,y]=y$, $[s,z]=2z$, $[u,x]=y$, $[u,y]=0$ et
$[u,z]=0$. On considère $\g$ comme une sous-algèbre de Lie de
$\hat{\g}$ en posant $h=s+u$.  Alors, $\hat{\g}$ est algébrique tandis
que $\g$ ne l'est pas. En fait, $h=s+u$ est la décomposition de
Chevalley de $h$ en ses parties semi-simple et unipotente. Par suite,
$\hat{\g}$ est l'enveloppe algébrique de $\g$.

Le stabilisateur $\hat{\g}(g)$ de $g$ est égal à $\K u$, de sorte que
$\hat{\mathbold{G}}(g)$ est un groupe unipotent de dimension $1$.

Soit $\hat{g}\in\hat{\g}^{*}$ une forme linéaire dont la restriction à
$\g$ est égale à $g$. On a également $\hat{\g}(\hat{g})=\K u$. On en
déduit facilement que $\hat{\g}$ n'est pas quasi-réductive.

%{\color{red} Votre article a JFA 235 n'est pas à sa place dans la
%bibliographie}

\bibliographystyle{smfplain}
\bibliography{biblio}

\providecommand{\bysame}{\leavevmode ---\ }
\providecommand{\og}{``}
\providecommand{\fg}{''}
\providecommand{\smfandname}{et}
\providecommand{\smfedsname}{\'eds.}
\providecommand{\smfedname}{\'ed.}
\providecommand{\smfmastersthesisname}{M\'emoire}
\providecommand{\smfphdthesisname}{Th\`ese}
\begin{thebibliography}{10}

\bibitem{andler-1985}
{\scshape M.~Andler} -- {\og {La formule de Plancherel pour les groupes
  alg\'ebriques complexes unimodulaires}\fg}, \emph{Acta. Math.} \textbf{154}
  (1985), p.~1--104.

\bibitem{anh-1974}
{\scshape N.~H. Anh} -- {\og {Algebraic groups with square-integrable
  representations}\fg}, \emph{Bull. Amer. Math. Soc.} \textbf{80} (1974),
  p.~539--542.

\bibitem{anh-1978}
\bysame , {\og {Classification of unimodular algebraic groups with square
  integrable representations}\fg}, \emph{Acta Mathematica Vietnamica}
  \textbf{3} (1978), no.~2, p.~75--81.

\bibitem{anh-1980}
\bysame , {\og {Classification of connected unimodular Lie groups with discrete
  series}\fg}, \emph{Ann. Inst. Fourier, Grenoble} \textbf{30} (1980), no.~1,
  p.~159--192.

\bibitem{baur-moreau-2009}
{\scshape K.~Baur {\normalfont \smfandname} A.~Moreau} -- {\og {Quasi-reductive
  (bi)parabolic subalgebras in the reductive Lie algebras}\fg},
  arXiv:0812.4275.

\bibitem{borel-harish-chandra-1962}
{\scshape A.~Borel {\normalfont \smfandname} Harish-Chandra} -- {\og
  {Arithmetic subroups of algebraic groups}\fg}, \emph{Annals of Mathematics}
  \textbf{75} (1962), p.~485--535.

\bibitem{burde-2006}
{\scshape D.~Burde} -- {\og {Characteristically nilpotent Lie algebras and
  symplectic structures}\fg}, \emph{Forum Math.} \textbf{18} (2006), no.~5,
  p.~769--787.

\bibitem{charbonnel-1982}
{\scshape J.-Y. Charbonnel} -- {\og {Sur les orbites de la repr\'esentation
  coadjointe.}\fg}, \emph{Compositio Math.} \textbf{46} (1982), p.~273--305.

\bibitem{charbonnel-dixmier-1981}
{\scshape J.-Y. Charbonnel {\normalfont \smfandname} J.~Dixmier} -- {\og {Sur
  la représentation coadjointe des algèbres de Lie résolubles.}\fg}, \emph{J.
  Reine Angew. Math.} \textbf{324} (1981), p.~154--164.

\bibitem{drozd-ip-1992}
{\scshape Y.~A. Drozd} -- {\og Matrix problems, small reduction and
  representations of a class of mixed {L}ie groups\fg}, \emph{Representations
  of algebras and related topics} (Kyoto, 1990), London Math. Soc. Lecture Note
  Ser., no. 168, Cambridge Univ. Press, Cambridge, 1992, p.~225--249.

\bibitem{duflo-1982}
{\scshape M.~Duflo} -- {\og {Th\'{e}orie de Mackey pour les groupes de Lie
  alg\'{e}briques}\fg}, \emph{Acta. Math.} \textbf{149} (1982), p.~153--213.

\bibitem{dvorsky-1995}
{\scshape A.~Dvorsky} -- {\og {Generic representations of parabolic subgroups.
  Case I - {\bf SO(2n,C)} }\fg}, \emph{Comm. Algebra} \textbf{23} (1995),
  no.~4, p.~1369--1402.

\bibitem{dvorsky-1996}
\bysame , {\og {Generic representations of parabolic subgroups. Case II - {\bf
  SO(2n+1,C)}}\fg}, \emph{Comm. Algebra} \textbf{24} (1996), no.~1,
  p.~259--288.

\bibitem{dvorsky-2003}
\bysame , {\og {Index of parabolic and seaweed subalgebras of
  $\mathfrak{so}_{n}$}\fg}, \emph{Linear Algebra Appl.} \textbf{374} (2003),
  p.~127--142.

\bibitem{elashvili-1982}
{\scshape A.~Elashvili} -- {\og {Frobenius Lie algebras}\fg},
  \emph{Funktsional. Anal. i Prilozhen.} \textbf{16} (1982), no.~4, p.~94--95
  (Russian).

\bibitem{elashvili-ooms-2003}
{\scshape A.~Elashvili {\normalfont \smfandname} A.~I. Ooms} -- {\og {On
  commutative polarizations}\fg}, \emph{J. Algebra} \textbf{26} (2003), no.~1,
  p.~129--154.

\bibitem{joseph-1977}
{\scshape A.~Joseph} -- {\og {A preparation theorem for the prime spectrum of a
  semisimple Lie algebra}\fg}, \emph{J. Algebra} \textbf{48} (1977), no.~2,
  p.~241--289.

\bibitem{khakimdjanov-goze-medina-2004}
{\scshape Y.~Khakimdjanov, M.Goze {\normalfont \smfandname} A.~Medina} -- {\og
  {Symplectic or contact structures on Lie groups}\fg}, \emph{Differential
  Geom. Appl.} \textbf{21} (2004), no.~1, p.~41--54.

\bibitem{khalgui-torasso-1993}
{\scshape M.~S. Khalgui {\normalfont \smfandname} P.~Torasso} -- {\og {Formule
  de Poisson-Plancherel pour un groupe presque alg\'ebrique r\'eel. I.
  Transform\'ee de Fourier d'int\'egrales orbitales}\fg}, \emph{J. Funct.
  Anal.} \textbf{116} (1993), p.~359--440.

\bibitem{khalgui-torasso-2006}
\bysame , {\og {La formule de Plancherel pour les groupes de Lie presque
  alg\'ebriques r\'eels}\fg}, \emph{J. Funct. Anal.} \textbf{235} (2006),
  p.~449--542.

\bibitem{kimura-b-2003}
{\scshape T.~Kimura} -- \emph{{Introduction to prehomogeneous vector spaces}},
  Translations of Mathematical Monographs, vol. 215, American Mathematical
  Society, 2003, Translated from the 1998 Japanese original by Makoto Nagura
  and Tsuyoshi Niitani and revised by the author.

\bibitem{kirillov-1962}
{\scshape A.~A. Kirillov} -- {\og {Repr\'esentations unitaires des groupes de
  Lie nilpotents}\fg}, \emph{Uspehi Mat. Nauk} \textbf{17} (1962), p.~57--110,
  en Russe.

\bibitem{kosmann-stern-1974}
{\scshape Y.~Kosmann {\normalfont \smfandname} S.~Sternberg} -- {\og
  {Conjugaison des sous-alg\`ebres d'isotropie}\fg}, \emph{C. R. Acad. Sci.
  Paris S\'er. A} \textbf{279} (1974), p.~777--779.

\bibitem{kostant}
{\scshape B.~Kostant} -- {\og {The cascade of orthogonal roots and the center
  of $U(\n)$}\fg}, \emph{Preprint}.

\bibitem{moreau-yakimova-2011}
{\scshape A.~Moreau {\normalfont \smfandname} O.~Yakimova} -- {\og {Coadjoint
  orbits of reductive type of seaweed Lie algebras}\fg}, arXiv:1101.0902.

\bibitem{ooms-1980}
{\scshape A.~I. Ooms} -- {\og {On Frobenius Lie algebras}\fg}, \emph{Comm.
  Algebra} \textbf{8} (1980), no.~1, p.~13--52.

\bibitem{panyushev-2005}
{\scshape D.~Panyushev} -- {\og {An extension of Raïs' theorem and seaweed
  subalgebras of simple Lie algebras}\fg}, \emph{Ann. Inst. Fourier (Grenoble)}
  \textbf{55} (2005), no.~3, p.~693--715.

\bibitem{premet-skryabin-1999}
{\scshape A.~Premet {\normalfont \smfandname} S.~Skryabin} -- {\og
  {Representations of restricted Lie algebras and families of associative $\cur
  L$-algebras}\fg}, \emph{J. Reine Angew. Math.} \textbf{507} (1999),
  p.~189--218.

\bibitem{tauvel-yu-2004-b}
{\scshape P.~Tauvel {\normalfont \smfandname} R.~W.~T. Yu} -- {\og {Indice et
  formes linéaires stables dans les algèbres de Lie}\fg}, \emph{J. of Algebra}
  \textbf{273} (2004), p.~507--516.

\bibitem{tauvel-yu-2004-a}
\bysame , {\og {Sur l'indice de certaines algèbres de Lie}\fg}, \emph{Ann.
  Inst. Fourier} \textbf{54} (2004), no.~6, p.~1793--1810.

\bibitem{tauvel-yu-2005}
{\scshape P.~Tauvel {\normalfont \smfandname} R.~Yu} -- \emph{{Lie algebras and
  algebraic groups}}, Springer Monographs in Mathematics, Springer-Verlag,
  Berlin, 2005.

\end{thebibliography}
\end{document}